\documentclass[a4,10pt,twocolumn]{scrartcl}
\usepackage[top=2.5cm, bottom=2.5cm, left=1.6cm, right=1.6cm]{geometry}
\pdfoutput=1

\input{myStandardPreamble_arXiv.tex}
\usepackage[english]{babel}
\usepackage[utf8x]{inputenc}
\usepackage{enumitem}					
\usepackage[normalem]{ulem}

\usepackage{libertine}
\usepackage[T1]{fontenc}
\usepackage[noBBpl]{mathpazo}			

\usepackage[hang]{footmisc}
\setfnsymbol{wiley}

\makeatletter
\renewcommand{\@biblabel}[1]{[#1]\hfill}
\makeatother

\makeatletter
\let\NAT@parse\undefined
\makeatother

\usepackage[format=plain,			
            labelsep=period,		
            font=small,			
            labelfont=bf,			
            skip=5pt			
            ]{caption}
            
\usepackage[final,kerning=false,protrusion=false]{microtype}

\usepackage{authblk}

\newcommand{\newmarkedtheorem}[1]{%
  \newenvironment{#1}
    {\pushQED{\oprocend}\csname inner@#1\endcsname}
    {\popQED\csname endinner@#1\endcsname}%
  \newtheorem{inner@#1}%
}
\newmarkedtheorem{theorem}{Theorem}
\newmarkedtheorem{lemma}{Lemma}
\newmarkedtheorem{definition}{Definition}
\newmarkedtheorem{remark}{Remark}
\newmarkedtheorem{example}{Example}
\newmarkedtheorem{assumption}{Assumption}
\newmarkedtheorem{corollary}{Corollary}
\newmarkedtheorem{proposition}{Proposition}
\newmarkedtheorem{problemFormulation}{Problem}

\makeatletter
\def\th@plain{%
  \thm@notefont{}
  \normalfont 
}
\def\th@definition{%
  \thm@notefont{}
  \normalfont 
}
\makeatother

\addtokomafont{paragraph}{\itshape}

\colorlet{fictitiousColor}{green!45!black}    

\let\citep\cite
\let\citet\cite

\newcommand{\nat}{\mathbb{N}}								
\newcommand{\natPos}{\mathbb{N}_{>0}}						
\newcommand{\integers}{\mathbb{Z}} 						
\newcommand{\real}{\mathbb{R}}							
\newcommand{\realNonNeg}{\mathbb{R}_{\geq 0}}				
\newcommand{\C}{\mathcal{C}}								
\newcommand{\unitCircle}{\mathbb{T}}						

\newcommand{\perc}{\%}									

\newcommand{\kk}{_{k+1}}
\def\k{_{k}}

\newcommand{\lext}{\ell_e}								
\newcommand{\ltwo}{\ell_2}								
\newcommand{\ltworho}{\ell_{2,\rho}}						
\newcommand{\lf}{\ell_f}									

\newcommand{\fObj}{H}

\newcommand{\classObj}[2]{\mathcal{S}_{#1,#2}}

\newcommand{\ABCD}[4]{ \left(\begin{array}{c|c}  {#1} & {#2} \\[0.15em] \hline \\[-0.95em]  {#3} & {#4} \end{array}\right)}

\newcommand{\optSign}{\star}

\newcommand{\optVar}{z}

\newcommand{\xTrafo}{\xi}

\newcommand{\xOpt}{x}
\newcommand{\Aopt}{A}
\newcommand{\Bopt}{B}
\newcommand{\Copt}{C}
\newcommand{\Dopt}{D}

\newcommand{\gradDesc}{Gradient Descent}
\newcommand{\nesterov}{Nesterov's Method}
\newcommand{\heavyBall}{Heavy Ball Method}
\newcommand{\tripleMomentum}{Triple Momentum Method}

\newcommand{\Anom}{\Aopt_{\textup{nom}}}

\newcommand{\perf}{\textup{p}}
\newcommand{\yperf}{y_\perf}

\newcommand{\wperf}{w_\perf}
\newcommand{\Bperf}{B_\perf}
\newcommand{\Cperf}{C_\perf}
\newcommand{\Dperf}{D_\perf}

\newcommand{\Gwtoy}{G_{yw}}
\newcommand{\Gwperftoy}{G_{yw_\perf}}
\newcommand{\Gwtoyperf}{G_{y_{\perf}w}}
\newcommand{\Gwperftoyperf}{G_{y_{\perf}w_{\perf}}}

\newcommand{\fdi}[1]{\stackrel{\scriptscriptstyle\mathbb{T}}{#1}}

\newcommand{\yAdd}{y_{\textup{in}}}
\newcommand{\yAddTrafo}{\tilde{y}_{\textup{in}}}

\newcommand{\causal}{-}
\newcommand{\anticausal}{+}
\newcommand{\dimCausal}{\ell_\causal}
\newcommand{\dimAnticausal}{\ell_\anticausal}

\newcommand{\tfbasisCausal}{\psi_\causal}
\newcommand{\AbasisCausal}{A_\causal}
\newcommand{\BbasisCausal}{B_\causal}
\newcommand{\CbasisCausal}{C_\causal}
\newcommand{\DbasisCausal}{D_\causal}
\newcommand{\tfbasisAnticausal}{\psi_\anticausal}
\newcommand{\AbasisAnticausal}{A_\anticausal}
\newcommand{\BbasisAnticausal}{B_\anticausal}
\newcommand{\CbasisAnticausal}{C_\anticausal}
\newcommand{\DbasisAnticausal}{D_\anticausal}

\newcommand{\tfMultFac}{\psi_\Delta}
\newcommand{\AMultFac}{A_\Delta}
\newcommand{\BMultFac}{B_\Delta}
\newcommand{\CMultFac}{C_\Delta}
\newcommand{\DMultFac}{D_\Delta}

\newcommand{\complete}{\textup{c}}
\newcommand{\tfcompleteNoPerf}{\psi_\complete}
\newcommand{\AcompleteNoPerf}{A_\complete}
\newcommand{\BcompleteNoPerf}{B_\complete}
\newcommand{\CcompleteNoPerf}{C_\complete}
\newcommand{\DcompleteNoPerf}{D_\complete}

\newcommand{\tfcomplete}{\boldsymbol{\psi}_{\textbf{\complete}}}
\newcommand{\Acomplete}{\boldsymbol{A}_\complete}
\newcommand{\Bcomplete}{\boldsymbol{B}_\complete}
\newcommand{\Ccomplete}{\boldsymbol{C}_\complete}
\newcommand{\Dcomplete}{\boldsymbol{D}_\complete}

\newcommand{\tfzf}{E}

\def\IQC{\textup{IQC}}

\newcommand{\img}{\mathrm{j}} 			
\newcommand{\complConj}{*}				

\newcommand{\nonrep}{\textup{n-rep}}
\newcommand{\rep}{\textup{rep}}
\newcommand{\deltaUnstr}{\mathbf{\Delta}}
\newcommand{\deltaNonrep}{\mathbf{\Delta}_\nonrep}
\newcommand{\deltaRep}{\mathbf{\Delta}_\rep}

\newcommand{\toep}{\textup{Toep}}

\newcommand{\opTimeDomain}{\phi}

\newcommand{\setDoublyHyp}{\mathbf{M}}

\newcommand{\trunc}[1]{O_{#1}}

\newcommand{\id}{\textup{id}}

\newcommand{\trafoSector}{W}

\newcommand{\trafoStructured}{T}

\newcommand{\emphDef}[1]{\textit{#1}}

\newcommand{\sepSet}{\mid}

\newcommand{\expec}{\textup{E}}

\newcommand{\oprocendsymbol}{{$\bullet$}} 					
\def\oprocend {{
\parfillskip=0pt        
\widowpenalty=10000     
\displaywidowpenalty=10000  
\finalhyphendemerits=0  
\leavevmode             
\unskip                 
\nobreak                
\hfil                   
\penalty50              
\hskip.45em              
\null                   
\hfill                  
\oprocendsymbol
\par}}                  

\renewcommand{\sout}[1]{}

\LetLtxMacro\shortArray\array

\AtBeginEnvironment{shortArray}{\setlength{\arraycolsep}{2pt}\def\arraystretch{1}}

\LetLtxMacro\shorterArray\array

\AtBeginEnvironment{shorterArray}{\setlength{\arraycolsep}{3pt}\def\arraystretch{1}}

\setenumerate{leftmargin=*,wide=0pt}

\newboolean{proofsAppendix}
\setboolean{proofsAppendix}{true}

\newboolean{singleColumn}
\setboolean{singleColumn}{false}

\newboolean{longVersion}
\setboolean{longVersion}{true}

\newboolean{saveSpace}
\setboolean{saveSpace}{true}

\title{{\LARGE \textbf
Robust and structure exploiting optimization algorithms:} \\
\LARGE \textbf An integral quadratic constraint approach
}

\setkomafont{section}{\Large}
\setkomafont{subsection}{\large}
\setkomafont{subsubsection}{\itshape}

\begin{document}

\date{}
\author[1]{Simon Michalowsky}
\author[2]{Carsten Scherer}
\author[1]{Christian Ebenbauer}
\affil[1]{Institute for Systems Theory and Automatic Control, University of Stuttgart, Germany \protect\\ \texttt{\small $\lbrace$michalowsky,ce$\rbrace$@ist.uni-stuttgart.de}}
\affil[2]{Mathematical Systems Theory, University of Stuttgart, Germany  \protect\\ \texttt{\small carsten.scherer@mathematik.uni-stuttgart.de} \protect\\[3em]}

\maketitle

\begin{abstract}
\textbf{Abstract.} We consider the problem of analyzing and designing gradient-based
discrete-time optimization algorithms for a class of unconstrained
optimization problems having strongly convex objective functions with Lipschitz continuous
gradient. By formulating the problem as a robustness analysis problem
and making use of a suitable adaptation of the theory of integral quadratic constraints, we establish
a framework that allows to analyze convergence rates and robustness properties of existing algorithms and 
enables the design of novel {robust} optimization algorithms with prespecified guarantees
{capable of exploiting additional structure in the objective function.} 

\end{abstract}


\section{Introduction}
Optimization algorithms are of key importance in science and engineering.
{First order, i.e., gradient-based algorithms, are an important
subclass that have proved themselves in a range
of applications. In recent years, such algorithms}
have regained interest since they are particularly suitable for large-scale
optimization. While {many variants} of gradient-based optimization
algorithms {are} known in {the} literature, no {general} framework for
their analysis and design exists. {Still,} 
a lot of these algorithms fall into the class of Lur'e
systems~\citep{lure1944theory}, i.e., a given or to be 
designed linear system in feedback with a nonlinearity,
given by the gradient in {optimization}. Lur'e systems and
the corresponding absolute stability problem are classical control
problems leading to {celebrated} results such as the Popov or the circle criterion,
which can be seen as {a} pioneering {contribution} to robust control theory.
Quite astonishingly, while not being unknown~\citep{polyak1987introduction}, 
the apparent relation of systems and control theory, in 
particular {absolute stability and} robust control theory, to optimization algorithm
analysis and design has not yet been exploited heavily.
Only recently, several works relying on this systems theoretic
view on optimization algorithms have been published~\citep{bhaya2006control,lessard2016analysis,drori2014firstorder,wilson2016lyapunov,taylor2017pesto,mic2014heavy,mic2016extremum,duerr2012algorithms},
partly also providing different approaches for convergence rate analysis
using {techniques from} robust control.
However, convergence rates are only one side of the coin; 
in several applications, e.g., in a data-based setting, also
robustness with respect to various kinds of disturbances is
a key issue {and similar analysis tools have been developed~\citep{aybat2019robust,mohammadi2018variance}.} In addition, the problem of designing algorithms specifically 
tailored to classes of {structured} optimization problems 
has only been touched upon so far~\citep{drori2018efficient,fazlyab2018design,lessard2019synthesis}. 
{However, in situations where a class of optimization
problems needs to be solved repeatedly online{,} as {for example
in model predictive control or reinforcement learning},
well-performing algorithms are key.} 

The present paper addresses these {two} issues and contributes to the existing literature 
by providing a systematic {framework}
to the {analysis and} design of robust and structure exploiting optimization algorithms. In particular, 
we address the following problem: Given a class of objective
functions, design a gradient-based optimization algorithm
with a guaranteed convergence rate that {also} fulfills
certain $H_2$-performance specifications. We show that these performance specifications can
be related to noise rejection properties of the algorithm such as 
the effect of additive gradient noise on the {variance} {of the algorithm's output}. 
We further address the problem of how to incorporate
possible structural properties of the class of objective functions
in our framework, and how this can be exploited to design novel
optimization algorithms superior to standard ones in terms of convergence rates.
To this end, in the spirit of~\citet{mic2014heavy,lessard2016analysis},
we reformulate this design problem as a robust controller
synthesis problem and employ integral quadratic constraints theory.
By building upon and extending these well-established results, we are able to
provide a {general} framework for algorithm analysis and design.
To validate the practical applicability and relevance, 
we provide several numerical results illustrating our methodology. 

More specifically, our main contributions are as follows:
We propose a class of gradient-based {optimization} algorithms
{\sout{that generalizes existing algorithms}} and derive
necessary and sufficient conditions
for these algorithms to be capable of solving a class of optimization
problems~(\Cref{lemmaConditionsAlgorithm}){, thereby providing the basis for the design of novel algorithm that generalize those known from the literature}. 
Embedding the problem in the framework of robust control, we then derive
convex analysis tools by means of linear matrix inequalities (LMIs), both in 
regard to convergence rates~(\Cref{lemmaAnalysisConvRate}) and
robustness~(\Cref{lemmaAnalysisH2}).
{To this end, we provide a general procedure to obtain
multipliers for exponential stability results from standard ones (\Cref{lemmaMultipliersTransformedLoop})
and utilize this to derive a class of multipliers generalizing
those proposed in~\citet{boczar2015exponential}, \citet{lessard2016analysis}, \citet{freeman2018noncausal}
(\Cref{lemmaZFIQCrho}).} 
We further provide convex synthesis conditions allowing to design novel algorithms with
specified robustness properties~(\Cref{lemmaConservativeConvexSynthesisNew})
and {show} how to additionally exploit structural characteristics of the 
objective function~(\Cref{lemmaZFstructuredTimeDomain}).

\section{Preliminiaries}\label{secPreliminaries}
\subsection{Notation}\label{secNotation}
We let $ \nat $ denote the set of non-negative integers,
$ \natPos $ the set of positive integers and denote by
$ \mathbb{Z} $ the set of all integers. We write $ \C^p $, $ p \in \nat $,
for the set of $p$-times continuously differentiable functions.
For $ n, m \in \natPos $, we denote by $ I_{n} \in \real^{n \times n} $
the $ n \times n $ identity matrix, by $ 0_{n} \in \real^{n \times n} $
the $ n \times n $ matrix of zeros and by $ 0_{n \times m} \in \real^{n \times m} $
the $ n \times m $ matrix of zeros. Sometimes we omit the subscript
if the dimensions are obvious. We let $ \mathbf{1} $ denote the column
vector with all entries being equal to one.
For a square matrix $ A $ we let $ \textup{tr}(A) $ denote the trace of $A$.
For any two matrices $ A_1, A_2 $, we let $ A_1 \otimes A_2 $ denote
{their} Kronecker product.
If $ A_1, A_2 $ are square{, symmetric} and of {the} same dimension,
we write $ A_1 \prec A_2 $ ($ A_1 \preceq A_2 $) if $ A_1 - A_2 $ is negative definite (negative semi-definite) and 
$ A_1 \succ A_2 $ ($ A_1 \succeq A_2 $) if $ A_1 - A_2 $ is positive definite (positive semi-definite). 
In the same manner, we use the relations $ <, \leq, >, \geq $ for elementwise comparison. 
We further
denote by $ \mathbb{T} := \lbrace z \in \mathbb{C} \sepSet \vert z \vert = 1 \rbrace $
the unit circle in the complex plane.

{We let \mbox{$ \lext^p = \lbrace ( q_k )_{k \in \nat} \sepSet q_k \in \real^p \rbrace $} denote the subspace
of all one-sided sequences and denote by $ \lf^p = \lbrace q \in \lext^p \sepSet \exists T \in \nat \text{ s.t. } q_k = 0 \text{ for all } k > T \rbrace $
the set of all finitely supported sequences therein.}
We let $ \ltwo^p \subset \lext^p $ denote all square summable
sequences in $ \lext^p $, i.e.,
\begin{align}
	\ltwo^p = \big\lbrace  q \in \lext^p \sepSet \sum\limits_{k=0}^{\infty} {\Vert q_k \Vert^2} < \infty \big\rbrace.
\end{align}
We sometimes omit the dimension of the signals and write $ \lext \coloneqq \bigcup_{i \in \natPos} \lext^i $, {$ \lf \coloneqq \bigcup_{i \in \natPos} {\lf^i} $},
$ \ltwo \coloneqq \bigcup_{i \in \natPos} \ltwo^i $ for the collection of 
all sequences of the respective type. We denote the standard inner product
by $ \langle u, y \rangle = \sum_{i=-\infty}^\infty u_i^\top y_i $,
where $ u, y \in \ltwo^p $ and the superscript $ {}^\top $ denotes transposition.
We further let $ \Vert y \Vert_{\ltwo} = {\sqrt{\langle y, y \rangle}} $ denote the 
corresponding induced norm.
If it exists{, i.e., if it is well-defined,} we denote by $ \widehat{q} = \mathcal{Z}(q) $ the one-sided
$z$-transform of a signal $ q \in \ltwo^p $ which is defined by
\begin{align}
	\widehat{q}(z) = \mathcal{Z}\big(q\big)(z) = \sum\limits_{k=0}^{\infty} q_k z^{-k},
\end{align}
where $ \widehat{q}: \mathbb{U} \to \mathbb{C}^p $ for some set $ \mathbb{U} \subseteq \mathbb{C} $. 

{We say that an operator $ \phi: \mathcal{U} \to \mathcal{V} $,
$ \mathcal{U}, \mathcal{V} \subseteq \ltwo^p $,
is bounded (on $ \mathcal{U}$) if there exists $ \beta \geq 0 $ such that
$ \Vert \opTimeDomain(y) \Vert_{\ltwo}  \leq \beta \Vert y \Vert_{\ltwo} $
for all $ y \in \mathcal{U} $. {If $ \mathcal{U} = \mathcal{V} = \ltwo^p $, 
we call the infimal $ \beta $ that fulfills this inequality the $ \ltwo $-gain of $ \phi $.}
We further let $ \mathcal{L}(\mathcal{U},\mathcal{V}) $ denote the set of
operators mapping $ \mathcal{U} $ to $ \mathcal{V} $ that are bounded
on $ \mathcal{U} $. }
For two operators $ G_1: \lext^p \to \lext^q $, $ G_2: \lext^q \to \lext^r $,
we let $ G_2 \circ G_1 $ denote the composition, i.e., $ G_2 \circ G_1 $
maps any input $ u \in \lext^p $ to $ y_2 = G_2\big( G_1(u) \big) \in \lext^r $.
For linear operators {$G_1$}, we often omit the brackets and simply write
$ y_1 = G_1 u $ instead of $ y_1 = G_1(u) $. {We further denote by
$ \id: \lext^p \to \lext^p $, $ p \in \natPos $, the identity operator.}

We denote by $ \mathcal{RH}_{\infty}^{n \times m} $ the set of all 
real-rational and proper transfer matrices of dimension $ n \times m $
having all poles in the open unit disk. Similarly, we denote by
$ \mathcal{RL}_{\infty}^{n \times m} $ the set of all real-rational
and proper transfer matrices of dimension $ n \times m $ having
no poles on the unit circle. For a square transfer matrix 
$ G $ we write 
\ifthenelse{\boolean{saveSpace}}{{$ G \stackrel{\scriptscriptstyle\mathbb{S}}{\prec} 0 $}}{
\begin{align}
	G \stackrel{\scriptscriptstyle\mathbb{S}}{\prec} 0
\end{align}}
if the matrix {$ G(z) + G(z)^\complConj $} is negative
definite for all $ z \in \mathbb{S} \subseteq \mathbb{C} $.
For any transfer matrix $ G \in {\mathcal{RL}_{\infty}^{n \times m}} $, we let {$ \opTimeDomain_G $}
denote the corresponding linear {Toeplitz} operator in the time domain {acting on sequences in $ \ltwo^m $}, i.e.,
$ \opTimeDomain_G $ fulfills
\ifthenelse{\boolean{saveSpace}}{{$\mathcal{Z}( \opTimeDomain_G u ) = G \widehat{u}$}}{
\begin{align}
	\mathcal{Z}( \opTimeDomain_G u ) = G \widehat{u}
	\label{eqDefTimeDomainOperator}
\end{align}}
for any $ u \in \ltwo^m $. Note that when $ \opTimeDomain_G $ is interpreted as an infinite-dimensional matrix 
acting on $u$, then $ \opTimeDomain_G $
has a block Toeplitz structure, where the blocks are of size $ n \times m $.
For a proper transfer matrix $ G(z) = C ( zI - A )^{-1} B + D $
we write $ G \sim ( A, B, C, D ) $ {or
\begin{align}
	G \sim \ABCD{A}{B}{C}{D}
\end{align}}
to indicate that $ G $
admits the state-space realization $ (A,B,C,D) $.

For quadratic forms $ x^\complConj P x $, $ x \in \real^n $, $ P \in \real^{n \times n} $, we often shortly write
$ ( \star )^\complConj P x $, where the superscript $ {}^\complConj $
denotes complex conjugation. We additionally replace blocks in symmetric
matrices by the $ \star $ symbol if they can be inferred from symmetry.
We further denote by $ \textup{blkdiag}(M_1,M_2,\dots,M_k) $ the block 
diagonal matrix with blocks $M_1$, $M_2$, $\dots$, $M_k$ on its
diagonal.

\subsection{The class of objective functions}
Throughout the paper we consider the following class of {(strongly)} convex objective functions
with Lipschitz continuous gradient:
{
\begin{definition}[Class $ \classObj{m}{L} $]\label{defClassObjectiveFunctions}
	A function $ \fObj :  \real^p \to \real $ is said to belong
	to the class $ \classObj{m}{L} $, $ L \geq m \geq 0 $, if 
	$ \fObj \in \mathcal{C}^2 $ and 
	\begin{subequations}
	\begin{align}
		\big( \nabla \fObj(x) - \nabla \fObj(y) \big)^\top (x - y) &\geq m \Vert x - y \Vert^2  \label{eqStronglyConvex} \\
		\Vert \nabla \fObj(x) - \nabla \fObj(y) \Vert &\leq L \Vert x-y \Vert \label{eqLipschitzContinuity}
	\end{align}
	\end{subequations}
	for all $ x,y \in \real^p $.
\end{definition}
In view of a lighter notation, we do not explicitly specify the dimension of the domain but tacitly
assume it to be clear from the context.} Note that {if $m>0$, then $ \fObj $ is strongly convex {by~\eqref{eqStronglyConvex}}
and $\fObj$ has a unique minimum.}
The so defined class of functions is well-known in the optimization literature,
see, e.g.,~\citet{nesterov2004introductory} for further properties and details.

\subsection{Robust stability and performance via IQCs}\label{secIntroIQCs}
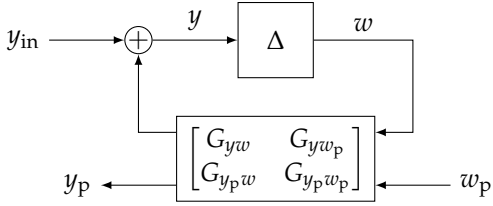
\begin{figure}[t]
	\begin{center}
		\begin{tikzpicture}[>=latex]
			\node[name=nodeUncertainty,draw,rectangle,minimum width=1cm,minimum height=1cm] {$\Delta$};
			\node[name=nodePlant,draw,rectangle,below=0.5cm of nodeUncertainty.south] {$ \begin{bmatrix} \Gwtoy & \Gwperftoy \\ \Gwtoyperf & \Gwperftoyperf \end{bmatrix} $};
			\draw[<-] ([yshift=10pt] nodePlant.east) -- ++(0.5cm,0) |- node[pos=0.75,anchor=south] {$w$} (nodeUncertainty.east);
			\node[name=nodeSum,draw,circle,inner sep=-0.2pt] at ([xshift=-0.5cm] nodePlant.west |- nodeUncertainty.west) {$+$};
			\draw[->] ([yshift=10pt] nodePlant.west) -| (nodeSum.south);
			\draw[->] (nodeSum.east) -- node[anchor=south] {$y$} (nodeUncertainty.west);
			\draw[<-] (nodeSum.west) -- node[pos=1,anchor=east] {$\yAdd$} ++(-1cm,0);
			\draw[<-] ([yshift=-10pt] nodePlant.east) -- node[pos=1,anchor=west] {$w_\perf$} ++(1cm,0);
			\draw[->] ([yshift=-10pt] nodePlant.west) -- node[pos=1,anchor=east] {$y_\perf$} ++(-1cm,0);
		\end{tikzpicture}
		\caption{The considered feedback interconnection.}
		\label{figFeedbackInterconnection}
	\end{center}
\end{figure}
In the following we briefly present some standard results on robust stability
and performance analysis via integral quadratic constraints (IQCs). For a 
more detailed treatment of the subject, we refer the 
reader to the classical paper~\citet{megretski1997iqcs} as well as~\citet{veenman2016robust}
in a continuous-time setup. These results literally carry over to the discrete-time
setup, see, e.g., \citet{kao2012discrete,fetzer2017absolute}.
In the following, we consider feedback interconnections of the form (cf.~\Cref{figFeedbackInterconnection})
\begin{subequations}
\begin{align}
	y &= \Gwtoy w + \Gwperftoy w_\perf + \yAdd \\
	y_\perf &= \Gwtoyperf w + \Gwperftoyperf w_\perf \\
	w &= \Delta(y),
\end{align}\label{eqFeedbackInterconnection}
\end{subequations}
consisting of a linear system with stable transfer matrix $ \Gwtoy \in \mathcal{RH}_{\infty}^{n_y \times n_w} $ in feedback
with an uncertainty \mbox{$ \Delta: \ltwo^{n_y} \to \ltwo^{n_w} $}
and an additional performance channel from $ w_\perf $ to $ y_\perf $,
where 
\mbox{$ \Gwperftoy \in \mathcal{RH}_{\infty}^{n_y \times n_{w_\perf}} $},
\mbox{$ \Gwtoyperf \in \mathcal{RH}_{\infty}^{n_{y_\perf} \times n_w} $} as well as
$ \Gwperftoyperf \in \mathcal{RH}_{\infty}^{n_{y_\perf} \times n_{w_\perf}} $.
Let $ \mathbf{\Delta} \subset \lbrace { {\Delta} \in \mathcal{L}(\ltwo^{n_y}, \ltwo^{n_w}), {\Delta} \textup{ causal} } \rbrace $ denote the set of uncertainties, i.e.,
$ \Delta \in \mathbf{\Delta} $. It is common to impose the following assumption:
\begin{assumption}\label{assStarShapedUncertainty}
If $ \Delta \in \mathbf{\Delta} $ then $ \tau \Delta \in \mathbf{\Delta} $
	for all $ \tau \in [0,1] $.
\end{assumption}
In many cases, this assumption can be ensured by a proper redefinition of
the uncertainty. In particular, this assumption will be met by the class of uncertainties 
we will consider in our problem setup.
Our goal is to show robust stability
of the feedback interconnection without the performance channel
for all $ \Delta \in \mathbf{\Delta} $ and robust performance
{concerning} the performance channel from $ w_\perf $ to $ y_\perf $ for all $ \Delta \in \mathbf{\Delta} $
with respect to some (integral) quadratic performance criterion.
In the remainder of this section, we briefly repeat the basic definitions 
and results from robust control and IQC theory.
We first give a proper definition of robust stability.
\begin{definition}[Robust stability]
The feedback interconnection~\eqref{eqFeedbackInterconnection} is said to be \emphDef{robustly stable against
	$ \mathbf{\Delta} $} if, for $ w_\perf = 0 $, it is well-posed
	and the $ \ltwo $-gain of the map from $ \yAdd $ to $ y $ is 
	bounded for all $ \Delta \in \mathbf{\Delta} $.
\end{definition}
Therein, well-posedness of the feedback interconnection means that 
the map $ q \mapsto (I - \Gwtoy \Delta) q $ has a causal inverse for all
$ \Delta \in \mathbf{\Delta} $. 
We next give the definition of an integral quadratic constraint (IQC).
\begin{definition}[IQC]
Let {a so-called multiplier $ \Pi:\mathcal{RL}_{\infty}^{(p_1+p_2)\times(p_1+p_2)} $, $ p_1, p_2 \in \natPos $,
	be given.}
	We say that \emphDef{two signals $ q_1 \in \ltwo^{p_1}, q_2 \in \ltwo^{p_2} $
	satisfy the IQC defined by $ \Pi $} if 
	{
	\ifthenelse{\boolean{singleColumn}}{
	\begin{align}
		\IQC(\Pi,q_1,q_2) &= \int\limits_{0}^{2\pi} \begin{bmatrix} \widehat{q_1}(e^{\img\omega}) \\ \widehat{q_2}(e^{\img\omega}) \end{bmatrix}^* \Pi(e^{\img\omega}) \begin{bmatrix} \widehat{q_1}(e^{\img\omega}) \\ \widehat{q_2}(e^{\img\omega})  \end{bmatrix} \mathrm{d}\omega \geq 0.
		\label{eqDefIQC}
	\end{align}
	}{
	\begin{flalign}
		\IQC(\Pi,q_1,q_2) &= \int\limits_{0}^{2\pi} \begin{bmatrix} \widehat{q_1}(e^{\img\omega}) \\ \widehat{q_2}(e^{\img\omega}) \end{bmatrix}^* \Pi(e^{\img\omega}) \begin{bmatrix} \widehat{q_1}(e^{\img\omega}) \\ \widehat{q_2}(e^{\img\omega})  \end{bmatrix} \mathrm{d}\omega \geq 0. \hspace*{-1em} &
		\label{eqDefIQC}
	\end{flalign}}}
	We further say that \emphDef{an operator $ \Delta: \ltwo^{p_1} \to \ltwo^{p_2} $ 
	satisfies the IQC defined by $ \Pi $} if $ \IQC\big( \Pi, q_1, \Delta(q_1) \big) \geq 0 $ holds
	for all $ q_1 \in \ltwo^{p_1} $.
\end{definition}
\begin{remark}\label{remarkIQCTimeDomain}
By Parseval's Theorem, {the IQC}~\eqref{eqDefIQC} can equivalently be formulated in the time domain as
	\begin{align}
		\left\langle \begin{bmatrix} q_1 \\ q_2 \end{bmatrix}, \opTimeDomain_\Pi \begin{bmatrix} q_1 \\ q_2 \end{bmatrix} \right\rangle \geq 0,
	\end{align}
	where $ \opTimeDomain_\Pi $ is the linear operator
	defined by $ \Pi $ as introduced in\ifthenelse{\boolean{saveSpace}}{~\Cref{secNotation}.}{~\eqref{eqDefTimeDomainOperator}.}
\end{remark}
In general, IQCs can be used to i) characterize the uncertain operator $ \Delta $ as well as 
to ii) describe performance criteria. For the case i), 
we try to find a set of suitable multipliers $ \mathbf{\Pi} $
such that each $ \Delta \in \mathbf{\Delta} $ satisfies the IQC defined by $ \Pi $ for all $ \Pi \in \mathbf{\Pi} $.
Similarly, for ii), we let $ q_1 = w_\perf $, $ q_2 = y_\perf $ and specify
a set of multipliers $ \mathbf{\Pi}_p \subset \mathcal{RL}_{\infty}^{(n_{w_\perf}+n_{y_\perf}) \times (n_{w_\perf}+n_{y_\perf})} $ that characterize the desired performance
criterion imposed on the performance channel. The corresponding performance IQC
is then given by 
\begin{align}
	\IQC\big( -\Pi_\perf , w_\perf, y_\perf \big) \geq 0,  \label{eqPerformanceIQC}
\end{align}
where $ \Pi_\perf \in \mathbf{\Pi}_\perf $.
We then define robust performance of~\eqref{eqFeedbackInterconnection} as follows.
\begin{definition}[Robust performance]
We say that the feedback interconnection~\eqref{eqFeedbackInterconnection}
	\textit{achieves robust performance with respect to $ \Pi_\perf $ against $ \mathbf{\Delta} $} if the 
	feedback interconnection is robustly stable against $ \mathbf{\Delta} $
	and, for $ \yAdd = 0 $,~\eqref{eqPerformanceIQC} holds for $ \Pi_\perf \in \mathbf{\Pi}_\perf $
	{and all $ \Delta \in \mathbf{\Delta} $.}
\end{definition}
In the following we limit ourselves to performance multipliers where the
block corresponding to the quadratic terms in $ y_\perf $ is positive semi-definite,
i.e., we assume
\ifthenelse{\boolean{singleColumn}}{
\begin{flalign}
	\mathbf{\Pi}_\perf \subset \bigg\lbrace \begin{bmatrix} \Pi_{\perf,11} & \Pi_{\perf,12} \\ \Pi_{\perf,12}^\complConj & \Pi_{\perf,22} \end{bmatrix} \in \mathcal{RL}_{\infty}^{(n_{w_\perf}+n_{y_\perf}) \times (n_{w_\perf}+n_{y_\perf})} \sepSet \Pi_{\perf,22} \fdi{\succeq} 0  \bigg\rbrace.  \label{eqAssSetPerformanceMultipliers}
\end{flalign}}{
\begin{flalign}
	\mathbf{\Pi}_\perf \subset \bigg\lbrace \begin{bmatrix} \Pi_{\perf,11} & \Pi_{\perf,12} \\ \Pi_{\perf,12}^\complConj & \Pi_{\perf,22} \end{bmatrix} \in \mathcal{RL}_{\infty}^{(n_{w_\perf}+n_{y_\perf}) \times (n_{w_\perf}+n_{y_\perf})} \sepSet  \hspace*{-2em} &  &  \nonumber  \\ \Pi_{\perf,22} \fdi{\succeq} 0  \bigg\rbrace.  & & \label{eqAssSetPerformanceMultipliers}
\end{flalign}}
This assumption is met by the most relevant performance criteria and, in particular,
it holds for all performance criteria we consider in the present manuscript. 
We next state the well-known IQC stability theorem.
\begin{theorem}\label{lemmaIQC}
Consider the feedback interconnection~\eqref{eqFeedbackInterconnection}.
	Let some set of multipliers $ \mathbf{\Pi} \subset \mathcal{RL}_{\infty}^{(n_y+n_w) \times (n_y + n_w)} $
	as well as some set of performance multipliers $ \mathbf{\Pi}_\perf $ 
	{with}~\eqref{eqAssSetPerformanceMultipliers} be given.
	Assume that $ \Gwtoy $ has all its poles in the open unit disk
	and assume {$\mathbf{\Delta} $ 
	fulfills~\Cref{assStarShapedUncertainty}}. Suppose further that
	\begin{enumerate}[leftmargin=*]
		\item the interconnection~\eqref{eqFeedbackInterconnection} is well-posed for all $ \Delta \in \mathbf{\Delta} $;
		\item each $ \Delta \in \mathbf{\Delta} $ satisfies the IQC defined by $ \Pi $ for all $ \Pi \in \mathbf{\Pi} $.
	\end{enumerate}
	If there exist $ \Pi \in \mathbf{\Pi} $ and $ \Pi_\perf \in \mathbf{\Pi}_\perf $
	such that
	\begin{align}
		\left[
		\vphantom{
		\begin{array}{cc}
			\Gwtoy & \Gwperftoy \\
			I & 0 \\ \hline
			0 & I \\
			\Gwtoyperf & \Gwperftoyperf 
		\end{array}}
		\star
		\right]^\complConj
		\left[
		\begin{array}{c|c}
			\Pi & 0 \\ \hline 0 & \Pi_\perf
		\end{array}\right]
		\left[
		\begin{array}{cc}
			\Gwtoy & \Gwperftoy \\
			I & 0 \\ \hline
			0 & I \\
			\Gwtoyperf & \Gwperftoyperf 
		\end{array}\right]
		\fdi{\prec} 0,
		\label{eqFDIRobustPerformance}
	\end{align}
	then the interconnection~\eqref{eqFeedbackInterconnection} is robustly 
	stable against $ \mathbf{\Delta} $ and it achieves robust performance w.r.t. $ \Pi_p $ against $ \mathbf{\Delta} $.
\end{theorem}
\ifthenelse{\boolean{proofsAppendix}}{
	The proof is literally the same as in the continuous-time setup 
	and is \ifthenelse{\boolean{longVersion}}{included in~\Cref{secProofLemmaIQC} for the sake of completeness.}{{hence omitted here.}}
}{
\begin{proof}
	We follow the lines of the proof of~\citet[Corollary 3]{veenman2016robust}
	and first show robust stability.
	Let $ \Pi_\perf $ be partioned according to~\eqref{eqAssSetPerformanceMultipliers}
	and observe that~\eqref{eqFDIRobustPerformance} implies that
	\begin{align}
		\begin{bmatrix} \Gwtoy \\ I \end{bmatrix}^\complConj 
		\Pi
		\begin{bmatrix} \Gwtoy \\ I \end{bmatrix}
		+
		\Gwtoyperf^\complConj \Pi_{\perf,22} \Gwtoyperf \stackrel{\scriptscriptstyle\unitCircle}{\prec} 0.
	\end{align}
	Since $ \Pi_{\perf,22} \succeq 0 $ according to~\eqref{eqAssSetPerformanceMultipliers},
	we have
	\begin{align}
		\begin{bmatrix} \Gwtoy \\ I \end{bmatrix}^\complConj 
		\Pi
		\begin{bmatrix} \Gwtoy \\ I \end{bmatrix}
		\stackrel{\scriptscriptstyle\unitCircle}{\prec} 0.
	\end{align}
	Thus, with~\Cref{assStarShapedUncertainty} and 1., 2., robust stability 
	follows from the standard IQC-Theorem, see, e.g.,~\citet{kao2012discrete,fetzer2017absolute}.
	For robust performance, let
	Let $ w_p \in \ell_2^{n_{w_\perf}} $
	and let $ w $ denote the corresponding signal resulting from the
	feedback interconnection~\eqref{eqFeedbackInterconnection}.
	Observe that $ w \in \ell_2^{n_w} $ due to robust 
	stability which implies that $ y \in \ell_2^{n_y} $. Hence,
	the $z$-transforms of $ w $ and $ w_\perf $ exist almost everywhere on the unit circle 
	and multiplying~\eqref{eqFDIRobustPerformance} from left by 
	{$ \begin{bmatrix} \widehat{w}(z)^\complConj & \widehat{w}_\perf(z)^\complConj \end{bmatrix} $}
	and from right by its transposed we obtain
	\ifthenelse{\boolean{singleColumn}}{
	\begin{align}
		\begin{bmatrix} \vphantom{\begin{bmatrix} \widehat{y} \\ \widehat{\Delta(y)} \end{bmatrix}} \star \end{bmatrix}^{{\complConj}}
		\Pi_2
		\begin{bmatrix} \widehat{y} \\ \widehat{\Delta(y)} \end{bmatrix}
		+
		\begin{bmatrix} \vphantom{\begin{bmatrix}\widehat{w_\perf} \\ \widehat{y_\perf} \end{bmatrix}} \star \end{bmatrix}^{{\complConj}}
		\Pi_p
		\begin{bmatrix} \widehat{w_\perf} \\ \widehat{y_\perf} \end{bmatrix}
		{{\stackrel{\scriptscriptstyle\unitCircle}{\preceq} 0}}.
	\end{align}	}{
	\begin{flalign}
		\begin{bmatrix} \vphantom{\begin{bmatrix} \widehat{y} \\ \widehat{\Delta(y)} \end{bmatrix}} \star \end{bmatrix}^{{\complConj}}
		\Pi_2
		\begin{bmatrix} \widehat{y} \\ \widehat{\Delta(y)} \end{bmatrix}
		+
		\begin{bmatrix} \vphantom{\begin{bmatrix}\widehat{w_\perf} \\ \widehat{y_\perf} \end{bmatrix}} \star \end{bmatrix}^{{\complConj}}
		\Pi_p
		\begin{bmatrix} \widehat{w_\perf} \\ \widehat{y_\perf} \end{bmatrix}
		\stackrel{\scriptscriptstyle\unitCircle}{{\preceq}} 0.
	\end{flalign}}
	Integrating on both sides over the set $ \unitCircle $ yields
	\ifthenelse{\boolean{singleColumn}}{
	\begin{align}
		\IQC\big(\Pi_2,y,\Delta(y)\big) + \IQC(\Pi_p,w_\perf,y_\perf) {\leq 0}.
	\end{align}}{
	\begin{align}
		\IQC\big(\Pi_2,y,\Delta(y)\big) + \IQC(\Pi_p,w_\perf,y_\perf) {\leq 0}.
	\end{align}}
	Since $ \IQC_{1}\big(\Pi_2,y,\Delta(y)\big) \geq 0 $ by assumption,
	we conclude that the performance criterion defined by $ \Pi_\perf $ is fulfilled.	

\end{proof}
}
The premise in applying~\Cref{lemmaIQC} is then to find a
class of multipliers valid for the class of {considered} uncertainties
as well as an IQC formulation of the desired performance
criterion such that all assumptions in the latter Theorem are met.

\section{Problem formulation}\label{secProblemFormulation}
Consider the following unconstrained convex optimization problem
\begin{align}
	\min\limits_{\optVar \in \real^p} \quad \fObj(\optVar), \label{eqOptProb}
\end{align}
where $ \fObj : \real^p \to \real $, $ \fObj \in \mathcal{C}^2 $,
$ \fObj $ is strongly convex with convexity modulus $m > 0$ and 
has Lipschitz continuous gradient with parameter {$L \geq  m $}, i.e.,
$ H $ is in the class {$ \classObj{m}{L} $} as in~\Cref{defClassObjectiveFunctions}.
Observe that, under these assumptions, for any fixed $\fObj$ in the class,~\eqref{eqOptProb}
has a unique global minimizer which will {subsequently} be denoted by $\optVar^\optSign$.

{Our goal is to analyze and design gradient-based optimization algorithms
that converge to that minimizer $ \optVar^\optSign $. 
In particular, we consider optimization algorithms of the form 
\begin{subequations}
\begin{align}
	\xOpt\kk  &= \Aopt \xOpt\k + \Bopt \nabla H ( \Copt \xOpt\k ) \\
	\optVar\k &= \Dopt \xOpt\k,
\end{align} \label{eqOptAlgo}
\end{subequations}
where $ \xOpt\k = \begin{bmatrix} \xOpt_{k,1}^\top & \dots & \xOpt_{k,n}^\top \end{bmatrix}^\top \in \real^{np} $,
$ \xOpt_{k,i} \in \real^p $, $ i~\in~\lbrace 1,\dots,n\rbrace$, $ \optVar\k \in \real^p $.
More formally, the design problem we want to address is
then as follows: 
Given the objective function parameters {$L \geq m > 0$}, the dimensions $n,p\geq1$,
and a convergence rate bound $ \rho\in (0,1) $, find matrices $ \Aopt \in \real^{np \times np} $, 
$ \Bopt \in \real^{np \times p} $, $ \Copt \in \real^{p \times np} $, $ \Dopt \in \real^{p \times np} $,
such that, for all $ \fObj \in {\classObj{m}{L}} $,~\eqref{eqOptAlgo} has a unique globally asymptotically stable equilbrium at $ \xOpt^\optSign $
with $ \Dopt \xOpt^\optSign = \optVar^\optSign $
and there exists $ \eta > 0 $ such that 
\begin{align}
	\Vert \optVar\k - \optVar^\optSign \Vert^2 \leq \eta \rho^{2k} \Vert \optVar_{0} - \optVar^\optSign \Vert^2
	\label{eqProblemFormulationConvergenceRate}
\end{align}
for each $ \xOpt_0 \in \real^{np} $ {with} $ \optVar_0 \coloneqq \Dopt \xOpt_0 $, i.e.,
$ \optVar_k $ converges exponentially with rate $ \rho $ to $ \optVar^\optSign $.
Note that~\eqref{eqOptAlgo} captures several popular optimization algorithms
with constant step size such as \gradDesc, the \heavyBall, or \nesterov, see also
\eqref{eqDefExistingAlgorithms} and \Cref{tableParametersAlgorithms}.
By solving the design problem, we also get a solution to the corresponding
analysis problem for these algorithms as a by-product{,} in which we aim
to find an as tight as possible convergence rate bound $ \rho $.
We emphasize that~\eqref{eqOptAlgo} is an overparametrization of
the class of algorithms; in fact, without loss of generality,
we can fix $ \Copt = \begin{bmatrix} \Copt_1 & \Copt_2 & \dots & \Copt_n \end{bmatrix} $
arbitrarily as long as $ \Copt_i $ is non-singular for some $i \in \lbrace 1,2,\dots,n \rbrace $.
In particular, we often let $ \Copt = \begin{bmatrix} I & 0 & \dots & 0 \end{bmatrix} $.
{A similar optimization algorithm design problem has been addressed 
in~\citet{fazlyab2018design,lessard2016analysis} with 
structured and parametrized matrices $ \Aopt, \Bopt, \Copt, \Dopt $ 
in~\eqref{eqOptAlgo} and $ n=2 $.
To the best of our knowledge, the only work considering 
general algorithms of the form~\eqref{eqOptAlgo} is~\citet{lessard2019synthesis}
which focuses on determining lower convergence rates bounds.}

While convergence rate bounds are an important performance measure,
also other performance specifications are key in efficiently solving
optimization problems, e.g., how well an algorithm performs in the 
presence of noise. In this {vein}, we consider~\eqref{eqOptAlgo}
together with an additional performance channel and address the
extended design problem aiming not only for minimizing the convergence
rate but also the bounds on the additional performance channel.
We formalize this problem in~\Cref{secFormulationRobustControl} and 
{propose an adapted $H_2$-performance measure. It has already
been recognized in different settings~\citep{mohammadi2018variance,aybat2019robust}
that $H_2$-performance is a measure for the robustness properties
of optimization algorithms. }

We are further interested in designing tailored algorithms for 
certain subclasses of the class of objective functions~{$ \classObj{m}{L} $} 
with additional structural properties such as {having a} diagonal Hessian.
Summing up, in rough words our goal is to design gradient-based algorithms of
the form~\eqref{eqOptAlgo} that are fast, robust and possibly 
exploit additional structural properties of the objective function. }

\section{Reformulation in the robust control framework}\label{secFormulationRobustControl}
Our approach relies on reformulating the problem as a robust control
problem {by} interpreting the gradient of the objective function $ \fObj $
as the uncertainty. {The reformulation as well as the adaptation to
the specific problem at hand is the main subject of this section.}

\subsection{{Equilibrium conditions}}
Before we proceed, we first discuss which properties
an algorithm of the form~\eqref{eqOptAlgo} must fulfill to be, in principle, a candidate
for the previously posed problems. This {important} question, which has not been addressed in the literature so far, is 
answered in the following Theorem {which provides necessary and
sufficient conditions for an algorithm of the form~\eqref{eqOptAlgo} to
possess an equilibrium at the global minimizer of $ \fObj $ for all {$ \fObj \in \classObj{m}{L} $}.} 
\begin{theorem}\label{lemmaConditionsAlgorithm}
Let an algorithm in the form~\eqref{eqOptAlgo} described by $ \Aopt, \Bopt, \Copt, \Dopt $ be given and assume that {the pair $ (\Aopt,\Copt) $ or the pair $ (\Aopt,\Dopt) $ is observable.}
	Then,~\eqref{eqOptAlgo} has
	an equilibrium $ \xOpt^\optSign $ with the property $ \Dopt \xOpt^\optSign = {\Copt \xOpt^\optSign} = \optVar^\optSign $
	for any $ \fObj \in {\classObj{m}{L}} $, {$ L \geq m > 0 $} arbitrary but fixed,
	if and only if there exists $ \Dopt^\dagger \in \real^{np \times p} $ such that
	\begin{subequations}
	\begin{align}
 		\Dopt \Dopt^\dagger &= I_p \label{eqConditionABCD1} \\
		\Copt \Dopt^\dagger &= I_p \label{eqConditionABCD2} \\
		( \Aopt - I ) \Dopt^\dagger &= 0_{np \times p}. \label{eqConditionABCD3}
	\end{align} \label{eqConditionABCD}
	\end{subequations}
	holds. 
\end{theorem}
\ifthenelse{\boolean{proofsAppendix}}{
	A proof {of~\Cref{lemmaConditionsAlgorithm}} is given in~\Cref{secProofConditionsAlgorithm}.
}{
\begin{proof}
	We first show sufficiency and then necessity. \\
\textit{A) Sufficiency.} Let some $ \Dopt^\dagger $ be given that fulfills~\eqref{eqConditionABCD} and let $ \xOpt^\optSign = \Dopt^\dagger \optVar^\optSign $.
Then $ \xOpt^\optSign $ fulfills $ \Dopt \xOpt^\optSign = {\Copt \xOpt^\optSign} = \optVar^\optSign $ by~\eqref{eqConditionABCD1}, \eqref{eqConditionABCD2}. The condition that~\eqref{eqOptAlgo} 
has an equilibrium at $ \xOpt^\optSign = \Dopt^\dagger \optVar^\optSign $ then amounts to
\begin{align}
	\Dopt^\dagger \optVar^\optSign = \Aopt \Dopt^\dagger \optVar^\optSign + \Bopt \nabla \fObj ( \Copt \Dopt^\dagger \optVar^\optSign ).
	\label{eqProofConditionEquilibrium}
\end{align}
Using~\eqref{eqConditionABCD2}, \eqref{eqConditionABCD3}, this
holds if
\ifthenelse{\boolean{saveSpace}}{{$ \Bopt \nabla \fObj(\optVar^\optSign) = 0 $,}}{
\begin{align}
	0 = \Bopt \nabla \fObj(\optVar^\optSign),
\end{align}}
which is fulfilled 
since $ \optVar^\optSign $ is the minimizer of $ \fObj $. \\
{\textit{B) Necessity.} Suppose that~\eqref{eqOptAlgo} has an equilbrium at $ \xOpt^\optSign $
with the property $ \Dopt \xOpt^\optSign = \Copt \xOpt^\optSign = \optVar^\optSign $
for any $ \fObj \in {\classObj{m}{L}}  $. This implies that 
\begin{align}
	\xOpt^\optSign 
	&= \Aopt \xOpt^\optSign + \Bopt \nabla\fObj( \Copt \xOpt^\optSign ) 
	= \Aopt \xOpt^\optSign \label{eqProofConditionEquilibrium1}
\end{align}
since $ \Copt \xOpt^\optSign = \optVar^\optSign $ is the unique minimizer of 
$ \fObj $. Hence,
\begin{align}
	\optVar^\optSign 
	= \Copt \Aopt^i \xOpt^\optSign = \Dopt \Aopt^i \xOpt^\optSign
\end{align}
for any $ i = 0,1,\dots $, and we infer that
\begin{align}
	\begin{bmatrix} I \\ I \\ I \\ \vdots \\ I \end{bmatrix} \optVar^\optSign
	&=
	\underbrace{\begin{bmatrix} \Copt \\ \Copt \Aopt \\ \Copt \Aopt^2 \\ \vdots \\ \Copt \Aopt^{{np-1}} \end{bmatrix}}_{=: Q_\Copt} \xOpt^\optSign 
	=
	\underbrace{\begin{bmatrix} \Dopt \\ \Dopt \Aopt \\ \Dopt \Aopt^2 \\ \vdots \\ \Dopt \Aopt^{{np-1}} \end{bmatrix}}_{=: Q_\Dopt} \xOpt^\optSign.
\end{align}
Since the pair $ (\Aopt,\Copt) $ or the pair $ ( \Aopt,\Dopt ) $ is observable by assumption, the matrix
$ Q_\Copt {\in \real^{n{p^2} \times p}} $ or $ Q_\Dopt {\in \real^{n{p^2} \times p}} $ 
{has full rank} and we infer that there exists some matrix $ U {\in \real^{np \times p}} $ 
independent of $ \xOpt^\optSign, \optVar^\optSign $ such that $ \xOpt^\optSign = U \optVar^\optSign $.
We next show that $ \Dopt^\dagger = U $ fulfills~\eqref{eqConditionABCD}.
To this end, note that $ \optVar^\optSign = \Copt \xOpt^\optSign = \Copt U \optVar^\optSign $
and{, equally well,} $ \optVar^\optSign = \Dopt \xOpt^\optSign = \Dopt U \optVar^\optSign $.
Since these two equation need to hold for all $ \optVar^\optSign $, this 
implies that $ \Copt U = \Dopt U = I $ and we infer that~\eqref{eqConditionABCD1}, \eqref{eqConditionABCD2} hold
for $ \Dopt^\dagger = U $. The third condition~\eqref{eqConditionABCD3}
follows directly from~\eqref{eqProofConditionEquilibrium1} using
$ \xOpt^\optSign = U \optVar^\optSign $.}

\end{proof}
}
{We emphasize that the conditions~\eqref{eqConditionABCD} are
necessary and sufficient, hence assuming that these conditions hold is
no restriction of the class of algorithms.}
Note that~\eqref{eqConditionABCD3} implies that $ \Aopt $ must have at least
$p$ eigenvalues at one. The proof also reveals that $ \xOpt^\optSign = \Dopt^\dagger \optVar^\optSign $
is the unique equilibrium of~\eqref{eqOptAlgo} with the desired property.
We further note that, under the assumption that $ \Dopt $ has full rank, the convergence
rate bound~\eqref{eqProblemFormulationConvergenceRate} holds if and only
if the full state $ \xOpt\k $ itself converges exponentially with rate $ \rho $
to $ \Dopt^\dagger \optVar^\optSign $. In view of this,
the choice of $ \Dopt $ does not alter the achievable convergence rate. Hence, remembering the previous
discussion that~\eqref{eqOptAlgo} is an overparametrization, 
we often let $ \Dopt = \Copt = \begin{bmatrix} I & 0 & \dots & 0 \end{bmatrix} $, $ \Dopt^\dagger = \Dopt^\top $,
such that~\eqref{eqConditionABCD1},~\eqref{eqConditionABCD2} are ensured to hold.
In the remainder of the paper we limit ourselves to algorithms
described by $ \Aopt, \Bopt, \Copt, \Dopt $ that fulfill all
conditions in~\Cref{lemmaConditionsAlgorithm}.

\subsection{Problem reformulation}
We next {embed} the design problem from~\Cref{secProblemFormulation} in
the standard setup as introduced in~\Cref{secIntroIQCs}.
To this end, consider~\eqref{eqOptAlgo} {initialized at any $ \xOpt_0 $} together with the state transformation
$ \xTrafo\k = \xOpt\k - \Dopt^\dagger \optVar^\optSign $, 
where $ \Dopt^\dagger \in \real^{np \times p} $ is {taken} such that~\eqref{eqConditionABCD} hold.
The transformed algorithm together with an additional performance
channel from $ \wperf $ to $ \yperf $ as described in~\Cref{secProblemFormulation} then takes the form
\begin{subequations}
\begin{align}
	\xTrafo\kk &= \Aopt \xTrafo\k + \Bopt \nabla \fObj ( \Copt \xTrafo\k + \optVar^\optSign ) + \Bperf \wperf{_,}{\k} \\
	\optVar\k &= \Dopt \xTrafo\k + \optVar^\optSign \label{eqTransformedOptAlgoOutput} \\
	\yperf{_,}{\k} &= \Cperf \xTrafo\k + \Dperf \wperf{_,}{\k},
\end{align}\label{eqTransformedOptAlgo}
\end{subequations}
{with initial condition $ \xTrafo_0 = \xOpt_0 - \Dopt^\dagger \optVar^\optSign $ and} 
where $ \Bperf \in \real^{np \times n_{\wperf}} $, 
$ \Cperf \in \real^{n_{\yperf} \times np } $ and
$ \Dperf \in \real^{n_{\yperf} \times n_{\yperf}} $,
$ n_{\wperf}, n_{\yperf} \in \natPos $.
Note that if $ \xTrafo $ converges exponentially to zero with rate $ \rho $,
then so does $ x $ converge to $ \Dopt^\dagger \optVar^\optSign $ in~\eqref{eqOptAlgo}; in other words~\eqref{eqProblemFormulationConvergenceRate}
holds. 

As we are aiming for analyzing and designing {$\Aopt, \Bopt, \Copt, \Dopt$ in}~\eqref{eqTransformedOptAlgo}
for a whole class of objective functions, we interpret $ \nabla H $
as an uncertainty and define the {causal and bounded} uncertain operator
$ \Delta_\fObj: {\lext^p} \to {\lext^p} $ as 
\begin{align}
	\Delta_\fObj(y)\k := \nabla \fObj ( y_k + \optVar^\optSign ) - m y_k
	\label{eqDefDelta}
\end{align}
for any $ y = {[ y_0, y_1, y_2, \dots ]} \in {\lext^p} $.
The corresponding set of all admissible operators is then defined as
\begin{align}
	\deltaUnstr(m,L) := \big\lbrace \Delta_\fObj : \fObj \in {\classObj{m}{L}} \big\rbrace.
	\label{eqDefSetOfDelta}
\end{align}
Note that the set $ {\deltaUnstr(m,L)} $ fulfills~\Cref{assStarShapedUncertainty}.
Observe further that each 
$ \Delta \in {\deltaUnstr(m,L)} $ is a slope-restricted operator with slope between
$ 0 $ and $ L-m $.
In the spirit of~\Cref{secIntroIQCs}
and \eqref{eqFeedbackInterconnection} we {drop the output~\eqref{eqTransformedOptAlgoOutput},} write~\eqref{eqTransformedOptAlgo} as
\begin{subequations}
\begin{align}
	\xTrafo\kk  &= ( \Aopt + m \Bopt \Copt ) \xTrafo\k + \Bopt w\k + \Bperf \wperf{_,}{\k} \\
	y\k &= \Copt \xTrafo\k \\
	\yperf{_,}{\k} &= \Cperf \xTrafo\k + \Dperf \wperf{_,}{\k} \\
	w\k &= {  \Delta_\fObj( y ) \k }
\end{align}\label{eqRobustStandardForm}
\end{subequations}
and {choose} {the transfer functions in~\eqref{eqFeedbackInterconnection} as}
{
\begin{align}
	\begin{bmatrix} \Gwtoy & \Gwperftoy \\ \Gwperftoyperf & \Gwtoyperf \end{bmatrix}
	\sim \ABCD{\Aopt + m \Bopt \Copt}{\begin{bmatrix} \Bopt & \Bperf \end{bmatrix}}{\begin{bmatrix} \Copt \\ \Cperf \end{bmatrix}}{\begin{bmatrix} 0 & 0 \\ 0 & \Dperf \end{bmatrix}}.
	\label{eqStateSpaceRepresentations}
\end{align}}
We note that the latter feedback interconnection is well-posed since $ \Gwtoy $ 
is strictly proper and $ \Delta $ is a static uncertainty.
{The original design problem can then be formulated as follows:
Given {a} set of operators~$\deltaUnstr(m,L)$ {(or a structured subclass thereof, see~\Cref{secStructureExploitingAlgorithms})} and a performance
channel described by $ \Bperf, \Cperf, \Dperf $, {our goal is to}
design $ \Aopt, \Bopt, \Copt $ in such a way that (a) 
the transformed algorithm dynamics~\eqref{eqRobustStandardForm}
have a globally asymptotically stable equilibrium at $ \xTrafo^\optSign = 0 $
for all $ \Delta_\fObj \in \deltaUnstr(m,L) $ and $ \xTrafo $
converges exponentially to $ 0 $ with rate $ \rho $ and (b)
the performance channel defined by {the map from $ \wperf $ to $ \yperf$ in~\eqref{eqTransformedOptAlgo}} fulfills 
a specified performance bound. {Note that we do no longer 
consider $\Dopt$ as a design variable since it does not have an effect
on the goals (a), (b).} }

Our approach to address the
problem at hand is then to make use of IQC theory, in particular~\Cref{lemmaIQC}.
{The key steps~\ref{itemStep1}--\ref{itemStep3} that pave the way to an implementable solution of the design problem
are then as follows:
\begin{enumerate}[leftmargin=*,label={$($S\arabic*$)$}]
	\item Extend~\Cref{lemmaIQC} to allow for exponential stability results (\Cref{subsecExponentialConvergenceViaIQCs}{, in particular~\Cref{lemmaIQCexponentialStability}}).\label{itemStep1}
	\item Derive IQCs valid for the class of uncertainties that are suitable for the {exponential} stability result (\Cref{subsecIQCsForTheClassOfUncertainties}{, in particular~\Cref{lemmaZFIQCrho}}).\label{itemStep2}
	\item Determine state-space representations of all transfer functions occurring in the frequency domain inequality (\Cref{subsecMultiplierParametrization}).\label{itemStep3}
\end{enumerate}
We note that the last step is more or less standard; still,
we include it here in view of our goal of presenting an easily implementable
framework.
}

\subsection{Exponential convergence via IQCs}\label{subsecExponentialConvergenceViaIQCs}
{We begin with the first step~\ref{itemStep1} of extending
\Cref{lemmaIQC} which} only allows to conclude $\ltwo$-stability results for~\eqref{eqRobustStandardForm}.
While it is noted in~\citet{megretski1997iqcs} that $\ltwo$-stability implies
exponential stability for some classes of systems, the resulting 
rate bounds are often conservative~\citep{boczar2017exponential}. 
In~\citet{boczar2015exponential,lessard2016analysis}, the concept of $\rho$-IQCs is introduced to derive exponential stability results.
In the present manuscript, we focus on obtaining exponential stability
results by an appropriate embedding 
into the existing IQC framework introduced in~\Cref{secIntroIQCs}.
As it will get apparent, this allows for a more systematic approach to derive 
IQCs for robust exponential convergence guarantees from standard IQCs.

We first give a proper definition of robust exponential stability.
To this end, we need to step from the input-output framework
introduced in~\Cref{secIntroIQCs} to state-space descriptions. 
Let a state-space representation of~\eqref{eqFeedbackInterconnection}
be given by
\begin{subequations}
\begin{align}
	x\kk        &= A x\k   + B_1 w\k    + B_2 w_{\perf,k} \\
	y\k         &= C_1 x\k + D_{11} w\k + D_{12} w_{\perf,k} \\
	y_{\perf,k} &= C_2 x\k + D_{21} w\k + D_{22} w_{\perf,k} \\
	w\k         &= {\Delta(y)\k},
\end{align}\label{eqSSgeneralFeedbackInterconnection}
\end{subequations}
{with initial condition $ x_0 \in \real^n $}
and $ \yAdd = 0 $. Under the assumption that $ \Delta(0) = 0 $,
i.e., the origin is an equilibrium of~\eqref{eqSSgeneralFeedbackInterconnection} for $ w_\perf = 0 $,
we then have the following definition:
\begin{definition}[Robust exponential stability]\label{eqDefRobustExponentialStability}
We say that the origin {of~\eqref{eqSSgeneralFeedbackInterconnection}} is \emphDef{robustly exponentially stable
against $ \mathbf{\Delta} $ with rate $ \rho \in (0,1) $} if, for $ w_\perf = 0 $, it is 
globally {(uniformly)} exponentially stable with rate $ \rho $ for all $ \Delta \in \mathbf{\Delta} $,
i.e., there exists $ \eta > 0 $ such that $ \Vert x\k \Vert^2 \leq \eta \rho^{2k} \Vert x_0 \Vert^2 $
for all $ k \in \nat $, for each $ x_0 \in \real^{n} $ and each $ \Delta \in \mathbf{\Delta} $.
\end{definition}
Hence, if we can show that the origin is robustly exponentially stable 
against $ \mathbf{\Delta}(m,L) $ as defined in~\eqref{eqDefSetOfDelta} with rate $ \rho \in (0,1) $ for the transformed
dynamics~\eqref{eqRobustStandardForm}, then we conclude that the algorithm~\eqref{eqOptAlgo}
converges with rate $ \rho $ to $ \xOpt^\optSign = \Dopt^\dagger \optVar^\optSign $
for all $ \fObj \in {{\classObj{m}{L}}} $ and~\eqref{eqProblemFormulationConvergenceRate}
holds. {Since robust exponential stability and additional performance
specifications are independent of each other, for the following discussions
we neglect the additional performance channel {in~\eqref{eqRobustStandardForm}}.}

The idea for obtaining exponential stability results from $\ltwo$-stability
statements is to make use of a proper time-varying transformation of the 
signals. Such signal transformations have a long history in convergence rate analysis,
e.g., {so-called exponential weightings {have already been mentioned} in~\citet{desoer1975feedback} 
and similar signal transfomations are used in~\citet{antipin1994minimization} in the context of optimization}.
In particular, as in~\citet{boczar2015exponential}, for any $ \rho \in (0,1] $ and any $ p \in \natPos $,
we define two operators $ \rho_+, \rho_-: \lext^p \to \lext^p $ as
\begin{align}
	\big( \rho_+(y) \big)_k = \rho^k y_k, \quad \big( \rho_-(y) \big)_k = \rho^{-k} y_k,
\end{align}
where {$ k \in \nat $}.
For $ \rho \in (0,1) $, we then introduce the signal space
\begin{align}
	\ltworho^p = \lbrace y \in \lext^p : \rho_-(y) \in \ltwo^p \rbrace \subset \ltwo^p.
\end{align}
Note that any $ y \in \ltworho^p $ decays exponentially with rate $\rho$ {in the future},
i.e., there exists $ \eta >0 $ such that $ \Vert y\k \Vert \leq \eta e^{-\rho k} $ 
for all $ k \in \nat $.
Note further that 
\ifthenelse{\boolean{saveSpace}}{{$ \rho_+ : \ltwo^p \to \ltworho^p$, $\rho_- : \ltworho^p \to \ltwo^p $}}{
\begin{align}
	\rho_+ : \ltwo^p \to \ltworho^p \qquad
	\rho_- : \ltworho^p \to \ltwo^p
\end{align}}
{are bijective and inverse to each other, i.e.,  
$ \rho_+ \circ \rho_- = \text{id}_{\ltworho} $ and $ \rho_- \circ \rho_+ = \text{id}_{\ltwo} $.}
{Hence, if a transformed signal $ \tilde{y} = \rho_-(y) $ resides in $ \ltwo $,
then we can conclude that the original signal $ y = \rho_+(\tilde{y}) $
is in $ \ltworho $ and thus decays exponentially.
Remembering that $\ltwo$-properties can be deduced employing
standard IQC methods, this captures the main idea of the signal transformation.

Consider~\eqref{eqRobustStandardForm} and let $ \tilde{y} = \rho_-(y) $, 
$ \tilde{w} = \rho_-(w) $, $ \tilde{\xTrafo} = \rho_-(\xTrafo) $.
Neglecting the performance channel, the modified feedback 
interconnection then has a state-space representation
\begin{subequations}
\begin{align}
	\tilde{\xTrafo}\kk &= \rho^{-1} (\Aopt + m\Bopt \Copt) \tilde{\xTrafo}\k + \rho^{-1} \Bopt \tilde{w}\k \\
	\tilde{y}\k &= \Copt \tilde{\xTrafo}\k \\
	\tilde{w}\k &= \rho^{-k} \Delta_\fObj( \rho^k \tilde{y}\k ),
\end{align}\label{eqRobustStandardFormAfterLoopTrafo}
\end{subequations}
with initial condition $ \tilde{\xTrafo}_0 = \xTrafo_0 $.
The transfer function $ G_{\tilde{y}\tilde{w}} $ from $ \tilde{w} $ to $ \tilde{y} $ 
is then given by
\begin{align}
	G_{\tilde{y}\tilde{w}} \sim \big( \rho^{-1} (\Aopt + m\Bopt \Copt), \rho^{-1} \Bopt, \Copt, 0 \big),
	\label{eqDefTransformedTransferFunction}
\end{align}
thus \ifthenelse{\boolean{saveSpace}}{{$ G_{\tilde{y}\tilde{w} }(z) = \Gwtoy(\rho z), $}}{
\begin{align}
	G_{\tilde{y}\tilde{w} }(z) = \Gwtoy(\rho z),
\end{align}}
and we have (see~\Cref{figLoopTrafo})
\begin{subequations}
\begin{align}
	\tilde{y} &= G_{\tilde{y}\tilde{w} } \tilde{w} \\
	\tilde{w} &= {\big( \rho_- \circ \Delta_\fObj \circ \rho_+ \big) (\tilde{y})}.
\end{align}\label{eqTransformedFeedbackInterconnection}
\end{subequations}}
\begin{figure}[t]
	\begin{center}
		\begin{tikzpicture}[>=latex]
			\draw[white] (-4,0) -- (4,0);
			\node[name=nodeUncertainty,draw,rectangle,minimum height=1.15cm,minimum width=1.15cm] { $ \Delta_\fObj $ };
			\node[name=nodeRhoMinus1,draw,rectangle,minimum height=1.15cm,minimum width=1.15cm,right=0.7cm of nodeUncertainty.east] { $ \rho_- $ };
			\node[name=nodeRhoPlus1,draw,rectangle,minimum height=1.15cm,minimum width=1.15cm,left=0.7cm of nodeUncertainty.west] { $ \rho_+ $ };
			\node[name=nodePlant,draw,rectangle,minimum height=1.15cm,minimum width=1.15cm,below=0.5cm of nodeUncertainty.south] { $ \Gwtoy $ };
			\node[name=nodeRhoPlus2,draw,rectangle,minimum height=1.15cm,minimum width=1.15cm,right=0.7cm of nodePlant.east] { $ \rho_+ $ };
			\node[name=nodeRhoMinus2,draw,rectangle,minimum height=1.15cm,minimum width=1.15cm,left=0.7cm of nodePlant.west] { $ \rho_- $ };
			\draw[->] (nodeRhoPlus1) -- (nodeUncertainty);
			\draw[->] (nodeUncertainty) -- (nodeRhoMinus1);
			\draw[->] (nodeRhoMinus1.east) -- ++(0.7cm,0) |- node[pos=0.24,anchor=west] {$\tilde{w}$} (nodeRhoPlus2);
			\draw[->] (nodeRhoPlus2) -- (nodePlant);
			\draw[->] (nodePlant) -- (nodeRhoMinus2);
			\path (nodeRhoMinus2.west) -- node[name=nodeSum,draw,circle,inner sep=0pt,pos=1] {$+$} ++(-0.7cm,0);
			\draw[->] (nodeRhoMinus2.west) -- (nodeSum);
			\draw[->] (nodeSum.north) |- node[pos=0.2,anchor=east] {$\tilde{y}$} (nodeRhoPlus1);
			\draw[<-] (nodeSum.west) -- node[pos=1,anchor=east] {$\yAddTrafo$} ++(-0.7cm,0);
			\draw[thick,dashed] ([xshift=-6pt,yshift=6pt] nodeRhoPlus1.north west) -| node[anchor=north west] {$\Delta_{\fObj,\rho}$} ([xshift=6pt,yshift=-6pt] nodeRhoMinus1.south east) -| ([xshift=-6pt,yshift=6pt] nodeRhoPlus1.north west);
			\draw[thick,dashed] ([xshift=-6pt,yshift=6pt] nodeRhoMinus2.north west) -| node[anchor=north west,pos=0.85] {$G_{\tilde{y}\tilde{w}}$} ([xshift=6pt,yshift=-6pt] nodeRhoPlus2.south east) -| ([xshift=-6pt,yshift=6pt] nodeRhoMinus2.north west);
		\end{tikzpicture}
	\end{center}
	\caption{Modified feedback interconnection.}\label{figLoopTrafo}
\end{figure}
{In accordance with~\eqref{eqDefSetOfDelta}, we define the set of
transformed uncertainties as
\begin{align}
	\mathbf{\Delta}_{\rho}(m,L) = \big\lbrace \Delta_{\fObj,\rho} = \rho_- \circ \Delta_\fObj \circ \rho_+ \sepSet & \Delta_\fObj \in \mathbf{\Delta}(m,L) \big\rbrace. \label{eqDefSetOfDeltaRhoUnstructured} 
\end{align}}
The interconnection~\eqref{eqTransformedFeedbackInterconnection} has the same structure as the
standard feedback interconnection~\eqref{eqFeedbackInterconnection}; thus, {given that $ \Delta_{\fObj,\rho} $ is a bounded operator 
mapping $ \ltwo^p $ to $ \ltwo^p $, we may apply~\Cref{lemmaIQC}
to show $\ltwo$-stability of~\eqref{eqTransformedFeedbackInterconnection}
and deduce exponential stability of~\eqref{eqRobustStandardForm}.
We make this more precise in the following result.}
\begin{theorem}\label{lemmaIQCexponentialStability}
Consider~\eqref{eqRobustStandardForm} and let the transfer functions be
	defined according to~\eqref{eqStateSpaceRepresentations}.
	Let {$ L \geq m > 0 $ } be given and fix $ \rho \in (0,1) $.
	Suppose that all eigenvalues of \mbox{$ \Aopt + m \Bopt \Copt $} are located
	in the open disk of radius $ \rho $. Let some set of
	multipliers $ {\mathbf{\Pi}(\rho)} \subset \mathcal{RL}_{\infty}^{2p \times 2p} $
	parametrized by $ \rho $ and some set of performance multipliers $ \mathbf{\Pi}_\perf $ 
	in the form~\eqref{eqAssSetPerformanceMultipliers} be given. Suppose that
	\begin{enumerate}[leftmargin=*]
		\item the interconnection depicted in~\Cref{figLoopTrafo} is well-posed for all $ \Delta_{\fObj,\rho} \in \mathbf{\Delta}_\rho(m,L) $;
		\item each $ \Delta_{\fObj,\rho} \in  \mathbf{\Delta}_\rho(m,L) $ satisfies the IQC defined by {$ \Pi $} for all $ {\Pi} \in {\mathbf{\Pi}(\rho)} $.
	\end{enumerate}
	If there exist $ \Pi_1 \in {\mathbf{\Pi}(\rho)} $, $ \Pi_2 \in {\mathbf{\Pi}(1)} $ and $ \Pi_\perf \in \mathbf{\Pi}_\perf $
	such that
	\begin{subequations}
	\begin{align}
		\left[
		\vphantom{\begin{array}{ll} G_{\tilde{y}\tilde{w}} \\ I \end{array}} \star \right]^\complConj 
		\Pi_1
		\begin{bmatrix} G_{\tilde{y}\tilde{w}} \\ I \end{bmatrix}
		&\fdi{\prec} 0 \label{eqFDIRobustExpStab} \\
		\left[
		\vphantom{
		\begin{array}{ll}
			\Gwtoy & \Gwperftoy \\
			I & 0 \\
			0 & I \\
			\Gwtoyperf & \Gwperftoyperf
		\end{array}} \star \right]^\complConj
		\begin{bmatrix}
			\Pi_2 & 0 \\ 0 & \Pi_\perf
		\end{bmatrix}
		\left[
		\begin{array}{ll}
			\Gwtoy & \Gwperftoy \\
			I & 0 \\ \hline
			0 & I \\
			\Gwtoyperf & \Gwperftoyperf
		\end{array}\right]
		&\fdi{\prec} 0,
		\label{eqFDIRobustPerformanceExpStab}
	\end{align}
	\end{subequations}
	where $ G_{\tilde{y} \tilde{w}} $ is defined in~\eqref{eqDefTransformedTransferFunction},
	then the origin is robustly exponentially
	stable against $ \mathbf{\Delta}(m,L) $ with rate $ \rho $ for~\eqref{eqRobustStandardForm}
	and it achieves robust performance w.r.t.~$ \Pi_p $ against $ \mathbf{\Delta}(m,L) $.
\end{theorem}
\ifthenelse{\boolean{proofsAppendix}}{A proof can be found in~\Cref{secProofIQCexponentialStability}.}{\begin{proof} \begin{figure}[t]
	\centering
	\newlength\boxSize
	\setlength\boxSize{0.8cm}
	\begin{tikzpicture}[>=latex]
		\node[name=nodeEll2rhoLeft,anchor=center,minimum width=\boxSize,minimum height=\boxSize] {$ \ltworho$};
		\node[name=nodeEll2rhoRight,right=1.5cm of nodeEll2rhoLeft,anchor=center,minimum width=\boxSize,minimum height=\boxSize] {$ \ltworho$};
		\node[name=nodeEll2Left,below=1.5cm of nodeEll2rhoLeft,anchor=center,minimum width=\boxSize,minimum height=\boxSize] {$\ltwo$};
		\node[name=nodeEll2Right,right=1.5cm of nodeEll2Left,anchor=center,minimum width=\boxSize,minimum height=\boxSize] {$\ltwo$};
		\draw[->] (nodeEll2rhoLeft) -- node[anchor=south,pos=0.5] {$ \Delta_\fObj $} (nodeEll2rhoRight);
		\draw[->] (nodeEll2rhoLeft) -- node[anchor=east,pos=0.5] {$\rho_+$} (nodeEll2Left);
		\draw[->] (nodeEll2Left) -- node[anchor=south,pos=0.5] {$ \Delta_{\fObj,\rho} $} (nodeEll2Right);
		\draw[->] (nodeEll2Right) -- node[anchor=west,pos=0.5] {$ \rho_- $} (nodeEll2rhoRight);
	\end{tikzpicture}
	\caption{A schematic overview of the maps and signal spaces in the transformed loop {depicted in}~\Cref{figLoopTrafo}.}
	\label{figSchemeSets}
\end{figure}
{In order to apply~\Cref{lemmaIQC} to the transformed 
feedback interconnection, we first need to show boundedness
of the transformed uncertainty as well as shift the non-zero 
initial conditions of~\eqref{eqRobustStandardFormAfterLoopTrafo}
to the signal $ \yAddTrafo $ in~\Cref{figLoopTrafo}. \\
\textbf{Boundedness of $\Delta_{\fObj,\rho}$.}
Consider the transformed uncertainty
\begin{align}
	\Delta_{\fObj,\rho} = \rho_- \circ \Delta_\fObj \circ \rho_+. 
	\label{eqDefDeltaRho}
\end{align}
We next show that $ \Delta_{\fObj,\rho} $ is a bounded operator 
mapping $ \ltwo^p $ to $ \ltwo^p $ since $ \Delta_\fObj $ itself
is bounded, static and time-invariant.
To this end, first note that, 
there exists $ \beta \in \realNonNeg $ such that $ \Vert \big( \Delta_\fObj(y) \big)_k \Vert = \Vert \Delta_\fObj(y_k) \Vert \leq \beta \Vert y_k \Vert $
for any $ k \in \nat $ and any $ y \in \ltwo^p $. With $ \ltworho^p \subset \ltwo^p $,
we then infer that for any $ y \in \ltworho^p $ we have
\begin{align}
	\Vert\big( \rho_-(\Delta_\fObj(y)) \big)_k \Vert = \rho^{-k} \Vert \Delta_\fObj(y_k) \Vert \leq \rho^{-k} \beta \Vert y_k \Vert.
	\label{eqBoundDelta}
\end{align}
With $ y \in \ltworho^p $, i.e., $ \sum_{k=0}^\infty \rho^{-2k} \Vert y_k \Vert^2 $ is finite,
we conclude that $ \rho_-\big( \Delta(y) \big) \in \ltwo^p $ for any 
$ y \in \ltworho^p $, hence $ \Delta_\fObj: \ltworho^p \to \ltworho^p $ and
$ \Delta_{\fObj,\rho}: \ltwo^p \to \ltwo^p $ as required, see also~\Cref{figSchemeSets}
for an overview of the relations. We next show that {$ \Delta_{\fObj,\rho} $}  
is bounded on $ \ltwo^p $. For $ \tilde{y} \in \ltwo^p $, we indeed have $ y \coloneqq \rho_+(\tilde{y}) \in \ltworho $
and thus, using~\eqref{eqBoundDelta},
\begin{align}
	\Vert \Delta_{\fObj,\rho}(\tilde{y}) \Vert &= \Vert \rho_{-} \big( \Delta_\fObj(y) \big) \Vert  \leq \beta \Vert \rho_-(y) \Vert = \beta \Vert \tilde{y} \Vert.
\end{align}
\textbf{Shift of initial conditions.}}
{We next need to shift the non-zero initial conditions of~\eqref{eqRobustStandardFormAfterLoopTrafo}
to the signal $ \yAddTrafo $ in~\Cref{figLoopTrafo}. To this end, let $ \tilde{y}_k $, $ k \in \nat $,
denote the output of~\eqref{eqRobustStandardFormAfterLoopTrafo} for $ \tilde{\xTrafo}_0 = \xTrafo_0 $.
We then have that $ \tilde{y}_k $ also is the output of 
\begin{subequations}
\begin{align}
	\tilde{x}\kk &= \rho^{-1} (\Aopt + m\Bopt \Copt) \tilde{x}\k + \rho^{-1} \Bopt \tilde{w}\k \\
	\tilde{y}\k &= \Copt \tilde{x}\k + \yAddTrafo{}_{,k} \\
	\tilde{w}\k &= \rho^{-k} \Delta_\fObj( \rho^k \tilde{y}\k )
\end{align}
\end{subequations}
with initial condition $ \tilde{x}_0 = 0 $ and
\ifthenelse{\boolean{saveSpace}}{{$ \yAddTrafo{}_{,k} = \rho^{-k} \Copt (\Aopt + m \Bopt \Copt)^k \tilde{\xTrafo}_0, $}}{
\begin{align}
	\yAddTrafo{}_{,k} = \rho^{-k} \Copt (\Aopt + m \Bopt \Copt)^k \tilde{\xTrafo}_0,
\end{align}}
$ k \in \nat $. We infer that $ \yAddTrafo \in \ltwo^p $ since
$ \rho^{-1} (\Aopt + m \Bopt \Copt) $ has all eigenvalues in the
open unit disk. \\
We are now ready to apply~\Cref{lemmaIQC}.}
Robust performance follows directly from~\Cref{lemmaIQC}. For robust exponential stability, we note that, by~\Cref{lemmaIQC},~\eqref{eqFDIRobustExpStab} implies robust stability of the transformed
loop from~\Cref{figLoopTrafo} against $ \mathbf{\Delta}_\rho(m,L) $.
Hence, $ \tilde{y} $ and $ \tilde{w} $ in~\eqref{eqRobustStandardFormAfterLoopTrafo}
reside in $ \ltwo $. Together with the assumption that all eigenvalues of $ \Aopt + m \Bopt \Copt $ are located
in the open disk of radius $ \rho $, we infer that also $ \tilde{\xi} \in \ell_2^{np} $
in~\eqref{eqRobustStandardFormAfterLoopTrafo}. With $ y = \rho_+(\tilde{y}) $,
$ w = \rho_+(\tilde{w})$, $ \xi = \rho_+(\tilde{\xi}) $, and since
$ \rho_+ $ maps $ \ell_2 $-signals to $ \ltworho $-signals,
we conclude that $ y, w, \xi $ in~\eqref{eqRobustStandardForm} reside
in $ \ltworho $, hence the exponential decay follows.
 \end{proof}}
In terms of the original problem described in~\Cref{secProblemFormulation},
the previous Theorem provides a result that allows to
analyze and design optimization algorithms of the form~\eqref{eqOptAlgo}
in terms of convergence rates as well as performance.
{We note that the {frequency domain inequality (FDI) for} performance~\eqref{eqFDIRobustPerformanceExpStab} is
equivalent to the FDI~\eqref{eqFDIRobustPerformance} in the classical result~\Cref{lemmaIQC}
and exponential convergence is captured by the FDI~\eqref{eqFDIRobustExpStab}.
Compared to~\Cref{lemmaIQC}, we then {need to have} an adapted set of IQCs
$ {\mathbf{\Pi}(\rho)} $ for the set of transformed uncertainties $ \mathbf{\Delta}_\rho(m,L) $;
determining such a set will be the main subject of the following section.
Following the steps~\ref{itemStep1}--\ref{itemStep3},
we then also discuss how to reformulate the two FDIs~\eqref{eqFDIRobustExpStab},~\eqref{eqFDIRobustPerformanceExpStab}
as matrix inequalities, thereby paving the way to an implementable solution.}
\begin{remark}
In principle, the approach is not limited to exponential convergence rates.
	By employing suitable time-varying transformations, it is expected that similar results
	can be obtained for other types of convergence rate specifications{, for example
	to obtain polynomial convergence rate guarantees.}
\end{remark}

\subsection{IQCs for the class $ \unboldmath{\classObj{m}{L}}$} 
\label{subsecIQCsForTheClassOfUncertainties}
{The main goal of this section is to determine a set of IQCs valid
for the class of transformed uncertainties $ \mathbf{\Delta}_{\rho}(m,L) $ (step~\ref{itemStep2}).}
In the following we show how to systematically
derive such IQCs from {well-known} {classical} IQCs for the original set of uncertainties
$ \deltaUnstr(m,L) $. This is in contrast to {existing
results providing similar IQCs such as~\citet{lessard2016analysis,boczar2017exponential,fazlyab2018iqc}.
Since our approach (in particular~\Cref{lemmaMultipliersTransformedLoop}) allows to directly employ results from classical
IQC theory, we are able to
not only recover the modified IQCs from the {latter} references
but also easily extend them. In particular, we provide
a novel class of anticausal Zames-Falb multipliers
and show in~\Cref{secStructureExploitingAlgorithms}
how additional structural properties of the uncertain operator
can be captured by IQCs.}

With $ \deltaUnstr(m,L) $ being a set of slope-restricted operators, we can employ 
Zames-Falb IQCs that have been extensively studied in the literature, see, e.g., 
\citet{zames1969stability}, \citet{willems1971analysis}, \citet{mancera2005multipliers}.
While there exist other related IQCs for slope-restricted operators based
on discrete-time variants of the Popov criterion such as the Tsypkin or the Jury-Lee criteria,
quite recently it has been shown in~\citet{fetzer2017absolute} that, in the discrete-time 
setup, these are included in the set of Zames-Falb IQCs{;} thus we concentrate
on this class here.
{We give a brief wrap up of those classical results in the sequel.
In the literature, Zames-Falb IQC are commonly stated in the time domain;
following~\Cref{remarkIQCTimeDomain} this amounts to classifying the corresponding
operator $ \opTimeDomain $. }
To this end, we {utilize} the following definition introduced in~\citet{willems1971analysis}.
{
\begin{definition}[Doubly hyperdominant matrix]\label{defDoublyHyperdominant}
A matrix $ M  = [ m_{ij} ]_{i,j \in \lbrace 0,1,\dots,r \rbrace } $, $ r \in \nat $,
	is said to be \emphDef{doubly hyperdominant} if 
	\begin{align}
		m_{ij} \leq 0 \text{ for } i \neq j \text{ and } \sum\limits_{k=0}^{r} m_{kj} \geq 0, \sum\limits_{k=0}^{r} m_{ik} \geq 0
		\label{eqDefDoublyHyperdominant}
	\end{align}
	for all $ i,j \in \lbrace 0,1,\dots,r \rbrace $.
\end{definition}}
{Let $ \trunc{T}: \lext^p \to \lf^p $, $ T \in \nat $, denote the truncation
operator, i.e., $ \trunc{T}(y) = ( y_0, y_1, \dots, y_T, 0, 0, \dots ) $ for $ y \in \lext^p $.
{An infinite (block) matrix $ M = ( M_{ij} )_{i,j \in \nat} $, $ M_{ij} \in \real^{p \times p } $,
defines a linear operator $ M: \lf^p \to \lext^p $ in a natural fashion; the adjoint operator $ M^\top: \lf^p \to \lext^p $ is then defined by the matrix transpose. Doubly hyperdominance for such 
operators is then defined as follows:}
{
\begin{definition}[Doubly hyperdominant operator]
We call a linear operator $ M: \lf^p \to \lext^p $, 
	$ p \in \natPos $, \emph{doubly hyperdominant} if {the associated matrix of the truncated operator}
	$ M_T = \trunc{T} M \trunc{T} $ is doubly hyperdominant
	according to~\Cref{defDoublyHyperdominant} for each $ T \in \nat $.
\end{definition}}
Let $ \setDoublyHyp $ denote the set of all {infinite matrices defining} doubly hyperdominant operators.}
{The following result then holds}\ifthenelse{\boolean{proofsAppendix}}{, a proof is provided in~\Cref{secProofZFGeneralTimeDomain}.}{}
\begin{lemma}\label{lemmaZFGeneralTimeDomain}
Let {$ L \geq m > 0 $} be given and let $ \deltaUnstr(m,L) $ be defined according to~\eqref{eqDefSetOfDelta}.
	Let $ M : {\lf^p} \to {\lext^p} $ be {a linear operator {defined by the infinite matrix} $ M = \bar{M} \otimes I_p $, $ \bar{M} \in \setDoublyHyp $.} 
	Then, for all $ y \in {\ltwo^p} $, {all $ \Delta_\fObj \in \deltaUnstr(m,L) $ {and all $ T \in \nat $}, we have}
	\begin{align}
		\left\langle \begin{bmatrix} y \\ \Delta_\fObj(y) \end{bmatrix}, {\trunc{T}} \trafoSector^\top \begin{bmatrix} 0 & {M}^\top \\ {M} & 0 \end{bmatrix} \trafoSector {\trunc{T}} \begin{bmatrix} y \\ \Delta_\fObj(y) \end{bmatrix} \right\rangle \geq 0 \label{eqIneqZamesFalbTimeDomain} 
	\end{align}
	{with}
	\begin{align}
		\qquad \trafoSector = {\begin{bmatrix} (L-m) \id & -\id \\ 0 & \id \end{bmatrix}}. \label{eqDefTrafoSector}
	\end{align}
\end{lemma}
\ifthenelse{\boolean{proofsAppendix}}{
}
{\begin{proof} {The proof boils down to extending the 
result for monotone nonlinearities from~\citet[Lemma 2]{mancera2005multipliers} 
to the case of slope-restricted nonlinearities.} 
The following {argumentations are} similar to {those} provided 
in the proof of~\cite[{Theorem 1}]{amato2001repeated} where only 
scalar nonlinearities are considered. 
We first repeat the definition of a monotone function.
\begin{definition}[Monotone function]
A continuously differentiable function $ f: \real^p \to \real^p $ is said to be \emph{monotone}
	if (1) it is conservative, i.e., there exists $ F: \real^p \to \real $
	such that $ \nabla F(x) = f(x) $ for all $ x \in \real^p $ and (2)
	it fulfills
	\begin{align}
		\langle f(x) - f(y), x - y \rangle \geq 0
		\label{eqIneqMonotoneFunction}
	\end{align}
	for all $ x, y \in \real^p $.
\end{definition}
{It is well-known that a function is monotone
if it is the gradient of some convex function.}
Monotone operators $ \varphi: \lext^p \to \lext^p $ are then defined as
operators of the form
\begin{align}
	\big( \varphi(y) \big)_k = f(y_k),
\end{align}
where $ y = {[y_0,y_1,\dots]} \in \lext^p $, $ k \in {\nat} $,
and $ f: \real^p \to \real^p $ is a monotone function.
If $M$ adheres to the conditions given in~\Cref{lemmaZFGeneralTimeDomain},
by~\citet[Lemma 2]{mancera2005multipliers} it is then known that
for any monotone operator $\varphi: \ltwo^{p} \to \ltwo^p $ we have
\begin{align}
	\langle {\trunc{T} y}, M^\top \varphi( {\trunc{T} y} ) \rangle \geq 0
	\label{eqIneqZamesFalbTimeDomainMonotoneNonl}
\end{align}
for all $ y \in \ltwo^p $, $ T \in \nat $. {To use this result, we next rewrite~\eqref{eqIneqZamesFalbTimeDomain} in the form~\eqref{eqIneqZamesFalbTimeDomainMonotoneNonl}.}
To this end, {consider~\eqref{eqIneqZamesFalbTimeDomain}, \eqref{eqDefTrafoSector}
with $L$ replaced by $ \tilde{L} \coloneqq L + \varepsilon $, $ \varepsilon > 0 ${. By}
simple calculations as well as the fact that $ \Delta_\fObj $ is static and time-invariant, {\eqref{eqIneqZamesFalbTimeDomain} is}}
\begin{align}
	2 \langle {(\tilde{L}-m)\trunc{T} y - \Delta_\fObj(\trunc{T} y)}, M^\top \Delta_\fObj({\trunc{T} y}) \rangle \geq 0.
	\label{eqIneqZamesFalbTimeDomainRewritten}
\end{align}
{In view of~\eqref{eqIneqZamesFalbTimeDomainMonotoneNonl}, for any 
$ y \in \ltwo^p $ we define the auxiliary sequence $ \tilde{y} \in \ltwo^p $ by
\begin{align}
	\tilde{y}\k \coloneqq (\tilde{L}-m) y\k - \Delta_\fObj(y\k), \quad k \in \nat.
\end{align}
Let $ \phi: \real^p \to \real $ be defined as
\begin{align}
	\phi(y) 
	&= \tfrac{1}{2} (\tilde{L}-m) y^\top y - \big( \fObj(y) - \tfrac{1}{2} m y^\top y \big) \ifthenelse{\boolean{singleColumn}}{}{\nonumber \\
	&}= \tfrac{1}{2} \tilde{L} y^\top y - \fObj(y),
\end{align}
{which gives} $ \tilde{y} = \nabla \phi(y) $ {and
\begin{align}
	\Delta_\fObj(y) = \nabla\phi(y) - ( \tilde{L}-m ) y.
	\label{eqDeltaHrewritten}
\end{align}
Then} \eqref{eqIneqZamesFalbTimeDomainRewritten}
reads as
\begin{align}
	2 \langle \trunc{T} \tilde{y}, M^\top \Delta_\fObj(\trunc{T}y) \rangle \geq 0.
	\label{eqIneqZamesFalbTimeDomainRewritten2}
\end{align}
Note that $ \fObj \in {\classObj{m}{L}} $ implies $ \phi \in {\classObj{\varepsilon}{\tilde{L}-m}} $; we then infer from~\citet[Theorem 26.6, Lemma 26.7]{rockafellar1970convex} 
that $ \nabla \phi $ has a well-defined inverse $ (\nabla \phi)^{-1} $
and $ (\nabla \phi)^{-1} $ is itself the gradient of a strictly convex function,
namely the Legendre conjugate of $ \phi $, cf.~\citet[Theorem 26.5]{rockafellar1970convex}.
Thus, we can write~\eqref{eqIneqZamesFalbTimeDomainRewritten2} as
\begin{align}
	2 \langle \trunc{T} \tilde{y}, M^\top \Delta_\fObj\big( (\nabla \phi)^{-1}( \trunc{T} \tilde{y}) \big) \rangle \geq 0.
\end{align}
Comparing with~\eqref{eqIneqZamesFalbTimeDomainMonotoneNonl}, we aim to show that
the operator $ \varphi \coloneqq \Delta_\fObj \circ (\nabla \phi)^{-1} $ is monotone.
To this end, }
first note that, since $ (\nabla \phi)^{-1} $ is the gradient of a strictly convex function,
there exists some function $ F: \real^p \to \real $ such that 
$ \big({\varphi}(\tilde{y})\big)_k = \nabla F(\tilde{y}_k) $.
We next show that $ \nabla F $ fulfills~\eqref{eqIneqMonotoneFunction},
hence {$ \tilde{w} $} is a monotone operator. 
{By~\eqref{eqDeltaHrewritten} we infer that 
\begin{align}
	\varphi(\tilde{y}) = \Delta_\fObj\Big( (\nabla \phi)^{-1}(\tilde{y}) \Big) = \tilde{y} - (\tilde{L}-m) \big(\nabla \phi\big)^{-1}(\tilde{y}).
\end{align}
Thus,} the Hessian of $F$ is given by
\begin{align}
\nabla^2 F( \tilde{y}_k ) = (\tilde{L}-m) \Big( \nabla^2 \phi\big( (\nabla\phi)^{-1}(\tilde{y}_k) \big) \Big)^{-1} - I.
\end{align}
Now, since $ \phi \in {\classObj{\varepsilon}{\tilde{L}-m}} $, we have
by~\citet[Theorem 2.1.6, Theorem 2.1.11]{nesterov2004introductory} that
$ \varepsilon I \preceq \nabla^2 \phi(y) \preceq (L-m+\varepsilon) I $
for all $y \in \real^p $, and hence $ \nabla^2 F(\tilde{y}_k) \succeq (\tfrac{\tilde{L}-m}{\tilde{L}-m} - 1) I = 0 $
for all $\tilde{y}_k \in \real^p $. Consequently, $F$ is convex \citep[Theorem 2.1.4]{nesterov2004introductory}
and we conclude that {$ \varphi $} is a monotone operator,
thus~\eqref{eqIneqZamesFalbTimeDomainRewritten} holds
for all $ y \in \ltwo^p $, $ T \in \nat $, { $\varepsilon > 0 $}.
{For arbitrary but fixed $ y \in \ltwo^p $, $ T \in \nat $, we can then 
take the limit $ \varepsilon \to 0 $ in~\eqref{eqIneqZamesFalbTimeDomainRewritten}
since the left-hand side is continuous in $ \varepsilon $ and infer
that~\eqref{eqIneqZamesFalbTimeDomain} holds, thus concluding the proof.}

 \end{proof} }
\begin{remark}
If $ \deltaUnstr(m,L) $ has the additional property that {$ \Delta(y) = -\Delta(-y) $} 
	for any $ y \in \ltwo^p $, $ \Delta \in \deltaUnstr(m,L) $, i.e., each $ \Delta $
	is an odd operator, then~\eqref{eqIneqZamesFalbTimeDomain} persists to hold
	for $ \bar{M} $ being
	doubly dominant, see~\citet{willems1971analysis} for a definition.
\end{remark}
{The latter result is more or less known but, to the best of our knowledge,
not exactly to be found in literature. In particular,
besides generalizing~\citet[Lemma 2]{mancera2005multipliers} 
to slope-restricted nonlinearities, with~\eqref{eqIneqZamesFalbTimeDomain}
we provide a so-called hard IQC (instead a soft IQC) and thus can 
omit the boundedness requirements on the operator $M$.
If $M$ is a bounded operator, then the limit of~\eqref{eqIneqZamesFalbTimeDomain} as $ T \to \infty $
is well-defined and the inequality corresponds to a soft IQC{, i.e.,
the following inequality holds for all $ \Delta_\fObj \in \deltaUnstr(m,L) $
{and all $ y \in \ltwo^p $:}
\begin{align}
	\left\langle \begin{bmatrix} y \\ \Delta_\fObj(y) \end{bmatrix}, \trafoSector^\top \begin{bmatrix} 0 & {M}^\top \\ {M} & 0 \end{bmatrix} \trafoSector \begin{bmatrix} y \\ \Delta_\fObj(y) \end{bmatrix} \right\rangle \geq 0. \label{eqZamesFalbIQCTimeDomain} 
\end{align}}
These extensions are important to derive an analog
of the soft IQC~\eqref{eqZamesFalbIQCTimeDomain}
for the transformed uncertainty $ \Delta_{\fObj,\rho} $.}
In the next Lemma we
present a general result building the basis for deriving IQCs for
transformed uncertainties from standard ones.
\ifthenelse{\boolean{proofsAppendix}}{A proof is provided in~\Cref{secProofLemmaMultipliersTransformedLoop}.}{}
{
\begin{lemma}\label{lemmaMultipliersTransformedLoop}
Let $ \mathbf{\Delta} \subset \lbrace \Delta \in \mathcal{L}(\ltwo^p,\ltwo^p) \sepSet \Delta(\ltworho^p) \subset \ltworho^p \rbrace $
	be some set of operators. 
	Suppose we have given $ \opTimeDomain: \lf^{{2p}} \to \lext^{{2p}} $ 
	such that $ \tilde{\opTimeDomain} = \rho_+ \circ \opTimeDomain \circ \rho_+ $ is bounded on $ \ltwo^{{2p}} $ and
	\begin{align}
		\left\langle \begin{bmatrix} y \\ \Delta(y) \end{bmatrix}, \trunc{T} \opTimeDomain \trunc{T} \begin{bmatrix} y \\ \Delta(y) \end{bmatrix} \right\rangle \geq 0 
		\label{eqIQClemmaMultipliersTransformedLoop}
	\end{align}
	for all $ y \in \ltwo^p $, $ \Delta \in \mathbf{\Delta} $, $ T \in \nat $.
	Then
	\begin{align}
		\left\langle \begin{bmatrix} y \\ \Delta_\rho(y) \end{bmatrix}, \tilde{\opTimeDomain} \begin{bmatrix} y \\ \Delta_\rho(y) \end{bmatrix} \right\rangle \geq 0
		\label{eqIQCdeltaRho}
	\end{align}
	holds for all $ y \in \ltwo^p $ and any $ \Delta_\rho = \rho_- \circ \Delta \circ \rho_+ $, $ \Delta \in \mathbf{\Delta} $.
\end{lemma}}
\ifthenelse{\boolean{proofsAppendix}}{
}
{\begin{proof} {Let $ \opTimeDomain_T = \trunc{T} \opTimeDomain \trunc{T} $, $ T \in \nat $.} {
Since $ \rho_+(\ltwo^p) = \ltworho^p $, we infer from~\eqref{eqIQClemmaMultipliersTransformedLoop} {that}
\begin{align}
	\left\langle \begin{bmatrix} \rho_+(y) \\ \Delta(\rho_+(y)) \end{bmatrix}, {\opTimeDomain_T} \begin{bmatrix} \rho_+(y) \\ \Delta(\rho_+(y)) \end{bmatrix} \right\rangle \geq 0 
\end{align}
for all $ y \in \ltwo^p $, $ \Delta \in \mathbf{\Delta} $, 
{$T \in \nat$}.
With $ \rho_+ \circ \rho_- = \text{id}_{\ltworho} $ this is 
equivalently formulated as
\begin{align}
	&\left\langle \begin{bmatrix} \rho_+(y) \\ \rho_+\Big( \rho_- \big(\Delta(\rho_+(y)) \big)\Big) \end{bmatrix}, {\opTimeDomain_T} \begin{bmatrix} \rho_+(y) \\ \rho_+\Big( \rho_- \big(\Delta(\rho_+(y)) \big)\Big) \end{bmatrix} \right\rangle \nonumber \\
	=~& \left\langle  \begin{bmatrix} y \\ \Delta_\rho(y) \end{bmatrix}, \rho_+\Big( {\opTimeDomain_T} \big( \rho_+ (  \begin{bmatrix} y \\ \Delta_\rho(y) \end{bmatrix} ) \big) \Big) \right\rangle \nonumber \\
	=~& \left\langle  \begin{bmatrix} y \\ \Delta_\rho(y) \end{bmatrix}, {\trunc{T} \Big( \tilde{\opTimeDomain} \big( \trunc{T} (}  \begin{bmatrix} y \\ \Delta_\rho(y) \end{bmatrix} ) \big) \Big) \right\rangle
	\geq 0 
	\label{eqInequalityProofMultipliersTransformedLoop}
\end{align}
for all $ y \in \ltwo^p $, $ \Delta \in \mathbf{\Delta} $, {$T \in \nat$}.
{Since $ \Delta (\ltworho^p) \subset \ltworho^p $, we have $ \Delta_\rho: \ltwo^p \to \ltwo^p $;
hence, by the assumption that $ \tilde{\opTimeDomain} $ is bounded on $ \ltwo^p $,
{we can take the limit $ T \to \infty $ in~\eqref{eqInequalityProofMultipliersTransformedLoop}
to obtain that}~\eqref{eqIQCdeltaRho} holds.}}
 \end{proof} }
{Note that~\eqref{eqIQClemmaMultipliersTransformedLoop} {provides a 
hard IQC for the transformed uncertainty as required.} If $ \opTimeDomain: \ltwo^{2p} \to \ltwo^{2p} $ is bounded on $ \ltwo^{2p} $ 
and $ \tilde{\opTimeDomain} $ is bounded, {then we can replace 
$ \trunc{T} \opTimeDomain \trunc{T} $ by $ \opTimeDomain $ itself in~\eqref{eqIQClemmaMultipliersTransformedLoop},
i.e.,
\begin{align}	
	\left\langle \begin{bmatrix} y \\ \Delta(y) \end{bmatrix}, \opTimeDomain \begin{bmatrix} y \\ \Delta(y) \end{bmatrix} \right\rangle \geq 0 
\end{align}
implies~\eqref{eqIQCdeltaRho} {for $y\in\ltwo^p$.}}
However, typically either only $ \opTimeDomain $ or the transformed
operator $ \tilde{\opTimeDomain} $ is bounded.

We next employ~\Cref{lemmaMultipliersTransformedLoop} to derive IQCs for the set of transformed uncertainties $ \mathbf{\Delta}_{\rho}(m,L) $. 
In the sequel we limit ourselves to Toeplitz type operators $ \bar{M} $. 
To this end, we {introduce the shorthand notation $ \toep(R): \lf^p \to \lext^p $
defined as
\ifthenelse{\boolean{singleColumn}}{
\begin{align}
	\toep(R) 
	:= 
	\left[
	\begin{array}{ccccccc}
		R_0             & \dots  & R_{\dimAnticausal} & 0 & \dots \\ 
		R_{-1}          & R_0    & \dots              & R_{\dimAnticausal} & 0 & \dots \\
		R_{-2}          & R_{-1} & R_0                & \dots & R_{\dimAnticausal} & 0 & \dots \\
		\vdots          & \ddots & \ddots             & \ddots & \dots              & \ddots & \ddots 
	\end{array}\right], 
\end{align}}{
\begin{flalign}
	\toep(R) 
	:= 
	\left[
	\begin{shorterArray}{ccccccc}
		R_0             & \dots  & R_{\dimAnticausal} & 0 & \dots \\ 
		R_{-1}          & R_0    & \dots              & R_{\dimAnticausal} & 0 & \dots \\
		R_{-2}          & R_{-1} & R_0                & \dots & R_{\dimAnticausal} & 0 & \dots \\
		\vdots          & \ddots & \ddots             & \ddots & \dots              & \ddots & \ddots 
	\end{shorterArray}\right], \hspace*{-2em} &&
\end{flalign}}
where $ R = \begin{bmatrix} R_{-\dimCausal} & R_{-\dimCausal+1} & \dots & R_{\dimAnticausal} \end{bmatrix} \in \real^{p \times p(\dimCausal+1+\dimAnticausal)} $,
$ R_k \in \real^{p \times p} $, $ k = -\dimCausal,-\dimCausal+1,\dots,\dimAnticausal $.}
Note that $ \toep(R) $ is a causal operator if and only if $ \dimAnticausal = 0$.
A combination of~\Cref{lemmaZFGeneralTimeDomain} with~\Cref{lemmaMultipliersTransformedLoop}
then yields the following result\ifthenelse{\boolean{proofsAppendix}}{, a proof can be found in~\Cref{secProofZFRhoTimeDomain}.}{.}
\begin{lemma}\label{lemmaZFRhoTimeDomain}
Let {$ L \geq m > 0 $} and $ \dimCausal, \dimAnticausal \in \nat $ be given. Let $ \rho \in (0,1] $.
Let $ \deltaUnstr_\rho(m,L) $ be defined according to~\eqref{eqDefSetOfDeltaRhoUnstructured} and 
	let $ {M_{\rho}}: {\ltwo^p} \to {\ltwo^p} $ be {the} {bounded} linear operator with
	{
	\begin{align}
		{M_{\rho}} = \toep\big(\begin{bmatrix} M_{-\dimCausal} & M_{-\dimCausal+1} & \dots & M_{\dimAnticausal} \end{bmatrix} \big) \label{eqDefOperatorM}
	\end{align}}
	where $ M_{i} \in \real^{p \times p} $. Then, for all $ y \in \ltwo^p $, the
	{inequality~\eqref{eqZamesFalbIQCTimeDomain}}
	holds for all $ \Delta_{\fObj} \in \deltaUnstr_\rho(m,L) $ if {$ M = M_{\rho} $}, $ M_i = m_i I_p $ and
	\begin{align}
		\toep(\begin{bmatrix} \rho^{-\dimCausal} m_{-\dimCausal} & \dots & m_{0} & \dots & \rho^{\dimAnticausal} m_{\dimAnticausal} \end{bmatrix}) \label{eqConditionHyperdominantToeplitz} 
	\end{align}
	is doubly hyperdominant, i.e., 
	\begin{align}
		\sum\limits_{i=-\dimCausal}^{\dimAnticausal} m_i \rho^{-i} \geq 0, \qquad \sum\limits_{i=-\dimCausal}^{\dimAnticausal} m_i \rho^{i} \geq 0,
	\end{align}
	and $ m_i \leq 0 $ for all {$ i \in \lbrace -\dimCausal, \dots, \dimAnticausal \rbrace $, $ i \neq 0 $.}
\end{lemma}
\ifthenelse{\boolean{proofsAppendix}}{
}{
\begin{proof}
	{
Observe that $ M_\rho = \bar{M}_\rho \otimes I_p $ with
\ifthenelse{\boolean{saveSpace}}{{$ \bar{M}_{\rho} = \toep\big(\begin{bmatrix} m_{-\dimCausal} & m_{-\dimCausal+1} & \dots & m_{\dimAnticausal} \end{bmatrix} \big) $}}{
\begin{align}
	\bar{M}_{\rho} = \toep\big(\begin{bmatrix} m_{-\dimCausal} & m_{-\dimCausal+1} & \dots & m_{\dimAnticausal} \end{bmatrix} \big)
\end{align}}
and define $ M = \rho_- M_\rho \rho_- = ( \rho_- \bar{M}_\rho \rho_- ) \otimes I_p $.
We compute
\begin{align}
	\rho_- \bar{M}_\rho \rho_-
	&= 
	\left[
	\begin{shorterArray}{cccc}
		\bar{m}_{0}            & \rho^{-1} \bar{m}_{1} & \rho^{-2} \bar{m}_{2} & \dots \\
		\rho^{-1} \bar{m}_{-1} & \rho^{-2} \bar{m}_{0} & \rho^{-3} \bar{m}_{1} & \dots \\
		\ddots           & \ddots          & \ddots          & \ddots
	\end{shorterArray}
	\right].
\end{align}
Thus, $ \rho_- \bar{M}_\rho \rho_- \in \setDoublyHyp $ if, for each finite $ k \in \nat $,
\begin{subequations}
\begin{align}
	\rho^{-k-i} \bar{m}_i \leq 0 \text{ for } i \in \lbrace -\dimCausal, \dots, \dimAnticausal \rbrace, i \neq 0, \\
	\rho^{-2k} \sum\limits_{i=-\dimCausal}^{\dimAnticausal} \bar{m}_i \rho^{-i} \geq 0, \rho^{-2k} \sum\limits_{i=-\dimCausal}^{\dimAnticausal} \bar{m}_i \rho^{i} \geq 0.
\end{align}
\end{subequations}
With $ \rho > 0 $, we can get rid of the dependency
on $k$ such that the above conditions are equivalent to requiring that
the Toeplitz operator in~\eqref{eqConditionHyperdominantToeplitz}
is doubly hyperdominant.
We hence conclude by~\Cref{lemmaZFGeneralTimeDomain} that under these 
conditions the inequality~\eqref{eqIQClemmaMultipliersTransformedLoop} holds with
\begin{align}
	\opTimeDomain = \trafoSector^\top \begin{bmatrix} 0 & {M}^\top \\ {M} & 0 \end{bmatrix} \trafoSector,
\end{align}
{where $ \trafoSector $ is defined in~\eqref{eqDefTrafoSector} and}
$ M = ( \rho_- \bar{M}_\rho \rho_- ) \otimes I_p $. Additionally, $ \rho_+ \circ \opTimeDomain \circ \rho_+ $ 
is bounded on $ \ltwo^{2p} $, since $ M $ is bounded on $ \ltwo^p $ by assumption.
Hence, using also that $ \Delta_{\fObj,\rho}: \ltwo^p \to \ltwo^p $
as discussed in~\Cref{subsecExponentialConvergenceViaIQCs}, we may apply~\Cref{lemmaMultipliersTransformedLoop} and conclude
that 
\begin{align}
	\left\langle \begin{bmatrix} y \\ \Delta_{\fObj,\rho}(y) \end{bmatrix}, \rho_+ \trafoSector^\top \begin{bmatrix} 0 & {M}^\top \\ {M} & 0 \end{bmatrix} \trafoSector \rho_+ \begin{bmatrix} y \\ \Delta_{\fObj,\rho}(y) \end{bmatrix} \right\rangle \nonumber \\
	= \left\langle \begin{bmatrix} y \\ \Delta_{\fObj,\rho}(y) \end{bmatrix}, \trafoSector^\top \begin{bmatrix} 0 & {M_\rho}^\top \\ {M_\rho} & 0 \end{bmatrix} \trafoSector \begin{bmatrix} y \\ \Delta_{\fObj,\rho}(y) \end{bmatrix} \right\rangle
	\geq 0 
\end{align}
for any $ y \in \ltwo^p $, $ \Delta_{\fObj,\rho} \in \mathbf{\Delta}_\rho(m,L) $,
thus concluding the proof.}

\end{proof}
}
The Toeplitz operator {$M_\rho$} from~\eqref{eqDefOperatorM} is fully
described by {the (finite) set of matrices $M_i$, $ i \in \lbrace -\dimCausal, -\dimCausal+1, \dots, \dimAnticausal \rbrace $}, and, in view of the previous result,
we then define the set of admissible matrices as
\begingroup
\arraycolsep=1pt\def\arraystretch{1}
\ifthenelse{\boolean{singleColumn}}{
\begin{align}
	\mathbb{M}(\rho,\dimAnticausal,\dimCausal,p) =\big\lbrace \begin{bmatrix} M_{-\dimCausal} & M_{-\dimCausal +1} & \dots & M_{\dimAnticausal} \end{bmatrix} \sepSet 
	~&M_i = m_i I_p, \; m_i \leq 0 \text{ for all } i \neq 0 \nonumber \\
	& ( \textstyle{\sum_{i=-\dimCausal}^{\dimAnticausal} M_i \rho^{-i}} ) \mathbf{1} \geq 0, \nonumber \\
	& \mathbf{1}^\top ( \textstyle{\sum_{i=-\dimCausal}^{\dimAnticausal} M_i \rho^{i}} ) \geq 0 \big\rbrace. \label{eqDefSetMunstructured}
\end{align}}{
\begin{align}
	\mathbb{M}(\rho,\dimAnticausal,\dimCausal,p) =\big\lbrace~& \begin{bmatrix} M_{-\dimCausal} & M_{-\dimCausal +1} & \dots & M_{\dimAnticausal} \end{bmatrix} \sepSet \nonumber \\
	&M_i = m_i I_p, \; m_i \leq 0 \text{ for all } i \neq 0 \nonumber \\
	& ( \textstyle{\sum_{i=-\dimCausal}^{\dimAnticausal} M_i \rho^{-i}} ) \mathbf{1} \geq 0, \nonumber \\
	& \mathbf{1}^\top ( \textstyle{\sum_{i=-\dimCausal}^{\dimAnticausal} M_i \rho^{i}} ) \geq 0 \big\rbrace. \label{eqDefSetMunstructured}
\end{align}}
\endgroup
{We note that the above constraints can be formulated in a
simpler fashion in terms of the parameters $ m_i $; we deliberately keep this
{more general version since it is better suited for the extension to the structured case} 
discussed in~\Cref{secStructureExploitingAlgorithms},
cf.~\eqref{eqDefSetMrepeated}, \eqref{eqDefSetMnonRepeated}.}

{With frequency domain representations being more common,
we next briefly sketch how to obtain such a
representation of the time domain IQC~\eqref{eqZamesFalbIQCTimeDomain}.
To this end, note that a Toeplitz operator $ M = M_\rho: {\ltwo^p} \to {\ltwo^p} $ as in~\eqref{eqDefOperatorM}
defines a corresponding transfer matrix $ \tfzf_M(z) = \sum_{j=-\dimAnticausal}^{\dimCausal} M_{-j} z^{-j}$ by
\begin{align}
    \widehat{{M}y}(z) = \tfzf_M(z) \widehat{y}(z).
\end{align}
By Parseval's Theorem, it is then straightforward to formulate~\eqref{eqZamesFalbIQCTimeDomain}
equally well in the frequency domain and the conditions on
the operator $M$ (or $M_\rho$) then translate to conditions
on the corresponding transfer matrix $ \tfzf_M $. 
More precisely, the corresponding class of Zames-Falb multipliers is {given as}
\ifthenelse{\boolean{singleColumn}}{
\begin{align}
	\mathbf{\Pi}_{\Delta,\rho}^p(m,L) = 
	\big\lbrace \Pi = {\widehat{\trafoSector}}^\top \begin{bmatrix} 0 & \tfzf^\complConj \\ \tfzf & 0 \end{bmatrix} {\widehat{\trafoSector}} \sepSet 
	~&\tfzf(z) = \textstyle{\sum_{j=-\dimCausal}^{\dimAnticausal} M_{j} z^{{j}}}, \dimAnticausal, \dimCausal \in \mathbb{N}, \label{eqDefZamesFalbMultipliers} \\
	&\left[\begin{shortArray}{ccc} M_{-\dimCausal} & \dots & M_{\dimAnticausal} \end{shortArray}\right] \in \mathbb{M}(\rho,\dimAnticausal,\dimCausal,p) \big\rbrace, \nonumber
\end{align}}{
\begin{align}
	\mathbf{\Pi}_{\Delta,\rho}^p(m,L) = 
	\big\lbrace ~& \Pi = {\widehat{\trafoSector}}^\top \begin{bmatrix} 0 & \tfzf^\complConj \\ \tfzf & 0 \end{bmatrix} {\widehat{\trafoSector}} \sepSet \label{eqDefZamesFalbMultipliers} \\
	&\tfzf(z) = \textstyle{\sum_{j=-\dimCausal}^{\dimAnticausal} M_{j} z^{{j}}}, \dimAnticausal, \dimCausal \in \mathbb{N}, \nonumber \\
	&\left[\begin{shortArray}{ccc} M_{-\dimCausal} & \dots & M_{\dimAnticausal} \end{shortArray}\right] \in \mathbb{M}(\rho,\dimAnticausal,\dimCausal,p) \big\rbrace, \nonumber
\end{align}}
where $ \mathbb{M} $ is defined in~\eqref{eqDefSetMunstructured}
and $ \widehat{\trafoSector} \in \real^{2p \times 2p} $ given by 
\begin{align}
	\widehat{W} = \begin{bmatrix} (L-m) I & -I \\ 0  & I \end{bmatrix}
\end{align}
is the $z$-transform corresponding to the operator defined by~\eqref{eqDefTrafoSector}.
}
Summing up, we then have the following result
{concluding Step~\ref{itemStep2} in {our} procedure.}
\begin{theorem}\label{lemmaZFIQCrho}
Let $ L \geq m > 0 $ and let $ \mathbf{\Delta}_{\rho}(m,L) $
	be defined as in~\eqref{eqDefSetOfDeltaRhoUnstructured}.
	Then, for each $ \rho \in (0,1] $, $ \Delta_\rho $ satisfies the IQC defined by $ \Pi $
	for each $ \Delta_\rho \in \mathbf{\Delta}_{\rho}(m,L) $ and each $ \Pi \in \mathbf{\Pi}_{\Delta,\rho}^p(m,L) $
	{as defined in~\eqref{eqDefZamesFalbMultipliers}}.
\end{theorem}
\Cref{lemmaZFIQCrho} is a generalization of the IQCs presented in~\citet{boczar2015exponential}, \citet{lessard2016analysis}
also allowing for anticausal multipliers $ \Pi $, i.e., anticausal
transfer matrices $ \tfzf $ in~\eqref{eqDefZamesFalbMultipliers}.
The class of multipliers from~\citet{boczar2015exponential}, \citet{lessard2016analysis}
is then obtained by letting $ \dimAnticausal = 0 $.
{We note that an extension to anticausal multipliers has also
been proposed in~\citet{freeman2018noncausal} using a slightly different approach.
{However, the conditions derived in~\citet{freeman2018noncausal} are
more restrictive, thus leading to a smaller class of multipliers compared
to the one proposed in~\Cref{lemmaZFIQCrho}, see~\Cref{secAppendixDiscussionComparison} for
a {detailed} discussion. As another advantage, the presented}	
approach is based on standard results and
allows for easy extensions such as, e.g., incorporating structural properties
of the objective function, see~\Cref{secStructureExploitingAlgorithms}.}

\subsection{{Multiplier parametrization}}\label{subsecMultiplierParametrization}
{Having a suitable set of {IQCs} available {as it is provided by~\Cref{lemmaZFIQCrho}}, applying~\Cref{lemmaIQCexponentialStability}
then basically amounts to checking the FDIs~\eqref{eqFDIRobustExpStab}, \eqref{eqFDIRobustPerformanceExpStab}.
In order to derive efficiently implementable conditions in terms of matrix 
inequalities employing the KYP-Lemma, we {need to factorize the multipliers
and determine} state-space representations
of all transfer functions in these FDIs (Step~\ref{itemStep3}).
{In particular, utilizing the class of multipliers~\eqref{eqDefZamesFalbMultipliers},
for a factorization of $ \Pi \in \mathbf{\Pi}_{\Delta,\rho}^p(m,L) $ we need to 
find a proper and stable transfer matrix $\psi_{\Delta} $ and
a matrix $ M_\Delta $ such that
\begin{align}
	\Pi
	=
	\psi_{\Delta}^\complConj M_{\Delta} \psi_{\Delta}. \label{eqFactorizationZF}
\end{align}
The following}
procedure is standard and included here for the sake of completeness{; for the multipliers proposed in~\citet{freeman2018noncausal} this step is also discussed in detail in~\citet{zhang2019noncausal}.}}
{Let $ \left[\begin{shortArray}{ccc} M_{-\dimCausal} & \dots & M_{\dimAnticausal} \end{shortArray}\right] \in \mathbb{M}(\rho,\dimAnticausal,\dimCausal,p) $ and define
\begin{subequations}
\ifthenelse{\boolean{singleColumn}}{
\begin{align}
	M_\anticausal = {\left[\begin{shortArray}{cccc} M_1^\top & M_2^\top & \dots & M_{\dimAnticausal}^\top \end{shortArray}\right]}, \qquad 
	M_\causal = \left[\begin{shortArray}{cccc} M_{-\dimCausal} & M_{-\dimCausal+1} & \dots & M_{-1} \end{shortArray}\right],
\end{align}}{
\begin{align}
	M_\anticausal &= {\left[\begin{shortArray}{cccc} M_1^\top & M_2^\top & \dots & M_{\dimAnticausal}^\top \end{shortArray}\right]} \in \real^{p \times \dimAnticausal p}, \\
	M_\causal &= \left[\begin{shortArray}{cccc} M_{-\dimCausal} & M_{-\dimCausal+1} & \dots & M_{-1} \end{shortArray}\right] \in \real^{p \times \dimCausal p},
\end{align}}\label{eqDefCoeffsZF}
\end{subequations}
\ifthenelse{\boolean{singleColumn}}{$ M_\anticausal \in \real^{p \times \dimAnticausal p}, M_\causal \in \real^{p \times \dimCausal p} $, }{}
as well as the vectors of strictly proper and 
stable basis functions $ \tfbasisAnticausal \in \mathcal{RH}_{\infty}^{\dimAnticausal \times 1} $, 
$ \tfbasisCausal \in \mathcal{RH}_{\infty}^{\dimCausal \times 1} $ given by
\begin{subequations}
\ifthenelse{\boolean{singleColumn}}{
\begin{align}
	\tfbasisAnticausal(z) = \begin{bmatrix} z^{-1} & z^{-2} & \dots & z^{-\dimAnticausal} \end{bmatrix}^\top, \qquad
	\tfbasisCausal(z) = \begin{bmatrix} z^{-\dimCausal} & z^{-\dimCausal+1} & \dots & z^{-1}  \end{bmatrix}^\top.
\end{align}}{
\begin{align}
	\tfbasisAnticausal(z) &= \begin{bmatrix} z^{-1} & z^{-2} & \dots & z^{-\dimAnticausal} \end{bmatrix}^\top \label{eqDefBasisZFAnticausal} \\
	\tfbasisCausal(z) &= \begin{bmatrix} z^{-\dimCausal} & z^{-\dimCausal+1} & \dots & z^{-1}  \end{bmatrix}^\top. \label{eqDefBasisZFCausal}
\end{align}}\label{eqDefBasisZF}
\end{subequations}
Using these definitions, it then follows by simple calculations
that $ \Pi \in \mathbf{\Pi}_{\Delta,\rho}^p(m,L) $ is represented as in~\eqref{eqFactorizationZF} if choosing} 
\begin{align}
	M_{\Delta} &=
	\left[
	\begin{array}{cc|cc|cc}
		0        & M_0^\top & 0   & 0        & 0        & 0 \\
		M_0      & 0        & 0   & 0        & 0        & 0 \\ \hline
		0        & 0        & 0   & {M_\causal^\top} & 0        & 0 \\
		0        & 0        & {M_\causal} & 0        & 0        & 0 \\ \hline 
		0        & 0        & 0   & 0        & 0        & {M_\anticausal} \\
		0        & 0        & 0   & 0        & {M_\anticausal^\top} & 0
	\end{array}
	\right], \label{eqDefMDelta} \\
	\psi_{\Delta}(z) &= 
	\left[
	\begin{array}{cc}
		I_p & 0 \\ 0 & I_p \\ \hline
		{\tfbasisCausal}(z) \otimes I_p & 0 \\ 0 & I_p  \\ \hline
		I_p & 0 \\ 0 & {\tfbasisAnticausal}(z) \otimes I_p
	\end{array}
	\right] {\widehat{\trafoSector}}.
	\label{eqDefPsiDelta}
\end{align}
Note that $ \psi_{\Delta} $ is proper {and stable} and, in the following, we denote 
by $ ( A_\Delta, B_\Delta, C_\Delta, D_\Delta ) $ its state-space realization
when using the state-space realizations of $ \psi_+, \psi_- $ provided in~\Cref{secAppendixStateSpaceRealizations}.

\section{Robust optimization algorithms}\label{secRobustOptimizationAlgorithms}
In the following we employ the previously developed tools to 
analyze optimization algorithms as well as design novel algorithms
with prespecified guarantees. 
In particular, we apply~\Cref{lemmaIQCexponentialStability}
utilizing the IQCs derived in~\Cref{subsecIQCsForTheClassOfUncertainties}
and the corresponding multiplier parametrizations from~\Cref{subsecMultiplierParametrization}. 
\subsection{Analysis}
\subsubsection{Convergence rate analysis}
{We begin with deriving convex conditions in terms of linear matrix inequalities (LMIs) for
determining convergence rate bounds for a given algorithm in the form~\eqref{eqOptAlgo}.
In a nutshell, this amount to reformulating the FDI~\eqref{eqFDIRobustExpStab} using the KYP-Lemma\ifthenelse{\boolean{longVersion}}{ (see \Cref{lemmaKYP})}{}.}
Let $ G_{\tilde{y}\tilde{w}} $ be defined according to~\eqref{eqDefTransformedTransferFunction}.
We then {introduce} the following state-space realization
\begin{align}
	\tfcompleteNoPerf = \psi_\Delta \begin{bmatrix} G_{\tilde{y}\tilde{w}} \\ I_p \end{bmatrix}
	\sim 
	( {\AcompleteNoPerf(\rho)}, \BcompleteNoPerf, {\CcompleteNoPerf(\rho)}, \DcompleteNoPerf ),
	\label{eqSSCompleteSystemNoPerformance}
\end{align}
where $ {\AcompleteNoPerf(\rho)} \in \real^{n_\complete \times n_\complete} $, $ \BcompleteNoPerf \in \real^{n_\complete \times p_\complete} $,
$ {\CcompleteNoPerf(\rho)} \in \real^{q_\complete \times n_\complete} $, $ \DcompleteNoPerf \in \real^{q_\complete \times p_\complete} $
with $ n_\complete = p(n+\dimCausal+\dimAnticausal) $, $ p_\complete = p $, $ q_\complete = {p(4+\dimCausal+\dimAnticausal)} $, {are as given in~\Cref{secAppendixStateSpaceRealizations}}. 
The following result then virtually follows immediately from~\Cref{lemmaIQCexponentialStability}
employing the class of multipliers introduced in~\Cref{lemmaZFIQCrho}.
\ifthenelse{\boolean{proofsAppendix}}{A proof is given in~\Cref{secProofLemmaAnalysisConvRate}.}{}
\begin{theorem}\label{lemmaAnalysisConvRate}
Consider \eqref{eqOptAlgo}.
	Suppose $ \Aopt \in \real^{np \times np} $, 
	$ \Bopt \in \real^{np \times p} $, $ \Copt \in \real^{p \times np} $, $ \Dopt \in \real^{np \times p} $,
	$ \Dopt^\dagger \in \real^{np \times p} $ are given and fulfill~\eqref{eqConditionABCD}.
	Let {$ L \geq m > 0 $} be given.
	Fix $ \rho \in (0,1) $ and assume that $ \Aopt + m \Bopt \Copt $ has all eigenvalues in the open disk of radius $ \rho $.
	Let some $ \dimAnticausal, \dimCausal \in \nat $ be given and let $ M_\Delta $, $ \mathbb{M} $ be defined
	according to~\eqref{eqDefMDelta}, \eqref{eqDefSetMunstructured}. If there exist
	$ P = P^\top \in \real^{n_\complete \times n_\complete}, M_+ \in \mathbb {R}^{p \times \dimAnticausal p}, M_- \in \real^{p \times \dimCausal p}, M_0 \in \real^{p \times p} $ such that
	\begin{subequations}
	\begin{align}
		\left[ 
		\vphantom{
		\begin{array}{cc}
			{\AcompleteNoPerf(\rho)} & \BcompleteNoPerf \\
			I                & 0 \\
			{\CcompleteNoPerf(\rho)} & \DcompleteNoPerf
		\end{array}} \star
		\right]^\top
		\left[ 
		\begin{array}{ccc}
			P & 0  & 0 \\
			0 & -P & 0 \\
			0 & 0  & M_\Delta
		\end{array}
		\right] 
		\left[ 
		\begin{array}{cc}
			{\AcompleteNoPerf(\rho)} & \BcompleteNoPerf \\
			I                & 0 \\
			{\CcompleteNoPerf(\rho)} & \DcompleteNoPerf
		\end{array}
		\right] 
		\prec 0,\label{eqLMIConvAnalysis} \\
		\begin{bmatrix}  M_- & M_0 &  M_+ \end{bmatrix} \in \mathbb{M}(\rho,\dimAnticausal,\dimCausal,p), \label{eqMultiplierConstraintConvAnalysis}
	\end{align}\label{eqAllConstraintsConvAnalysis}
	\end{subequations}
	then the origin is robustly exponentially stable against
	$ \mathbf{\Delta}(m,L) $ for~\eqref{eqRobustStandardForm}.
\end{theorem}
\ifthenelse{\boolean{proofsAppendix}}{
}{
\begin{proof}
	The result virtually is a direct consequence of~\Cref{lemmaIQCexponentialStability}
employing the class of multipliers introduced in~\Cref{lemmaZFIQCrho}.
We first note that under the assumption that $ \Aopt + m \Bopt \Copt $ has all
eigenvalues in the open disk of radius $ \rho $, the same holds for $ \AcompleteNoPerf $
by~\eqref{eqSSCompleteSystemNoPerformance} since $ \psi_\Delta $ has all its poles
at zero. With $ \psi_\Delta^\complConj M_\Delta (M_+,M_-,M_0) \psi_\Delta = \Pi $
and by the KYP-Lemma we then observe that~\eqref{eqLMIConvAnalysis} is equivalent
to the frequency domain inequality~\eqref{eqFDIRobustExpStab} for exponential
stability in~\Cref{lemmaIQCexponentialStability}. Further, by~\Cref{lemmaZFIQCrho}
and~\eqref{eqMultiplierConstraintConvAnalysis}, we infer that
$ \IQC\big(\Pi_\rho,\tilde{y},\Delta_{\fObj,\rho}(\tilde{y})\big) \geq 0 $ holds
for all $ \Pi_\rho \in \mathbf{\Pi}_{\Delta,\rho}^p $ and all $ \Delta_{\fObj,\rho} \in \mathbf{\Delta}_\rho(m,L) $.
Hence, by~\Cref{lemmaIQCexponentialStability}, we conclude that the origin is globally robustly exponentially
stable against $ \mathbf{\Delta}(m,L) $ with rate $ \rho $ for~\eqref{eqRobustStandardForm},
which in turn implies the last claim in~\Cref{lemmaAnalysisConvRate}, thus concluding
the proof.

\end{proof}
}
{If the LMIs~\eqref{eqAllConstraintsConvAnalysis} are feasible
for a given optimization algorithm~\eqref{eqOptAlgo} and some $\rho$, then, for all
$ \fObj \in \classObj{m}{L} $, its output 
converges exponentially with rate $ \rho $ to the unique minimizer $ \optVar^\optSign $ of~\eqref{eqOptProb}.}
In convergence rate analysis we are typically interested in finding the smallest $ \rho $
such that~\eqref{eqAllConstraintsConvAnalysis} is feasible. If $ \rho $ 
is kept as a free variable, the inequalities are no longer linear
in the unknown parameters. Still, if~\eqref{eqAllConstraintsConvAnalysis} is feasible for some $ \rho \in (0,1) $,
then it also feasible for all $ \bar{\rho} \in [\rho,1] $. 	
Therefore, it is convenient to do a bisection search
over $ \rho $ to find the optimal convergence rate.

Many optimization algorithms admit the structure
$ \Aopt = \overline{\Aopt} \otimes I_p $, $ \Bopt = \overline{\Bopt} \otimes I_p $,
$ \Copt = \bar{\Copt} \otimes I_p $, $ \Dopt = \overline{\Dopt} \otimes I_p $.
In that case we can also {choose} a state-space realization $ ({\AcompleteNoPerf(\rho)},\BcompleteNoPerf,{\CcompleteNoPerf(\rho)},\DcompleteNoPerf) $
such that $ \AcompleteNoPerf = {\overline{\AcompleteNoPerf}(\rho)} \otimes I_p $,
$ \BcompleteNoPerf = \overline{\BcompleteNoPerf} \otimes I_p $, $ \CcompleteNoPerf = {\overline{\CcompleteNoPerf}(\rho)} \otimes I_p $,
$ \DcompleteNoPerf = \overline{\DcompleteNoPerf} \otimes I_p $. {Recall} that also 
$ M_\Delta(M_+,M_-,M_0) =  M_\Delta(m_+,m_-,m_0) \otimes I_p $.
In this situation we can take $ P = \overline{P} \otimes I_p $ without loss of generality, i.e.,
instead of~\eqref{eqLMIConvAnalysis} we can solve the LMI
\begin{subequations}
\begin{align}
	\left[ 
	\vphantom{
	\begin{array}{cc}
		\AcompleteNoPerf & \BcompleteNoPerf \\
		I                & 0 \\
		\CcompleteNoPerf & \DcompleteNoPerf
	\end{array}} \star
	\right]^\top
	\left[ 
	\begin{array}{ccc}
		{\overline{P}} & 0  & 0 \\
		0 & -{\overline{P}} & 0 \\
		0 & 0  & M_\Delta
	\end{array}
	\right] 
	\left[ 
	\begin{array}{cc}
		{\overline{\AcompleteNoPerf}(\rho)} & \overline{\BcompleteNoPerf} \\
		I                      & 0 \\
		{\overline{\CcompleteNoPerf}(\rho)} & \overline{\DcompleteNoPerf}
	\end{array}
	\right] 
	\prec 0,\label{eqLMIConvAnalysisReduced} \\
	\begin{bmatrix}  m_- & m_0 & m_+ \end{bmatrix} \in \mathbb{M}(\rho,\dimAnticausal,\dimCausal,1),
\end{align}\label{eqAllConstraintsConvAnalysisReduced}
\end{subequations}
where now the dimension is reduced by a factor of $p$ compared to~\eqref{eqLMIConvAnalysis}.
This dimensionality reduction is lossless in the sense that~\eqref{eqAllConstraintsConvAnalysis}
are feasible if and only if~\eqref{eqAllConstraintsConvAnalysisReduced} are feasible.
\ifthenelse{\boolean{longVersion}}{We provide a proof of this statement in~\Cref{secAppendixProofDimensionalityReduction}.}{{A proof of this statement is provided, e.g., in~\citet{mic2019algoDesignArxiv,lessard2016analysis}}.}

\subsubsection{Performance analysis}
{We next derive similar convex conditions for analyzing
the performance of a given algorithm by analogously reformulating
the FDI~\eqref{eqFDIRobustPerformanceExpStab}. While the methodology
applies to a {wide} class of performance specifications 
{\sout{such as $H_\infty$-performance}}, 
we concentrate on $H_2$-performance here and
show that this can be related to noise rejection properties of
the optimization algorithm.} 

We first recall that the $H_2$-norm of a stable linear time-invariant
system with {a} strictly proper transfer matrix $ G $ is defined as
\begin{align}
	\Vert G \Vert_2^2 = \tfrac{1}{2\pi} \textup{tr} \big( { \int\limits_{0}^{2\pi} G(e^{\img \omega})^\complConj G(e^{\img \omega}) \mathrm{d}\omega } \big).
	\label{eqDefH2linearSystem}
\end{align}
This definition of the $H_2$-norm is mainly motivated by the following
two interpretations: first, the $H_2$-norm is a measure for the energy
of the system's impulse response; second, it can as well be interpreted
as the asymptotic variance of the output when the system is driven by white noise.
However, these interpretations do not directly carry over to nonlinear or
time-varying systems~\citep{paganini2000h2}. Different extensions
have been proposed, which are rather based
on the actual desired performance measure instead of a certain system norm.
In the following we build upon a stochastic interpretation; this {will be
motivated in more
detail} in~\Cref{secAnalysisExistingAlgorithms}. We next give a precise
definition of the proposed performance measure.
{
\begin{definition}[{Asymptotic} $H_2$-performance]\label{defAveragedH2Performance}
Consider the state-space representation~\eqref{eqSSgeneralFeedbackInterconnection}
of the feedback interconnection~\eqref{eqFeedbackInterconnection} and
let {$ w_{\perf} = ( W_k )_{k \in \nat} $ be a 
discrete-time white noise process}. Let $ ( y_{\perf,k} )_{k\in \nat} $ denote the corresponding
response of~\eqref{eqSSgeneralFeedbackInterconnection} with initial condition $ x_0 = 0 $.
We then say that the feedback interconnection~\eqref{eqFeedbackInterconnection}
\emphDef{achieves a robust $H_2$-performance level of $ \gamma > 0 $ against $ \mathbf{\Delta} $}
if~\eqref{eqFeedbackInterconnection} is robustly stable against $ \mathbf{\Delta} $ 
and if, for all $ \Delta \in \mathbf{\Delta} $, the performance channel fulfills
{$ \Vert \Gwperftoyperf \Vert_{2,\textup{av}} \leq \gamma $}, where
\begin{align}
	\Vert \Gwperftoyperf \Vert_{2,\textup{av}} =
	\limsup\limits_{k_{\textup{max}} \to\infty} \sqrt{ \tfrac{1}{k_{\textup{max}}} \sum\limits_{k=0}^{k_{\textup{max}}} \expec ( y_{\perf,k}^\top y_{\perf,k}\vphantom{^\top} ) }
\end{align}
and $ \expec(Y) $ {is} the expected value of a random variable $Y$.
\end{definition}
}
{A similar measure has been proposed in~\citet{paganini2000h2}. 
The motivation for this definition is that a small value of this {asymptotic} $H_2$-measure
can be related to the noise rejection properties of the system under consideration.}
If the feedback interconnection~\eqref{eqFeedbackInterconnection} is linear
time invariant, then this definition {specializes} to a standard 
$ H_2 $-norm constraint {involving}~\eqref{eqDefH2linearSystem}.

We next discuss how the so-defined robust $H_2$-performance {properties} can be 
analyzed within the IQC framework.
{We consider} performance multipliers that admit a parametrization of the 
form
\ifthenelse{\boolean{singleColumn}}{
\begin{align}
	\mathbf{\Pi}_\perf = \big\lbrace \psi_\perf^\complConj M_\perf \psi_\perf : M_\perf = \begin{bmatrix} M_{\perf,11} & M_{\perf,12} \\ M_{\perf,12}^\top & M_{\perf,22} \end{bmatrix} \in \mathbf{M}_\perf, 
	M_{\perf,22} \succeq 0 \big\rbrace,
	\label{eqFactorizationPerformanceMultiplier}
\end{align}}{
\begin{align}
	\mathbf{\Pi}_\perf = \big\lbrace \psi_\perf^\complConj M_\perf \psi_\perf : M_\perf = \begin{bmatrix} M_{\perf,11} & M_{\perf,12} \\ M_{\perf,12}^\top & M_{\perf,22} \end{bmatrix} \in \mathbf{M}_\perf, \nonumber \\
	M_{\perf,22} \succeq 0 \big\rbrace,
	\label{eqFactorizationPerformanceMultiplier}
\end{align}}
where $ \mathbf{M}_\perf \subset \real^{ {q \times q} } $ is some convex set
and $ \psi_\perf \in {\mathcal{RH}_{\infty}^{q \times (n_{\wperf} + n_{\yperf})} } $, $ q \in \natPos $, is a proper {and stable} transfer function.
{We emphasize that the most relevant performance measures {fall into that class}, e.g., $H_\infty$- or the $H_2$-performance measure {as considered here}.}
Using the factorization of the Zames-Falb multipliers from~\eqref{eqFactorizationZF}, the FDI for performance~\eqref{eqFDIRobustPerformanceExpStab}
is {then} given by
\begin{align}
	\tfcomplete^\complConj
	\left[ 
	\begin{array}{c|c}
		M_\Delta & 0 \\ \hline 
		0        & M_\perf 
	\end{array}
	\right] 
	\tfcomplete
	\fdi{\prec} 0,
	\label{eqFDIPerformanceParticularProblem}
\end{align}
where
\begin{align}
	\tfcomplete =
	\left[ 
	\begin{array}{c|c}
		\psi_{\Delta} & 0 \\ \hline 
		0             & \psi_\perf
	\end{array}
	\right]
	\left[
	\begin{array}{c|c}
		{\Gwtoy} & {\Gwperftoy} \\
		I_p & 0 \\ \hline
		0 & I_{n_{\wperf}} \\
		{\Gwtoyperf} & {\Gwperftoyperf} 
	\end{array}
	\right].
	\label{eqDefCompleteTransferFunctionFDI}
\end{align}
Let a state-space realization be given by
\begin{align}
	\tfcomplete
	\sim
	\left( 
	\begin{array}{cc}
		\Acomplete                                                           & \begin{bmatrix} \Bcomplete{}_{,1} & \Bcomplete{}_{,2} \end{bmatrix} \\
		\begin{bmatrix} \Ccomplete{}_{,1} \\ \Ccomplete{}_{,2} \end{bmatrix} & \begin{bmatrix} \Dcomplete{}_{,11} & \Dcomplete{}_{,12} \\ \Dcomplete{}_{,21} & \Dcomplete{}_{,22} \end{bmatrix}
	\end{array}
	\right).
	\label{eqSSCompleteSystemWithPerformance}
\end{align}
{Deriving such a state-space realization is standard.}
For implementation purpose, this is particularly easy making 
use of numerical tools such as the control systems toolbox
in \textsc{Matlab}. For completeness, we give an explicit
state-space realization in~\Cref{secAppendixStateSpaceRealizations}.

In the particular case of the considered $H_2$-performace measure, the performance
multiplier {is} defined {with}
\begin{align}
	\psi_{\perf}(z) = \begin{bmatrix} 0_{n_{\yperf} \times n_{\wperf}} & I_{n_{\yperf}} \end{bmatrix}, \quad
	\mathbf{M}_\perf = \big\lbrace M = I_{n_{\yperf}} \big\rbrace
	\label{eqDefH2IQC}
\end{align}
and {we} denote the resulting state-space realization by~\eqref{eqSSCompleteSystemWithPerformance}.
{\sout{In contrast to $H_\infty$-performance,
$H_2$-performance cannot be directly captured by means of a performance IQC,
but an additional condition based on state-space representations is required.}}
Note that $ \mathbf{M}_\perf $ 
is a singleton and the $H_2$-performance level is set by
the aforementioned additional condition {resulting in}~\eqref{eqTraceH2Analysis}.
We then have the following result for $H_2$-performance\ifthenelse{\boolean{proofsAppendix}}{; a proof is given in~\Cref{secProofLemmaAnalysisH2}.}{.}
\begin{theorem}[Robust $H_2$-performance]\label{lemmaAnalysisH2}
Consider \eqref{eqRobustStandardForm} and suppose $ \Dperf = 0 $.
	Suppose $ \Aopt \in \real^{np \times np} $, 
	$ \Bopt \in \real^{np \times p} $, $ \Copt \in \real^{p \times np} $,
	$ \Dopt \in \real^{p \times np} $, $ \Dopt^\dagger \in \real^{np \times p} $ are given and fulfill~\eqref{eqConditionABCD}.
	Let {$ L \geq m > 0 $} be given.
	Assume that \mbox{$ \Aopt + m \Bopt \Copt $} has all eigenvalues in the open unit disk.
	Let some $ \dimAnticausal, \dimCausal \in \nat $ be given and let $ M_\Delta $ be defined
	according to~\eqref{eqDefMDelta}. If there exist $ P_\perf = P_\perf^\top {\in \real^{n_\complete \times n_\complete}} $,
	$ M_+ \in \mathbb {R}^{p \times p\dimAnticausal} $, $M_- \in \real^{p \times p\dimCausal}, M_0 \in \real^{p \times p} $, $ \gamma > 0 $,
	such that
	\begin{subequations}
	\begin{align}
		\left[
		\vphantom{
		\begin{array}{ccc}
			\Acomplete        & \Bcomplete{}_{,1}  & \Bcomplete{}_{,2} \\
			I                 & 0                  & 0 \\ \hline 
			\Ccomplete{}_{,1} & \Dcomplete{}_{,11} & 0 \\ 
			\Ccomplete{}_{,2} & 0                  & 0
		\end{array}} \star
		\right]^\top
		\left[
		\begin{shortArray}{cc|cc}
			P_\perf & 0         & 0        & 0 \\
			0       & - P_\perf & 0        & 0 \\ \hline
			0       & 0         & M_\Delta & 0 \\
			0       & 0         & 0        & I_{n_{\yperf}} 
		\end{shortArray}
		\right]
		\left[
		\begin{shortArray}{ccc}
			\Acomplete        & \Bcomplete{}_{,1}   \\
			I                 & 0                   \\ \hline 
			\Ccomplete{}_{,1} & \Dcomplete{}_{,11}  \\ 
			\Ccomplete{}_{,2} & 0                  
		\end{shortArray}
		\right] \prec 0, & \label{eqLMIH2Analysis} \\
		\textup{tr}( \Bcomplete{}_{,2}^\top N N^\top P_\perf N N^\top \Bcomplete{}_{,2} ) \leq \gamma^2, & \label{eqTraceH2Analysis} \\
		{ P_\perf } \succ 0, & \label{eqPosDefH2Analysis} \\
		N^\top = {\begin{bmatrix} 0_{np \times {(\dimCausal+\dimAnticausal)p}} & I_{np} \end{bmatrix}}, & \label{eqDefN} \\
		\begin{bmatrix}  M_- & M_0 &  M_+ \end{bmatrix} \in \mathbb{M}({1},\dimAnticausal,\dimCausal,p), & 
	\end{align}\label{eqLMIPerfAnalysis}
	\end{subequations}
	then{\sout{, for all $ \fObj \in {\classObj{m}{L}} $,}} the interconnection~\eqref{eqRobustStandardForm}
	is robustly stable {against $ \deltaUnstr(m,L) $} and it achieves a robust $ H_2 $-performance level {of $ \gamma $ against $ \deltaUnstr(m,L) $}.
\end{theorem}
\ifthenelse{\boolean{proofsAppendix}}{
}{
\begin{proof}
	\input{proofH2.tex}
\end{proof}
}
\begin{remark}\label{remarkReformulationTrace}
    {Under~\eqref{eqPosDefH2Analysis}, it is well-known that}~\eqref{eqTraceH2Analysis} holds
	if {and only if} there exists a matrix $Z$ of suitable {dimensions} such that
	\begin{align}
		\begin{bmatrix} N^\top P_\perf N & N^\top P_\perf N N^\top \Bcomplete{}_{,2} \\ \star & Z \end{bmatrix} \succ 0, \quad 
		\textup{tr}(Z) < \gamma^2.
	\end{align}
	\ifthenelse{\boolean{longVersion}}{To see this, note that, by taking Schur complements, the above
	conditions hold if and only if 
	\begin{align}
		N^\top P_\perf N \succ 0, \quad Z \succ \Bcomplete{}_{,2}^\top N N^\top P_\perf N N^\top \Bcomplete{}_{,2}  
	\end{align}
	and $ \textup{tr}(Z) < \gamma^2 $.
	If the latter inequalities are feasible, this implies that there exists an $ \varepsilon > 0 $ such that
	$ Z = \Bcomplete{}_{,2}^\top N N^\top P_\perf N N^\top \Bcomplete{}_{,2} + \varepsilon I $
	is a valid solution. Then this implies that~\eqref{eqTraceH2Analysis} holds.}{}
\end{remark}
{If the LMIs~\eqref{eqLMIPerfAnalysis} are feasible
for a given optimization algorithm~\eqref{eqTransformedOptAlgo} and some $ \gamma $,
then, for all $ \fObj \in \classObj{m}{L} $, its performance
channel has an $H_2$-performance level of $ \gamma $.}
{Similar as} for the convergence
rate test from~\Cref{lemmaAnalysisConvRate}, an analogous lossless dimensionality
reduction is possible for $H_2$-performance analysis
for equally structured optimization algorithms whenever the matrices 
defining the performance channel admit the same structure, i.e.,
$ \Bperf = \overline{\Bperf} \otimes I_p $, $ \Cperf = \overline{\Cperf} \otimes I_p $, $ \Dperf = \overline{\Dperf} \otimes I_p $.
This structure then also transfers to the matrices $ \Acomplete, \Bcomplete, \Ccomplete, \Dcomplete $.

\subsubsection{Numerical results}\label{secAnalysisExistingAlgorithms}
In the following we use the presented results to numerically analyze existing
optimization algorithms in terms of their convergence rate as well as
their properties in the presence of additive gradient noise.
Here, we consider the Gradient Descent algorithm (GD), Nesterov's Method (NM) with constant step
size~\citep{nesterov2004introductory}, the Triple Momentum Method (TMM)~\citep{vanScoy2018fastest} {and}
the Heavy Ball Method (HB)~\citep{polyak1987introduction}.
In all cases, the matrices $ \Aopt, \Bopt, \Copt, \Dopt $ in~\eqref{eqOptAlgo}
take the form
\begin{subequations}
\begin{align}
	\Aopt &= \begin{bmatrix} 1+\nu_2 & -\nu_2 \\ 1 & 0 \end{bmatrix} \otimes I_p, &
	\Bopt &= \begin{bmatrix} - \nu_1 \\ 0 \end{bmatrix} \otimes I_p, \\ 
	\Copt &= \begin{bmatrix} 1+\nu_3 & -\nu_3 \end{bmatrix} \otimes I_p, &
	\Dopt &= \begin{bmatrix} 1 & 0 \end{bmatrix} \otimes I_p,
\end{align}\label{eqDefExistingAlgorithms}
\end{subequations}
where the scalar parameters $ \nu_1, \nu_2, \nu_3 $ are as given 
in \Cref{tableParametersAlgorithms}. Note that with $ \Dopt^\dagger = \begin{bmatrix} 1 & 1 \end{bmatrix}^\top \otimes I_p $,
the conditions~\eqref{eqConditionABCD} are fulfilled.
Due to the Kronecker structure and by our previous discussion, 
we set $p=1$ without loss of generality.
{The following {and} all other numerical results in the paper were obtained 
using \textsc{Matlab} together with YALMIP \citep{lofberg2004yalmip}.}
\ifthenelse{\boolean{singleColumn}}{
\begin{table*}[ht!]
	\begin{center}
	\begin{tabular}{ccccc}
		\toprule
		& \multicolumn{3}{c}{\textbf{Parameters}}                                                                                        & \multicolumn{1}{c}{\textbf{Convergence rate}} \\
		& $ \nu_1 $                            & $ \nu_2 $                                        & $ \nu_3 $                            & $ \rho $ \\ \midrule
		$ \begin{matrix} \text{\textbf{Gradient}} \\ \text{\textbf{Descent}} \end{matrix} $       & $ \tfrac{2}{m+L} $                   & $ 0 $                                            & $ 0 $                                & $ \tfrac{\kappa-1}{\kappa+1} $ \\[0.8em]
		$ \begin{matrix} \text{\textbf{Nesterov's}} \\ \text{\textbf{Method}} \end{matrix} $      & $ \tfrac{1}{L} $                     & $ \tfrac{\sqrt{L}-\sqrt{m}}{\sqrt{L}+\sqrt{m}} $ & $ \nu_2 $                            & {$ \sqrt{ 1 - \tfrac{1}{\sqrt{\kappa}} } $} 
		\\[0.8em]		
		{$ \begin{matrix} \text{\textbf{Nesterov's}} \\ \text{\textbf{Method (modified)}} \end{matrix} $}       & {$ \tfrac{1}{L} $}                     & {$ \tfrac{2\kappa - \sqrt{2\kappa-1}-1}{2(\kappa+\sqrt{\kappa-1})} $} & {$ \nu_2 $}                            & {$ \sqrt{ 1- \tfrac{1}{\kappa} \sqrt{2\kappa-1} } $} \\[0.8em]
		$ \begin{matrix} \text{\textbf{Triple Momentum}} \\ \text{\textbf{Method}} \end{matrix} $ & $ \tfrac{1+\rho}{L} $                & $ \tfrac{\rho^2}{2-\rho} $                       & $ \tfrac{\rho^2}{(1+\rho)(2-\rho)} $ & $ 1 - \tfrac{1}{\sqrt{\kappa}} $ \\[0.8em] 
		$ \begin{matrix} \text{\textbf{Heavy Ball}} \\ \text{\textbf{Method}} \end{matrix} $      & $ (\tfrac{2}{\sqrt{L}+\sqrt{m}})^2 $ & $ {(\tfrac{\sqrt{L}-\sqrt{m}}{\sqrt{L}+\sqrt{m}})^2} $ & $ 0 $                                & -- \\ \bottomrule
	\end{tabular}
	\end{center}
	\caption{Parameters of several popular algorithms of the form~\eqref{eqDefExistingAlgorithms}
	and the corresponding known upper bounds on the convergence rate in dependency of the condition ratio $ \kappa = \tfrac{L}{m} $.
	For details, see~\citet{polyak1987introduction} in case of the Gradient Descent, {\citet{nesterov2004introductory} in case of Nesterov's Method and~\citet{safavi2018explicit} in case of the modified version thereof}, \citet{vanScoy2018fastest} for the Triple Momentum Method and
	\citet{polyak1987introduction} for the Heavy Ball Method.	
	Note that {for}
	Gradient Descent algorithm the parameter choice reduces~\eqref{eqDefExistingAlgorithms} to a
	first order algorithm. While an explicit convergence rate can be given for the Heavy Ball Method
	in case of quadratic objective functions, it is known not to be globally convergent for general
	$ \fObj \in {\classObj{m}{L}} $, see~\cite{lessard2016analysis} for a counterexample. }\label{tableParametersAlgorithms}
\end{table*}}{
\begin{table*}[ht!]
	\begin{minipage}[T]{0.634\textwidth}
	\begin{tabular}{cccccc}
		\toprule
		& \multicolumn{3}{c}{\textbf{Parameters}}                                                                                        & \multicolumn{2}{c}{\textbf{Convergence rate}} \\
		& $ \nu_1 $                            & $ \nu_2 $                                        & $ \nu_3 $                            & $ \rho $ \\ \midrule \\[-0.5em]
		\textbf{GD}       & $ \tfrac{2}{m+L} $                   & $ 0 $                                            & $ 0 $                                & $ \tfrac{\kappa-1}{\kappa+1} $ & \citet{polyak1987introduction} \\[0.8em]
		\textbf{NM}      & $ \tfrac{1}{L} $                     & $ \tfrac{\sqrt{L}-\sqrt{m}}{\sqrt{L}+\sqrt{m}} $ & $ \nu_2 $                            & $ \sqrt{ 1- \tfrac{1}{\kappa} \sqrt{2\kappa-1} } $  
		& \citet{safavi2018explicit} \\[0.8em]		
		\textbf{TMM} & $ \tfrac{1+\rho}{L} $                & $ \tfrac{\rho^2}{2-\rho} $                       & $ \tfrac{\rho^2}{(1+\rho)(2-\rho)} $ & $ 1 - \tfrac{1}{\sqrt{\kappa}} $ & \citet{vanScoy2018fastest} \\[0.8em] 
		\textbf{HB}  & $ (\tfrac{2}{\sqrt{L}+\sqrt{m}})^2 $ & $ \tfrac{\sqrt{L}-\sqrt{m}}{\sqrt{L}+\sqrt{m}} $ & $ 0 $                                & -- & \citet{polyak1987introduction} \\ \bottomrule
	\end{tabular}
	\end{minipage}
	\begin{minipage}[T]{0.36\textwidth}
		\caption{Parameters of several popular algorithms of the form~\eqref{eqDefExistingAlgorithms}
		and the corresponding known upper bounds on the convergence rate in dependency of the condition ratio $ \kappa = \tfrac{L}{m} $. Note that in the case of the 
		Gradient Descent algorithm the parameter choice reduces~\eqref{eqDefExistingAlgorithms} to a
		first order algorithm. While an explicit convergence rate can be given for the Heavy Ball Method
		in case of quadratic objective functions, it is known not to be globally convergent for general
		$ \fObj \in {\classObj{m}{L}} $, see~\citet{lessard2016analysis} for a counterexample. }\label{tableParametersAlgorithms}
	\end{minipage}
\end{table*}}
\paragraph{Convergence rate analysis}
For convergence rate analysis, we use~\Cref{lemmaAnalysisConvRate} together with a bisection search over
$ \rho $ to determine upper bounds on the convergence rates for different 
condition ratios $ \kappa = L/m $, see~\Cref{figConvRatesAnalysis}. For Nesterov's
Method, the Gradient Descent algorithm and the Triple Momentum Method, we reproduce
the known convergence rate bounds. For the Heavy Ball Method, global 
convergence can be guaranteed for small condition numbers; for larger
condition numbers the LMIs in~\Cref{lemmaAnalysisConvRate} turn out
to be infeasible, even for $ \rho = 1 $. These results are in concordance
with~\citet{lessard2016analysis}. 
{For all considered algorithms, }
the addition of anticausal
Zames-Falb multipliers does not lead to any improvement of the upper bound
on the convergence rate and it is even sufficient to choose {$ \dimAnticausal = 0, \dimCausal = 1 $}.
{\sout{In {contrast}, for the Heavy Ball Method employing anticausal multipliers
with $ \dimAnticausal = 1, \dimCausal = 5 $ leads to slightly better convergence rate
estimates. Using higher dimensions for both the causal and the anticausal
part does not significantly improve the bounds.}}
Except for the Gradient
Descent algorithm, it is still an open 
question whether the known as well as the numerically determined
bounds are tight in the sense that there exists some objective
function $ \fObj \in {\classObj{m}{L}} $ such that the resulting
algorithm does not converge faster than presumed by the given bound.
\newlength\figurewidth
\newlength\figureheight
\ifthenelse{\boolean{singleColumn}}{
\setlength\figurewidth{0.5\textwidth}
\setlength\figureheight{0.5\textwidth}
}{
\setlength\figurewidth{0.85\columnwidth}
\setlength\figureheight{0.7\columnwidth}}
\begin{figure}[t]
	\begin{center}
		\input{../figs/analysisConvRateExistingAlgos.tex}
	\end{center}
	\caption{Upper bounds on the convergence rates of different algorithms 
	obtained from~\Cref{lemmaAnalysisConvRate} for condition ratios between $1.02$ and $1000$.
	The corresponding dashed lines indicate the known {analytical} convergence rate bounds,
	cf.~\Cref{tableParametersAlgorithms}{, where, for Nesterov's Method, 
	the analytical bound for the modified version is depicted}. The black line indicates
	the fundamental lower bound for any first order optimization algorithm
	for objective functions in $ {\classObj{m}{L}} $,
	see~\citet{nesterov2004introductory}.}
	\label{figConvRatesAnalysis}
\end{figure}

\paragraph{$H_2$-performance} In the following we investigate the properties
of the algorithms from~\Cref{tableParametersAlgorithms} when the
gradient is affected by white noise in an additive fashion, i.e., in rough words, in~\eqref{eqOptAlgo}
we have $ \nabla \fObj(\Copt \xOpt\k) + W\k $ instead of $ \nabla \fObj(\Copt \xOpt\k) $,
where $ ( W_k )_{k \in \nat} $ is a {discrete-time white noise process}. Such situations occur in several applications, 
e.g., when the objective function and its gradient are not evaluated by means of 
numeric calculations but rather by measurements~\citep{romer2017sampling}.
Likewise, in empirical risk minimization as it is utilized, e.g., in the context of 
learning algorithms~\citep{murphy2012machine}, the objective function is given by an expected value,
which, however, cannot be evaluated since the underlying probability 
distribution is unknown. It is therefore common to use a sample-based
approximation of the expected value, the so-called empirical risk.
According to the central limit theorem, this approximation differs from
the original expected value by additive random noise, see~\citet{polyak1987introduction}.
Similar situations appear when employing Monte-Carlo methods.
The $H_2$-performance as defined in~\Cref{defAveragedH2Performance} 
is a measure for the noise attenuation; 
hence we choose the performance channel as $ \Bperf = \Bopt $,
$ \Cperf = \Dopt $, $ \Dperf = 0 $ such that the corresponding
$ H_2 $-performance level is a measure how additive gradient noise
affects the resulting optimizer. We then use~\Cref{lemmaAnalysisH2} to 
determine upper bounds on the corresponding $H_2$-performance.
The results are depicted \ifthenelse{\boolean{singleColumn}}{on the left-hand side of~\Cref{figH2Analysis}}{in~\Cref{figH2Analysis}}.
\ifthenelse{\boolean{singleColumn}}{
\setlength\figurewidth{0.39\textwidth}
\setlength\figureheight{0.39\textwidth}
\begin{figure}[ht!]
	\begin{center}
		\begin{footnotesize}
%
%
\definecolor{mycolor1}{rgb}{0.00000,0.44700,0.74100}%
\definecolor{mycolor2}{rgb}{0.85000,0.32500,0.09800}%
\definecolor{mycolor3}{rgb}{0.92900,0.69400,0.12500}%
\definecolor{mycolor4}{rgb}{0.49400,0.18400,0.55600}%
\begin{tikzpicture}

\begin{axis}[%
width=0.951\figurewidth,
height=\figureheight,
at={(0\figurewidth,0\figureheight)},
scale only axis,
unbounded coords=jump,
xmode=log,
xmin=1,
xmax=1000,
xminorticks=true,
xlabel={Condition ratio $L/m$},
ymode=log,
ymin=0,
ymax=1.25,
yminorticks=true,
ylabel={$H_2$-performance},
axis background/.style={fill=white},
legend style={at={(0.03,0.03)},anchor=south west,legend cell align=left,align=left,draw=white!15!black,font=\footnotesize}
]
\addplot [color=mycolor1,mark size=2.5pt,only marks,mark=*,mark options={solid,fill=mycolor1}]
  table[row sep=crcr]{%
1.02329299228075	0.977493376291846\\
1.29748829388885	0.794133790294052\\
1.64515528346031	0.670678852888087\\
2.08598098298468	0.582215773499571\\
2.64492762787805	0.515463220425393\\
3.35364618075428	0.462898801369478\\
4.25226860165957	0.420044228240986\\
5.39168036402487	0.384116188593178\\
6.83640190002716	0.353316754485212\\
8.66824214034203	0.326444559173021\\
10.9909310339556	0.302682904709018\\
13.9359933695162	0.282472928578736\\
17.6701965097587	0.26552673926366\\
22.4049937750743	0.250843495285018\\
28.4084982181091	0.23774009171415\\
36.020664817419	0.225757193748679\\
45.6725407984347	0.214597575147028\\
57.9106741521309	0.204069026093036\\
73.4280625103569	0.194052381949024\\
93.1033948916044	0.184470543087233\\
118.050808423815	0.175281911465653\\
149.68297757285	0.166459882680632\\
189.791108372913	0.157989896948332\\
240.646367419354	0.149861855359282\\
305.128489151053	0.142069570399091\\
386.888843949848	0.134607601335331\\
490.557200965753	0.127471748449694\\
622.003893838172	0.120656101188857\\
788.672234732637	0.114154177697952\\
1000	0.107960106602726\\
};
\addlegendentry{Nesterov's Method};

\addplot [color=mycolor2,mark size=2.5pt,only marks,mark=*,mark options={solid,fill=mycolor2}]
  table[row sep=crcr]{%
1.02329299228075	0.988489821423798\\
1.29748829388885	0.873022871188138\\
1.64515528346031	0.769569100153459\\
2.08598098298468	0.681329302949565\\
2.64492762787805	0.60756934437046\\
3.35364618075428	0.546134249262287\\
4.25226860165957	0.49635445769231\\
5.39168036402487	0.458065320975845\\
6.83640190002716	0.429466565712428\\
8.66824214034203	0.408607542616401\\
10.9909310339556	0.39357315741378\\
13.9359933695162	0.382640276098979\\
17.6701965097587	0.374372168487981\\
22.4049937750743	0.3677149626777\\
28.4084982181091	0.362056526617709\\
36.020664817419	0.356937604717942\\
45.6725407984347	0.352068741781101\\
57.9106741521309	0.347280871032798\\
73.4280625103569	0.342490069814027\\
93.1033948916044	0.3376670474672\\
118.050808423815	0.332817375430929\\
149.68297757285	0.327959034042448\\
189.791108372913	0.323129019173182\\
240.646367419354	0.318358723728575\\
305.128489151053	0.313683092209288\\
386.888843949848	0.309131117965972\\
490.557200965753	0.304734214053568\\
622.003893838172	0.300515186892446\\
788.672234732637	0.296487948212373\\
1000	0.292662483580413\\
};
\addlegendentry{Triple Momentum Method};

\addplot [color=mycolor3,mark size=2.5pt,only marks,mark=*,mark options={solid,fill=mycolor3}]
  table[row sep=crcr]{%
1.02329299228075	0.988553097160194\\
1.29748829388885	0.877906522889538\\
1.64515528346031	0.779644382409624\\
2.08598098298468	0.692380499306918\\
2.64492762787805	0.614883876641706\\
3.35364618075428	0.546061275724849\\
4.25226860165957	0.484941839571305\\
5.39168036402487	0.430663370544359\\
6.83640190002716	0.382460164168735\\
8.66824214034203	0.339652040855451\\
10.9909310339556	0.301635690247641\\
13.9359933695162	0.267874290818168\\
17.6701965097587	0.237891710972655\\
22.4049937750743	0.211264871788207\\
28.4084982181091	0.187618589115092\\
36.020664817419	0.166618543326396\\
45.6725407984347	0.147969573265824\\
57.9106741521309	0.131407663411517\\
73.4280625103569	0.116697670757833\\
93.1033948916044	0.103637296297367\\
118.050808423815	0.0920376619090705\\
149.68297757285	0.0817234332700407\\
189.791108372913	0.0725822172100945\\
240.646367419354	0.0644493432874639\\
305.128489151053	0.0572478164339619\\
386.888843949848	0.0508064444669131\\
490.557200965753	0.0451498143095696\\
622.003893838172	0.0400797822030554\\
788.672234732637	0.0356084634945976\\
1000	0.0316230835082188\\
};
\addlegendentry{Gradient Descent};

\addplot [color=mycolor4,mark size=2.5pt,only marks,mark=*,mark options={solid,fill=mycolor4}]
  table[row sep=crcr]{%
1.02329299228075	0.988586148657721\\
1.29748829388885	0.882027882436043\\
1.64515528346031	0.794370124366101\\
2.08598098298468	0.723883803057016\\
2.64492762787805	0.669137931908033\\
3.35364618075428	0.629037877388318\\
4.25226860165957	0.602996067034194\\
5.39168036402487	0.591353739445011\\
6.83640190002716	0.596342410163422\\
8.66824214034203	0.624585395768024\\
10.9909310339556	0.695623382015381\\
13.9359933695162	0.888832387391917\\
17.6701965097587	3.31327873506929\\
22.4049937750743	inf\\
28.4084982181091	inf\\
36.020664817419	inf\\
45.6725407984347	inf\\
57.9106741521309	inf\\
73.4280625103569	inf\\
93.1033948916044	inf\\
118.050808423815	inf\\
149.68297757285	inf\\
189.791108372913	inf\\
240.646367419354	inf\\
305.128489151053	inf\\
386.888843949848	inf\\
490.557200965753	inf\\
622.003893838172	inf\\
788.672234732637	inf\\
1000	inf\\
};
\addlegendentry{Heavy Ball Method};

\end{axis}
\end{tikzpicture}%
%
%
\definecolor{mycolor1}{rgb}{0.00000,0.44700,0.74100}%
\definecolor{mycolor2}{rgb}{0.85000,0.32500,0.09800}%
\definecolor{mycolor3}{rgb}{0.92900,0.69400,0.12500}%
\definecolor{mycolor4}{rgb}{0.49400,0.18400,0.55600}%
\begin{tikzpicture}

\begin{axis}[%
width=0.951\figurewidth,
height=\figureheight,
at={(0\figurewidth,0\figureheight)},
scale only axis,
xmode=log,
xmin=1,
xmax=1000,
xminorticks=true,
xlabel={Condition ratio $L/m$},
ylabel={$H_2$-performance},
ymode=log,
ymin=1e-05,
ymax=1,
yminorticks=true,
axis background/.style={fill=white},
legend style={at={(0.03,0.03)},anchor=south west,legend cell align=left,align=left,draw=white!15!black}
]
\addplot [color=mycolor1,only marks,mark=*,mark options={solid,fill=mycolor1}]
  table[row sep=crcr]{%
1.02329299228075	0.954744019462147\\
1.29748829388885	0.605963216621298\\
1.64515528346031	0.3946619248306\\
2.08598098298468	0.261569720237658\\
2.64492762787805	0.176790055445649\\
3.35364618075428	0.120370896922846\\
4.25226860165957	0.0831860543200114\\
5.39168036402487	0.0572941093581222\\
6.83640190002716	0.0402043961725882\\
8.66824214034203	0.0283710141308075\\
10.9909310339556	0.0198249493011395\\
13.9359933695162	0.0139468298904699\\
17.6701965097587	0.0100053357145601\\
22.4049937750743	0.00713664706083574\\
28.4084982181091	0.00510590471086536\\
36.020664817419	0.00361217041373601\\
45.6725407984347	0.00258858782348511\\
57.9106741521309	0.00183999848180137\\
73.4280625103569	0.00133648377699781\\
93.1033948916044	0.000963970386103095\\
118.050808423815	0.000684323408333818\\
149.68297757285	0.000494231586023779\\
189.791108372913	0.000369488802415231\\
240.646367419354	0.000255497649440165\\
305.128489151053	0.000164891298346631\\
386.888843949848	0.000124431124484346\\
490.557200965753	8.87827152921937e-05\\
622.003893838172	6.82974237560819e-05\\
788.672234732637	4.50411826941585e-05\\
1000	3.5107314510134e-05\\
};
\addlegendentry{Nesterov's Method};

\addplot [color=mycolor2,only marks,mark=*,mark options={solid,fill=mycolor2}]
  table[row sep=crcr]{%
1.02329299228075	0.976600196935738\\
1.29748829388885	0.75229166988042\\
1.64515528346031	0.561675959847349\\
2.08598098298468	0.410650394645257\\
2.64492762787805	0.296731874116579\\
3.35364618075428	0.21245461519408\\
4.25226860165957	0.15110094639101\\
5.39168036402487	0.106525649488244\\
6.83640190002716	0.0753926944427708\\
8.66824214034203	0.0532905127155199\\
10.9909310339556	0.0373002030697877\\
13.9359933695162	0.0261703224443414\\
17.6701965097587	0.0185800343951864\\
22.4049937750743	0.0130668244401583\\
28.4084982181091	0.00922581936298606\\
36.020664817419	0.00645694741836424\\
45.6725407984347	0.00458237249591113\\
57.9106741521309	0.0032159251505746\\
73.4280625103569	0.00230631100462642\\
93.1033948916044	0.00162975286712203\\
118.050808423815	0.00115202458517901\\
149.68297757285	0.000818401890260592\\
189.791108372913	0.000597459915405302\\
240.646367419354	0.000413430718217041\\
305.128489151053	0.000268234092877216\\
386.888843949848	0.000198891172183039\\
490.557200965753	0.000140538765085174\\
622.003893838172	0.000106129847815099\\
788.672234732637	7.03768245292143e-05\\
1000	5.37505814829613e-05\\
};
\addlegendentry{Triple Momentum Method};

\addplot [color=mycolor3,only marks,mark=*,mark options={solid,fill=mycolor3}]
  table[row sep=crcr]{%
1.02329299228075	0.976726710738169\\
1.29748829388885	0.762207753849123\\
1.64515528346031	0.583245617406832\\
2.08598098298468	0.438142540113632\\
2.64492762787805	0.32403663427265\\
3.35364618075428	0.236278565467668\\
4.25226860165957	0.169390206897533\\
5.39168036402487	0.11970533633624\\
6.83640190002716	0.0839126605482172\\
8.66824214034203	0.0582305120134707\\
10.9909310339556	0.0398143969450539\\
13.9359933695162	0.0271286206232095\\
17.6701965097587	0.0184917079243119\\
22.4049937750743	0.0123292165381882\\
28.4084982181091	0.0082609174082255\\
36.020664817419	0.00544769456296845\\
45.6725407984347	0.00367620849605631\\
57.9106741521309	0.00241358580183607\\
73.4280625103569	0.00163450020562211\\
93.1033948916044	0.00104066398342944\\
118.050808423815	0.000706796321464197\\
149.68297757285	0.0004484204665634\\
189.791108372913	0.00029314573126722\\
240.646367419354	0.000189973175774232\\
305.128489151053	0.000118853761885499\\
386.888843949848	7.86953596965827e-05\\
490.557200965753	4.96105161379406e-05\\
622.003893838172	3.33912232106491e-05\\
788.672234732637	2.10135576872057e-05\\
1000	1.35685074148269e-05\\
};
\addlegendentry{Gradient Descent};

\addplot [color=mycolor4,only marks,mark=*,mark options={solid,fill=mycolor4}]
  table[row sep=crcr]{%
1.02329299228075	0.97679143209998\\
1.29748829388885	0.768679026421152\\
1.64515528346031	0.601405621507534\\
2.08598098298468	0.468088194480475\\
2.64492762787805	0.363121467575609\\
3.35364618075428	0.280862449078522\\
4.25226860165957	0.215712126421549\\
5.39168036402487	0.164740539879127\\
6.83640190002716	0.125739297002123\\
8.66824214034203	0.0956007410449833\\
10.9909310339556	0.0720237769567625\\
13.9359933695162	0.0542820838294688\\
17.6701965097587	0.0410784168675066\\
22.4049937750743	0.0304932522217804\\
28.4084982181091	0.0227905543670198\\
36.020664817419	0.0168053707125002\\
45.6725407984347	0.0126908181100361\\
57.9106741521309	0.00935165873092465\\
73.4280625103569	0.00709021074723033\\
93.1033948916044	0.00506785629756977\\
118.050808423815	0.00387237120880326\\
149.68297757285	0.00276509986546831\\
189.791108372913	0.00202163005295546\\
240.646367419354	0.00148059473150475\\
305.128489151053	0.00103979830988667\\
386.888843949848	0.000773421935630782\\
490.557200965753	0.000552241117167733\\
622.003893838172	0.000424561208453707\\
788.672234732637	0.000295110550742672\\
1000	0.000215639226629231\\
};
\addlegendentry{Heavy Ball Method};

\end{axis}
\end{tikzpicture}%
		\end{footnotesize}
	\end{center}
	\captionof{figure}{Upper bounds on the $H_2$-performance related to additive gradient noise
	for different optimization algorithms and condition ratios between $1.02$ and $1000$
	obtained using~\Cref{lemmaAnalysisH2} (left) and {lower bounds from} a sample-based approach (right).
	In all cases, the dimensions of the Zames-Falb multipliers were chosen as {$ \dimAnticausal = 0 $,
	$ \dimCausal = 4 $}.}
	\label{figH2Analysis}
 \end{figure}
}{
\setlength\figurewidth{0.85\columnwidth}
\setlength\figureheight{0.7\columnwidth}
\begin{figure}[ht!]
	\begin{center}
%
%
\definecolor{mycolor1}{rgb}{0.00000,0.44700,0.74100}%
\definecolor{mycolor2}{rgb}{0.85000,0.32500,0.09800}%
\definecolor{mycolor3}{rgb}{0.92900,0.69400,0.12500}%
\definecolor{mycolor4}{rgb}{0.49400,0.18400,0.55600}%
\begin{tikzpicture}

\begin{axis}[%
width=0.951\figurewidth,
height=\figureheight,
at={(0\figurewidth,0\figureheight)},
scale only axis,
unbounded coords=jump,
xmode=log,
xmin=1,
xmax=1000,
xminorticks=true,
xlabel={Condition ratio $L/m$},
ymode=log,
ymin=0,
ymax=1.25,
yminorticks=true,
ylabel={$H_2$-performance},
axis background/.style={fill=white},
legend style={at={(0.03,0.03)},anchor=south west,legend cell align=left,align=left,draw=white!15!black,font=\footnotesize}
]
\addplot [color=mycolor1,mark size=2.5pt,only marks,mark=*,mark options={solid,fill=mycolor1}]
  table[row sep=crcr]{%
1.02329299228075	0.977493376291846\\
1.29748829388885	0.794133790294052\\
1.64515528346031	0.670678852888087\\
2.08598098298468	0.582215773499571\\
2.64492762787805	0.515463220425393\\
3.35364618075428	0.462898801369478\\
4.25226860165957	0.420044228240986\\
5.39168036402487	0.384116188593178\\
6.83640190002716	0.353316754485212\\
8.66824214034203	0.326444559173021\\
10.9909310339556	0.302682904709018\\
13.9359933695162	0.282472928578736\\
17.6701965097587	0.26552673926366\\
22.4049937750743	0.250843495285018\\
28.4084982181091	0.23774009171415\\
36.020664817419	0.225757193748679\\
45.6725407984347	0.214597575147028\\
57.9106741521309	0.204069026093036\\
73.4280625103569	0.194052381949024\\
93.1033948916044	0.184470543087233\\
118.050808423815	0.175281911465653\\
149.68297757285	0.166459882680632\\
189.791108372913	0.157989896948332\\
240.646367419354	0.149861855359282\\
305.128489151053	0.142069570399091\\
386.888843949848	0.134607601335331\\
490.557200965753	0.127471748449694\\
622.003893838172	0.120656101188857\\
788.672234732637	0.114154177697952\\
1000	0.107960106602726\\
};
\addlegendentry{Nesterov's Method};

\addplot [color=mycolor2,mark size=2.5pt,only marks,mark=*,mark options={solid,fill=mycolor2}]
  table[row sep=crcr]{%
1.02329299228075	0.988489821423798\\
1.29748829388885	0.873022871188138\\
1.64515528346031	0.769569100153459\\
2.08598098298468	0.681329302949565\\
2.64492762787805	0.60756934437046\\
3.35364618075428	0.546134249262287\\
4.25226860165957	0.49635445769231\\
5.39168036402487	0.458065320975845\\
6.83640190002716	0.429466565712428\\
8.66824214034203	0.408607542616401\\
10.9909310339556	0.39357315741378\\
13.9359933695162	0.382640276098979\\
17.6701965097587	0.374372168487981\\
22.4049937750743	0.3677149626777\\
28.4084982181091	0.362056526617709\\
36.020664817419	0.356937604717942\\
45.6725407984347	0.352068741781101\\
57.9106741521309	0.347280871032798\\
73.4280625103569	0.342490069814027\\
93.1033948916044	0.3376670474672\\
118.050808423815	0.332817375430929\\
149.68297757285	0.327959034042448\\
189.791108372913	0.323129019173182\\
240.646367419354	0.318358723728575\\
305.128489151053	0.313683092209288\\
386.888843949848	0.309131117965972\\
490.557200965753	0.304734214053568\\
622.003893838172	0.300515186892446\\
788.672234732637	0.296487948212373\\
1000	0.292662483580413\\
};
\addlegendentry{Triple Momentum Method};

\addplot [color=mycolor3,mark size=2.5pt,only marks,mark=*,mark options={solid,fill=mycolor3}]
  table[row sep=crcr]{%
1.02329299228075	0.988553097160194\\
1.29748829388885	0.877906522889538\\
1.64515528346031	0.779644382409624\\
2.08598098298468	0.692380499306918\\
2.64492762787805	0.614883876641706\\
3.35364618075428	0.546061275724849\\
4.25226860165957	0.484941839571305\\
5.39168036402487	0.430663370544359\\
6.83640190002716	0.382460164168735\\
8.66824214034203	0.339652040855451\\
10.9909310339556	0.301635690247641\\
13.9359933695162	0.267874290818168\\
17.6701965097587	0.237891710972655\\
22.4049937750743	0.211264871788207\\
28.4084982181091	0.187618589115092\\
36.020664817419	0.166618543326396\\
45.6725407984347	0.147969573265824\\
57.9106741521309	0.131407663411517\\
73.4280625103569	0.116697670757833\\
93.1033948916044	0.103637296297367\\
118.050808423815	0.0920376619090705\\
149.68297757285	0.0817234332700407\\
189.791108372913	0.0725822172100945\\
240.646367419354	0.0644493432874639\\
305.128489151053	0.0572478164339619\\
386.888843949848	0.0508064444669131\\
490.557200965753	0.0451498143095696\\
622.003893838172	0.0400797822030554\\
788.672234732637	0.0356084634945976\\
1000	0.0316230835082188\\
};
\addlegendentry{Gradient Descent};

\addplot [color=mycolor4,mark size=2.5pt,only marks,mark=*,mark options={solid,fill=mycolor4}]
  table[row sep=crcr]{%
1.02329299228075	0.988586148657721\\
1.29748829388885	0.882027882436043\\
1.64515528346031	0.794370124366101\\
2.08598098298468	0.723883803057016\\
2.64492762787805	0.669137931908033\\
3.35364618075428	0.629037877388318\\
4.25226860165957	0.602996067034194\\
5.39168036402487	0.591353739445011\\
6.83640190002716	0.596342410163422\\
8.66824214034203	0.624585395768024\\
10.9909310339556	0.695623382015381\\
13.9359933695162	0.888832387391917\\
17.6701965097587	3.31327873506929\\
22.4049937750743	inf\\
28.4084982181091	inf\\
36.020664817419	inf\\
45.6725407984347	inf\\
57.9106741521309	inf\\
73.4280625103569	inf\\
93.1033948916044	inf\\
118.050808423815	inf\\
149.68297757285	inf\\
189.791108372913	inf\\
240.646367419354	inf\\
305.128489151053	inf\\
386.888843949848	inf\\
490.557200965753	inf\\
622.003893838172	inf\\
788.672234732637	inf\\
1000	inf\\
};
\addlegendentry{Heavy Ball Method};

\end{axis}
\end{tikzpicture}%
	\end{center}
	\caption{Upper bounds on the $H_2$-performance related to additive gradient noise
	for different optimization algorithms and condition ratios between $1.02$ and $1000$
	obtained using~\Cref{lemmaAnalysisH2}.
	In all cases, the dimensions of the Zames-Falb multipliers were chosen as {$ \dimAnticausal = 0 $,
	$ \dimCausal = 4 $}.}
	\label{figH2Analysis}
\end{figure}
\begin{figure}[ht!]
	\begin{center}
%
%
\definecolor{mycolor1}{rgb}{0.00000,0.44700,0.74100}%
\definecolor{mycolor2}{rgb}{0.85000,0.32500,0.09800}%
\definecolor{mycolor3}{rgb}{0.92900,0.69400,0.12500}%
\definecolor{mycolor4}{rgb}{0.49400,0.18400,0.55600}%
\begin{tikzpicture}

\begin{axis}[%
width=0.951\figurewidth,
height=\figureheight,
at={(0\figurewidth,0\figureheight)},
scale only axis,
xmode=log,
xmin=1,
xmax=1000,
xminorticks=true,
xlabel={Condition ratio $L/m$},
ylabel={$H_2$-performance},
ymode=log,
ymin=1e-05,
ymax=1,
yminorticks=true,
axis background/.style={fill=white},
legend style={at={(0.03,0.03)},anchor=south west,legend cell align=left,align=left,draw=white!15!black}
]
\addplot [color=mycolor1,only marks,mark=*,mark options={solid,fill=mycolor1}]
  table[row sep=crcr]{%
1.02329299228075	0.954744019462147\\
1.29748829388885	0.605963216621298\\
1.64515528346031	0.3946619248306\\
2.08598098298468	0.261569720237658\\
2.64492762787805	0.176790055445649\\
3.35364618075428	0.120370896922846\\
4.25226860165957	0.0831860543200114\\
5.39168036402487	0.0572941093581222\\
6.83640190002716	0.0402043961725882\\
8.66824214034203	0.0283710141308075\\
10.9909310339556	0.0198249493011395\\
13.9359933695162	0.0139468298904699\\
17.6701965097587	0.0100053357145601\\
22.4049937750743	0.00713664706083574\\
28.4084982181091	0.00510590471086536\\
36.020664817419	0.00361217041373601\\
45.6725407984347	0.00258858782348511\\
57.9106741521309	0.00183999848180137\\
73.4280625103569	0.00133648377699781\\
93.1033948916044	0.000963970386103095\\
118.050808423815	0.000684323408333818\\
149.68297757285	0.000494231586023779\\
189.791108372913	0.000369488802415231\\
240.646367419354	0.000255497649440165\\
305.128489151053	0.000164891298346631\\
386.888843949848	0.000124431124484346\\
490.557200965753	8.87827152921937e-05\\
622.003893838172	6.82974237560819e-05\\
788.672234732637	4.50411826941585e-05\\
1000	3.5107314510134e-05\\
};
\addlegendentry{Nesterov's Method};

\addplot [color=mycolor2,only marks,mark=*,mark options={solid,fill=mycolor2}]
  table[row sep=crcr]{%
1.02329299228075	0.976600196935738\\
1.29748829388885	0.75229166988042\\
1.64515528346031	0.561675959847349\\
2.08598098298468	0.410650394645257\\
2.64492762787805	0.296731874116579\\
3.35364618075428	0.21245461519408\\
4.25226860165957	0.15110094639101\\
5.39168036402487	0.106525649488244\\
6.83640190002716	0.0753926944427708\\
8.66824214034203	0.0532905127155199\\
10.9909310339556	0.0373002030697877\\
13.9359933695162	0.0261703224443414\\
17.6701965097587	0.0185800343951864\\
22.4049937750743	0.0130668244401583\\
28.4084982181091	0.00922581936298606\\
36.020664817419	0.00645694741836424\\
45.6725407984347	0.00458237249591113\\
57.9106741521309	0.0032159251505746\\
73.4280625103569	0.00230631100462642\\
93.1033948916044	0.00162975286712203\\
118.050808423815	0.00115202458517901\\
149.68297757285	0.000818401890260592\\
189.791108372913	0.000597459915405302\\
240.646367419354	0.000413430718217041\\
305.128489151053	0.000268234092877216\\
386.888843949848	0.000198891172183039\\
490.557200965753	0.000140538765085174\\
622.003893838172	0.000106129847815099\\
788.672234732637	7.03768245292143e-05\\
1000	5.37505814829613e-05\\
};
\addlegendentry{Triple Momentum Method};

\addplot [color=mycolor3,only marks,mark=*,mark options={solid,fill=mycolor3}]
  table[row sep=crcr]{%
1.02329299228075	0.976726710738169\\
1.29748829388885	0.762207753849123\\
1.64515528346031	0.583245617406832\\
2.08598098298468	0.438142540113632\\
2.64492762787805	0.32403663427265\\
3.35364618075428	0.236278565467668\\
4.25226860165957	0.169390206897533\\
5.39168036402487	0.11970533633624\\
6.83640190002716	0.0839126605482172\\
8.66824214034203	0.0582305120134707\\
10.9909310339556	0.0398143969450539\\
13.9359933695162	0.0271286206232095\\
17.6701965097587	0.0184917079243119\\
22.4049937750743	0.0123292165381882\\
28.4084982181091	0.0082609174082255\\
36.020664817419	0.00544769456296845\\
45.6725407984347	0.00367620849605631\\
57.9106741521309	0.00241358580183607\\
73.4280625103569	0.00163450020562211\\
93.1033948916044	0.00104066398342944\\
118.050808423815	0.000706796321464197\\
149.68297757285	0.0004484204665634\\
189.791108372913	0.00029314573126722\\
240.646367419354	0.000189973175774232\\
305.128489151053	0.000118853761885499\\
386.888843949848	7.86953596965827e-05\\
490.557200965753	4.96105161379406e-05\\
622.003893838172	3.33912232106491e-05\\
788.672234732637	2.10135576872057e-05\\
1000	1.35685074148269e-05\\
};
\addlegendentry{Gradient Descent};

\addplot [color=mycolor4,only marks,mark=*,mark options={solid,fill=mycolor4}]
  table[row sep=crcr]{%
1.02329299228075	0.97679143209998\\
1.29748829388885	0.768679026421152\\
1.64515528346031	0.601405621507534\\
2.08598098298468	0.468088194480475\\
2.64492762787805	0.363121467575609\\
3.35364618075428	0.280862449078522\\
4.25226860165957	0.215712126421549\\
5.39168036402487	0.164740539879127\\
6.83640190002716	0.125739297002123\\
8.66824214034203	0.0956007410449833\\
10.9909310339556	0.0720237769567625\\
13.9359933695162	0.0542820838294688\\
17.6701965097587	0.0410784168675066\\
22.4049937750743	0.0304932522217804\\
28.4084982181091	0.0227905543670198\\
36.020664817419	0.0168053707125002\\
45.6725407984347	0.0126908181100361\\
57.9106741521309	0.00935165873092465\\
73.4280625103569	0.00709021074723033\\
93.1033948916044	0.00506785629756977\\
118.050808423815	0.00387237120880326\\
149.68297757285	0.00276509986546831\\
189.791108372913	0.00202163005295546\\
240.646367419354	0.00148059473150475\\
305.128489151053	0.00103979830988667\\
386.888843949848	0.000773421935630782\\
490.557200965753	0.000552241117167733\\
622.003893838172	0.000424561208453707\\
788.672234732637	0.000295110550742672\\
1000	0.000215639226629231\\
};
\addlegendentry{Heavy Ball Method};

\end{axis}
\end{tikzpicture}%
	\end{center}
	\caption{{Lower} bounds on the $H_2$-performance related to additive gradient noise
	for different optimization algorithms and condition ratios between $1.02$ and $1000$
	obtained from a sample-based approach.}
	\label{figH2AnalysisSampleBased}
\end{figure}}
These numerical results suggest that the Gradient Descent algorithm 
has the best properties in terms of additive noise attenuation for
condition ratios larger than approximately $10$, followed
by Nesterov's Method. The fastest method in terms of convergence rates,
the Triple Momentum Method, however, has the worst noise attenuation.
{We note that these results have to be taken with care since they only provide upper bounds; still, they are}
{qualitatively in accordance with the results from~\citet{mohammadi2018variance},
where a similar performance channel has been analyzed for quadratic optimization problems.}
To underpin this statement, we also carried out a sample-based approach
to evaluate the $H_2$-performance. To this end, we randomly generated 
$10000$ functions from the class $ \classObj{m}{L} $ and simulated the four
considered optimization algorithms under additional noise. The 
corresponding {lower} $H_2$-performance bound estimates are depicted \ifthenelse{\boolean{singleColumn}}{on the right-hand side of~\Cref{figH2Analysis}}{in~\Cref{figH2AnalysisSampleBased}}.
While the results are qualitatively similar, still, there is quite a
gap quantitatively. We emphasize that this is not a contradiction 
since, first,~\Cref{lemmaAnalysisH2} does only provide an upper bound,
and second -- and probably more important -- the sampling of the 
class $ \classObj{m}{L} $ is not very dense. In fact, we only sample quadratic
functions $ \fObj \in \classObj{m}{L} $ as well as functions with a Hessian of 
the form $ \nabla^2 \fObj(\optVar) = c_1 + c_2 \cos(\omega \optVar) $,
where $ m \leq c_1 + c_2 \leq L $ {and} $ \omega > 0 $ {are chosen randomly.}

The previous results {suggest} that there is a trade-off between convergence speed and
robustness {against} noise, an observation that has also been made in the 
literature in different settings, see~\citet{polyak1987introduction} or \citet{lessard2016analysis} in the case of
relative deterministic noise. Our numerical results underpin 
and quantify these findings.
This observation becomes even more apparent when designing novel
optimization algorithms, see~\Cref{secSynthesisNumericalResults}.
Similarly to the convergence rate analysis, anticausal Zames-Falb
multipliers do not lead to improved bounds.

\subsection{Synthesis}
In the following we consider the problem of designing algorithms that have
desired properties specified in terms of convergence rates and $H_2$-performance 
bounds. {To this end, roughly speaking, we need to combine~\Cref{lemmaAnalysisConvRate}
and~\Cref{lemmaAnalysisH2} with the difference that the algorithm parameters 
$ \Aopt, \Bopt, \Copt, \Dopt $ are now to be determined as well.}
Motivated by our discussions in~\Cref{secProblemFormulation}, we 
fix $\Copt$, $ \Dopt $ and $\Dopt^\dagger$ to 
\begin{align}
	\Copt = \begin{bmatrix} 1 & 0 & \dots & 0 \end{bmatrix} \otimes I_p, \quad \Dopt = \Copt, \quad \Dopt^\dagger = \Copt^\top.
	\label{eqDefCoptDdagger}
\end{align}
{Still, having $ \Aopt, \Bopt $ as additional design variables, 
the following main difficulties compared to the analysis problem
need to be addressed: 
(1) The conditions~\eqref{eqConditionABCD} have to be ensured to
hold; (2) nominal exponential stability of the feedback interconnection
has to be guaranteed, i.e., $ \Aopt + m \Bopt \Copt $ must have all
its eigenvalues in the open disk of radius $ \rho $; and
(3) the introduction of new design variables leads to bilinear
instead of linear matrix inequalities in~\eqref{eqAllConstraintsConvAnalysis}, \eqref{eqLMIPerfAnalysis},
which, in general, cannot be solved efficiently. We emphasize that
in the end we are interested in finding LMI conditions, which is the 
core difficulty of all three problems. As it turns out, 
(1) can be handled without much effort and, as we show next,
(2) can be resolved utilizing the special structure of the 
problem at hand.}
\begin{lemma}[Nominal exponential stability]\label{lemmaNominalExponentialStability}
Let $ \rho \in (0,1) $ be given and consider $ \AcompleteNoPerf, \BcompleteNoPerf, \CcompleteNoPerf, \DcompleteNoPerf $
	as defined by~\eqref{eqSSCompleteSystemNoPerformance}, cf.~\eqref{secAppendixStateSpaceRealizationNoPerf}.
	Suppose that there exists $ P $ partitioned as
	\ifthenelse{\boolean{singleColumn}}{
	\begin{align}
		P 
		=
		\begin{bmatrix} P_{11} & P_{12} \\ P_{12}^\top & P_{22} \end{bmatrix},  
		P_{11} \in \real^{p(\dimCausal+\dimCausal+1) \times p(\dimCausal+\dimCausal+1)}, P_{22} \in \real^{np \times np}, \label{eqPartitioningPexpStab}
	\end{align}
	}{
	\begin{align}
		&P 
		=
		\begin{bmatrix} P_{11} & P_{12} \\ P_{12}^\top & P_{22} \end{bmatrix}, \label{eqPartitioningPexpStab} \\
		&P_{11} \in \real^{p(\dimCausal+\dimCausal+1) \times p(\dimCausal+\dimCausal+1)}, P_{22} \in \real^{np \times np}, \nonumber
	\end{align}}
	such that~\eqref{eqLMIConvAnalysis} holds.
	Then $ \Aopt + m \Bopt \Copt $ has all eigenvalues in the open
	disk of radius $ \rho $ {if and only if $ P_{22} \succ 0 $}.
\end{lemma}
\begin{proof}
Consider the lower right block of~\eqref{eqLMIConvAnalysis}.
	The inequality then implies that 
	\ifthenelse{\boolean{singleColumn}}{
	\begin{align}
		\rho^{-2} ( \Aopt + m \Bopt \Copt )^\top P_{22} ( \Aopt + m \Bopt \Copt ) - P_{22} 
		+ \rho^{-2} C_{\psi_1}^\top D_\Delta^\top M_\Delta(M_+, M_-, M_0) D_\Delta C_{\psi_1} \prec 0.
		\label{eqIneqNominalStab}
	\end{align}}{
	\begin{align}
		\begin{split}
		\rho^{-2} ( \Aopt + m \Bopt \Copt )^\top P_{22} ( \Aopt + m \Bopt \Copt ) - P_{22} \qquad \qquad \\
		+ \rho^{-2} C_{\psi_1}^\top D_\Delta^\top M_\Delta(M_+, M_-, M_0) D_\Delta C_{\psi_1} \prec 0.
		\end{split}
		\label{eqIneqNominalStab}
	\end{align}}
	By straightforward calculations, it is seen that $ C_{\psi_1}^\top D_\Delta^\top M_\Delta(M_+, M_-, M_0) D_\Delta C_{\psi_1} = 0$;
	thus, {by standard Lyapunov theory for linear systems}, we infer that $ \Aopt + m \Bopt \Copt $
	has all eigenvalues in the open disk of radius $ \rho $ {if and only if $ P_{22} \succ 0 $}.
\end{proof}
{In order to address~(3), it is convenient to first reformulate the matrix inequalities
by means of Schur complements.}
As apparent in~\eqref{eqIneqNominalStab}, with $ \Aopt $, $ \Bopt $ being design variables in 
the synthesis case, \eqref{eqLMIConvAnalysis} is a quadratic matrix inequality 
in the unknowns. The same holds true for the performance inequalities 
\eqref{eqLMIH2Analysis}. As a first step, we utilize the positive definiteness
condition for nominal exponential stability and a similar condition for performance
in order to reformulate these quadratic matrix inequalities to bilinear
matrix inequalities (BMIs). {Note that with~\eqref{eqDefN} and~\eqref{eqPartitioningPexpStab}
we have
\ifthenelse{\boolean{saveSpace}}{{$P_{22} = N^\top P N$;}}{
\begin{align}
	P_{22} = N^\top P N;
\end{align}}
hence}
the positive definiteness constraint
$ P_{22} \succ 0 $ allows us to equivalently formulate~\eqref{eqLMIConvAnalysis} 
as {(see~\Cref{secAppendixDerivationH2LMIForSynthesis} for a derivation)}
{
\begin{align}
	\begin{bmatrix}
		- P_{22} & P_{22} N^\top \begin{bmatrix} \AcompleteNoPerf(\rho) & \BcompleteNoPerf \end{bmatrix} \\
		\star    & U(P,M_\Delta) 		
	\end{bmatrix}
	\prec 0,
	\label{eqReformulationConvRateLMIForSynthesis}
\end{align}}
where {we introduced the shorthand notation}
\ifthenelse{\boolean{singleColumn}}{
\begin{align}
	U(P,M_\Delta) 
	= 
	\left[ 
	\vphantom{
	\begin{array}{cc}
		\AcompleteNoPerf & \BcompleteNoPerf \\
		I                & 0 \\
		\CcompleteNoPerf & \DcompleteNoPerf
	\end{array}} \star
	\right]^\top
	\left[ 
	\begin{array}{ccc}
		P-N P_{22} N^\top & 0  & 0 \\
		0   & -P & 0 \\
		0   & 0  & M_\Delta
	\end{array}
	\right] 
	\left[ 
	\begin{array}{cc}
		{\AcompleteNoPerf(\rho)} & \BcompleteNoPerf \\
		I                & 0 \\
		\CcompleteNoPerf & \DcompleteNoPerf
	\end{array}
	\right].  	\label{eqSynthConvRateLowerRightBlock}  
\end{align}}{
\begin{align}
	&U(P,M_\Delta) \label{eqSynthConvRateLowerRightBlock}  \\
	&= 
	\left[ 
	\vphantom{
	\begin{array}{cc}
		\AcompleteNoPerf & \BcompleteNoPerf \\
		I                & 0 \\
		\CcompleteNoPerf & \DcompleteNoPerf
	\end{array}} \star
	\right]^\top
	\left[ 
	\begin{shorterArray}{ccc}
		P-N P_{22} N^\top & 0  & 0 \\
		0   & -P & 0 \\
		0   & 0  & M_\Delta
	\end{shorterArray}
	\right] 
	\left[ 
	\begin{array}{cc}
		{\AcompleteNoPerf(\rho)} & \BcompleteNoPerf \\
		I                & 0 \\
		\CcompleteNoPerf & \DcompleteNoPerf
	\end{array}
	\right]. \nonumber 	
\end{align}}
{Note that, by the structure of the state-space
representation~\eqref{secAppendixStateSpaceRealizationNoPerf}, 
\eqref{eqReformulationConvRateLMIForSynthesis} is bilinear in 
the unknowns $ \Aopt, \Bopt, P, M_\Delta $.}

A similar procedure applies to the problem with additional performance specifications.
More precisely, {assuming $P_\perf$ to be structured in the same manner as $P$,}~\eqref{eqLMIH2Analysis}, \eqref{eqPosDefH2Analysis}
are equivalently formulated as
{
\begin{align}
	\begin{bmatrix}
		-P_{\perf,22} & P_{\perf,22} N^\top \begin{bmatrix} \Acomplete & \Bcomplete{}_{,1} \end{bmatrix} \\
		\star       & 
		U_\perf(P_\perf,M_\Delta)
	\end{bmatrix}
	\prec
	0,
	\label{eqReformulationH2LMIForSynthesis}
\end{align}}
where
\ifthenelse{\boolean{singleColumn}}{
\begin{align}
	U_\perf(P_\perf,M_\Delta) = 
	\left[
	\vphantom{
	\begin{array}{ccc}
		\Acomplete        & \Bcomplete{}_{,1}  & \Bcomplete{}_{,2} \\
		I                 & 0                  & 0 \\ \hline 
		\Ccomplete{}_{,1} & \Dcomplete{}_{,11} & \Dcomplete{}_{,12} \\ 
		\Ccomplete{}_{,2} & \Dcomplete{}_{,21} & \Dcomplete{}_{,22}
	\end{array}} \star
	\right]^\top
	\left[
	\begin{array}{cc|cc}
		P_\perf-NP_{\perf,22}N^\top & 0         & 0        & 0 \\
		0         & - P_\perf & 0        & 0 \\ \hline
		0         & 0         & M_\Delta & 0 \\
		0         & 0         & 0        & I_{n_{\yperf}} 
	\end{array}
	\right]
	\left[
	\begin{array}{ccc}
		\Acomplete        & \Bcomplete{}_{,1}  \\
		I                 & 0                  \\ \hline 
		\Ccomplete{}_{,1} & \Dcomplete{}_{,11} \\ 
		\Ccomplete{}_{,2} & 0                  
	\end{array}
	\right]. 
	\label{eqH2LMIForSynthesisLowerRightBlock}
\end{align}}{
\begin{flalign}
	&U_\perf(P_\perf,M_\Delta) \label{eqH2LMIForSynthesisLowerRightBlock} \\
	&= \left[
	\vphantom{
	\begin{array}{ccc}
		\Acomplete        & \Bcomplete{}_{,1}  & \Bcomplete{}_{,2} \\
		I                 & 0                  & 0 \\ \hline 
		\Ccomplete{}_{,1} & \Dcomplete{}_{,11} & \Dcomplete{}_{,12} \\ 
		\Ccomplete{}_{,2} & \Dcomplete{}_{,21} & \Dcomplete{}_{,22}
	\end{array}} \star
	\right]^\top
	\left[
	\begin{shortArray}{cc|cc}
		P_\perf-NP_{\perf,22}N^\top & 0         & 0        & 0 \\
		0         & - P_\perf & 0        & 0 \\ \hline
		0         & 0         & M_\Delta & 0 \\
		0         & 0         & 0        & I_{n_{\yperf}} 
	\end{shortArray}
	\right]
	\left[
	\begin{shortArray}{ccc}
		\Acomplete        & \Bcomplete{}_{,1}  \\
		I                 & 0                  \\ \hline 
		\Ccomplete{}_{,1} & \Dcomplete{}_{,11} \\ 
		\Ccomplete{}_{,2} & 0                  
	\end{shortArray}
	\right]. &&
	\nonumber
\end{flalign}}
In the following we consider a performance channel defined as in the example in~\Cref{secAnalysisExistingAlgorithms},
i.e., $ \Bperf = \Bopt $, $ \Cperf = \Copt $, $ \Dperf = 0 $. 
Note that, by {this} choice, $ \Bperf $ is as well a design variable, 
while $ \Cperf, \Dperf $ are fixed,
{and~\eqref{eqReformulationH2LMIForSynthesis} is bilinear in the unknowns
$ \Aopt, \Bopt, P, M_\Delta $.}
Additionally, keeping in mind~\Cref{remarkReformulationTrace} and employing
the specific state-space realization from~\Cref{secAppendixStateSpaceRealizations}, \eqref{eqTraceH2Analysis}
can be reformulated as
\begin{align}
	\begin{bmatrix}
		{P_{\perf,22}} & {P_{\perf,22}} N^\top \Bcomplete{}_{,2} \\ \star & Z 
	\end{bmatrix} \succ 0, \quad
	\textup{tr}(Z) < \gamma^2.
	\label{eqTraceH2reformulated}
\end{align}
Still, for the synthesis problem all of the former inequalities are bilinear. 
We next discuss two approaches to cope with 
that fact. The first one relies on rendering the conditions linear by imposing
certain restrictions on the design variables and employing a suitable variable transformation;
the second one is {an often used} rather hands-on heuristic based on alternately solving LMIs to obtain feasible solutions of the BMIs.

{We formulate the subsequent synthesis results for full 
matrices $ \Aopt, \Bopt, \Copt, \Dopt $. We emphasize that 
a Kronecker type structure can be enforced simply by desigining
the matrices for $p=1$ and then adapt the obtained
matrices to the required dimension $\bar{p}$ taking Kronecker products with $I_{\bar{p}}$.
This enables a design independent of the dimension of the optimization
variable.}

\subsubsection{{Convex solution}}\label{subsecConservativeConvexSynthesis}
We next discuss how to obtain LMI conditions from the BMIs~\eqref{eqReformulationConvRateLMIForSynthesis},
\eqref{eqReformulationH2LMIForSynthesis}, \eqref{eqTraceH2reformulated}.
Variable transformations have successfully been employed in many
situations to render the inequalities linear in the new variables.
However, since $ {\AcompleteNoPerf(\rho)}, \BcompleteNoPerf, \Acomplete, \Bcomplete $
are structured matrices here, the standard variable {transformations do}
not directly apply. To make this more vivid, consider~\eqref{eqReformulationConvRateLMIForSynthesis}.
With $ \Aopt, \Bopt, P $ being design variables, a typical approach
is to define $ Q_\Aopt = P_{22} \Aopt $, $ Q_\Bopt = P_{22} \Bopt $.
However, $U$ as defined in~\eqref{eqSynthConvRateLowerRightBlock}
also contains nonlinear terms of the form $ P_{12} \Aopt = P_{12} P_{22}^{-1} Q_\Aopt $, $ P_{12} \Bopt = P_{12} P_{22}^{-1} Q_\Bopt $,
{where $P_{12}$ is the right upper block of $P$, see~\eqref{eqPartitioningPexpStab}.}
These terms cannot be handled by that approach; a simple remedy is to let $ P_{12} = 0 $,
i.e., $ P $ is block diagonal, which is the idea behind the following Theorem.
\begin{theorem}\label{lemmaConservativeConvexSynthesisNew} 
Let {$ L \geq m > 0 $} as well as $ n \in \natPos $, $ p \in \natPos $ be given.
	Let $ \Copt \in \real^{p\times np} $, $ \Dopt \in \real^{np \times p} $, $ \Dopt^\dagger \in \real^{p \times np} $
	be defined as in~\eqref{eqDefCoptDdagger}.
	Fix $ \rho \in (0,1) $. Let some $ \dimCausal, \dimAnticausal \in \nat $ be given 
	and let $ M_\Delta $, $ \mathbb{M} $ be defined
	according to~\eqref{eqDefMDelta}, \eqref{eqDefSetMunstructured}.
	Suppose there exist $ P_{11} = P_{11}^\top {\in \real^{p(\dimCausal+\dimAnticausal) \times p(\dimCausal+\dimAnticausal)}} $, $P_{22} = P_{22}^\top {\in \real^{np \times np}} $, $Q_{\Aopt} {\in \real^{np \times np}} , Q_{\Bopt} {\in \real^{np \times p}}  $, 
	$ M_+ \in \real^{p \times \dimCausal p}$, $M_- \in \real^{p \times \dimCausal p}$, $M_0 \in \real^{p \times p} $ 
	such that
	\begin{subequations}
	\begin{align}
		\begin{bmatrix}
			- P_{22} & \begin{bmatrix} 0 & \rho^{-1} Q_{\Aopt} + m \rho^{-1} Q_{\Bopt} \Copt & Q_\Bopt \end{bmatrix} \\
			\star    & U({P},M_\Delta)
		\end{bmatrix}
		\prec 0 \label{eqSynthConvRate1} \\
		{P = \textup{blkdiag}(P_{11},P_{22})} \\
		( Q_{\Aopt} - P_{22} ) \Dopt^\dagger = 0 \label{eqSynthConvRate3} \\
		\begin{bmatrix}  M_- & M_0 &  M_+ \end{bmatrix} \in \mathbb{M}(\rho,\dimCausal,\dimCausal,p). \label{eqSynthConvRate4}
	\end{align}\label{eqSynthConvRate}
	\end{subequations}
	Then, with $ \Aopt = P_{22}^{-1} Q_{\Aopt} $, $ \Bopt = P_{22}^{-1} Q_{\Bopt} $,
	the equilibrium $ \xOpt^\optSign = \Dopt \optVar^\optSign $ is globally
	robustly exponentially stable against $ \mathbf{\Delta}(m,L) $ with rate $ \rho $
	for~\eqref{eqRobustStandardForm}.
\end{theorem}
\begin{proof}
The result follows directly from~\Cref{lemmaAnalysisConvRate} noting that~\eqref{eqSynthConvRate1}, \eqref{eqSynthConvRateLowerRightBlock}, \eqref{eqSynthConvRate4}
	imply~\eqref{eqLMIConvAnalysis}, \eqref{eqSynthConvRate1} implies $ P_{22} \succ 0 $ and
	thus nominal exponential stability by~\Cref{lemmaNominalExponentialStability},
	and~\eqref{eqSynthConvRate3} together with the choice of $ \Copt $, $ \Dopt^\dagger $
	implies~\eqref{eqConditionABCD}.
\end{proof}
We note that the latter result provides a convex solution to the original
design problem introduced in~\Cref{secProblemFormulation}; however, as it turns out
in numerical examples, the restrictions on $P$ are conservative
and we trade off convexity for slower convergence rates.  

When it comes to the synthesis problem with additional performance
specifications, similar restrictions on $P_\perf$ can be employed.
Building upon that, the following Theorem provides a convex solution 
to this extended problem. \ifthenelse{\boolean{proofsAppendix}}{A proof is provided in~\Cref{secProofLemmaConvexSynthesisWithPerf}.}{}
{
\begin{theorem}\label{lemmaConservativeConvexSynthesisWithPerf}
Let {$ L \geq m > 0 $} as well as $ n \in \natPos $, $ p \in \natPos $ be given.
	Let $ \Copt \in \real^{p\times np} $, $ \Dopt \in \real^{np \times p} $, $ \Dopt^\dagger \in \real^{p \times np} $
	be defined as in~\eqref{eqDefCoptDdagger}.
	Fix $ \rho \in (0,1) $. Let some $ \dimCausal, \dimAnticausal \in \nat $ be given 
	and let $ M_\Delta $, $ \mathbb{M} $ be defined
	according to~\eqref{eqDefMDelta}, \eqref{eqDefSetMunstructured}.
	Let further $ \Bperf = \Bopt $, $ \Cperf = \Copt $, $ \Dperf = 0 $.
	Suppose there exist 
	$ P_{11} = P_{11}^\top {\in \real^{p(\dimCausal+\dimAnticausal) \times p(\dimCausal+\dimAnticausal)}} $, $P_{22} = P_{22}^\top {\in \real^{np \times np}} $, $Q_{\Aopt} {\in \real^{np \times np}} , Q_{\Bopt} {\in \real^{np \times p}}  $,
	$ P_{\perf,11} = P_{\perf,11}^\top {\in \real^{p(\dimCausal+\dimAnticausal) \times p(\dimCausal+\dimAnticausal)}} $, $Z=Z^\top {\in \real^{p \times p}}$, 
	$ M_+, M_{\perf,+} \in \real^{p \times \dimCausal p}$, $M_-, M_{\perf,-} \in \real^{p \times \dimCausal p}$, $M_0, M_{\perf,0} \in \real^{p \times p} $ 
	such that~\eqref{eqSynthConvRate} holds and
	\begin{subequations}
	\begin{align}
		\left[
		\begin{shortArray}{cc}
			-P_{22} & \begin{bmatrix} 0 & Q_{\Aopt} + m Q_{\Bopt} \Copt & Q_\Bopt \end{bmatrix} \\
			\star   & U_\perf(P_\perf,M_{\perf,\Delta}) 
		\end{shortArray}\right]
		\prec 0 \label{eqConservativeConvexSynthesisWithPerf1} \\
		P_{\perf}
		=
		\textup{blkdiag}( P_{\perf,11}, P_{22} ) \label{eqConservativeConvexSynthesisWithPerf2} \\
		\begin{bmatrix}
			P_{22} & Q_\Bopt \\
			\star  & Z 
		\end{bmatrix} 
		\succ 0 \label{eqConservativeConvexSynthesisWithPerf3} \\
		\textup{tr}(Z) < \gamma^2, \label{eqConservativeConvexSynthesisWithPerf4} \\
		\begin{bmatrix}  M_{\perf,-} & M_{\perf,0} &  M_{\perf,+} \end{bmatrix} \in \mathbb{M}(1,\dimCausal,\dimCausal,p),
	\end{align}\label{eqSynthConvexConvRateWithPerf}
	\end{subequations}
	where $ M_{\perf,\Delta} $ is defined as $ M_\Delta $ in~\eqref{eqDefMDelta}
	replacing $ M_+, M_-, M_0 $ by their counterparts $ M_{\perf,+}, M_{\perf,-}, M_{\perf,0} $.
	Then, with $ \Aopt = P_{22}^{-1} Q_{\Aopt} $, $ \Bopt = P_{22}^{-1} Q_{\Bopt} $,
	the equilibrium $ \xOpt^\optSign = \Dopt \optVar^\optSign $ is globally
	robustly exponentially stable against $ \mathbf{\Delta}(m,L) $ with rate $ \rho $
	for~\eqref{eqRobustStandardForm} and~\eqref{eqRobustStandardForm} achieves
	a robust $ H_2 $-performance level {of} $ \gamma $ against $ \mathbf{\Delta}(m,L) $.
\end{theorem}}
\ifthenelse{\boolean{proofsAppendix}}{
}
{\begin{proof} Robust exponential stability follows directly from~\Cref{lemmaConservativeConvexSynthesisNew}.
For robust performance, we first note that~\eqref{eqConservativeConvexSynthesisWithPerf1}, \eqref{eqConservativeConvexSynthesisWithPerf2}
implies that~\eqref{eqReformulationH2LMIForSynthesis} holds with {$ P_{\perf,22} = P_{22} $},
$ Q_\Aopt = P_{22} \Aopt $, $ Q_{\Bopt} = P_{22} \Bopt $, and making use of the 
specific state-space realization from~\Cref{secAppendixStateSpaceRealizations},
thus~\eqref{eqConservativeConvexSynthesisWithPerf1}, \eqref{eqConservativeConvexSynthesisWithPerf2} implies~\eqref{eqLMIH2Analysis}.
Note further that $ U_\perf $ as defined in~\eqref{eqH2LMIForSynthesisLowerRightBlock}
is independent of $ \Aopt, \Bopt $ by the specific structure of $ P_\perf $.
We further note that~\eqref{eqConservativeConvexSynthesisWithPerf3} holds 
if and only if $ P_{22} \succ 0 $ and $ Z \succ Q_\Bopt^\top P_{22}^{-1} Q_{\Bopt} = \Bopt^\top P_{22} \Bopt $.
Together with that the trace condition~\eqref{eqConservativeConvexSynthesisWithPerf4} then implies
that $ \textup{tr}( \Bopt^\top P_{22} \Bopt ) < \gamma^2 $. Now note that
with the specific state-space realization~\Cref{secAppendixStateSpaceRealizations}
and by the definition of $ P_\perf $ we have that $ \textup{tr}( \Bcomplete{}_{,2}^\top P_\perf \Bcomplete{}_{,2} ) = \textup{tr}(\Bopt^\top P_{22} \Bopt ) $,
hence~\eqref{eqConservativeConvexSynthesisWithPerf3} together with~\eqref{eqConservativeConvexSynthesisWithPerf4}
implies~\eqref{eqTraceH2Analysis}, thus concluding the proof.
 \end{proof}
}

\subsubsection{Synthesis based on BMI optimization techniques}\label{secBMIoptimization}
While~\Cref{lemmaConservativeConvexSynthesisNew}, \Cref{lemmaConservativeConvexSynthesisWithPerf} provide a convex synthesis procedure,
the block diagonal structure of $P$ as well as the assumptions on $ P_\perf $ are restrictive and result in convergence 
rates inferior to what can be achieved, see~\Cref{secSynthesisNumericalResults}. As an alternative, BMI optimization
techniques, which directly try to solve~\eqref{eqReformulationConvRateLMIForSynthesis},
\eqref{eqReformulationH2LMIForSynthesis}, \eqref{eqTraceH2reformulated}
or variants thereof, can be employed. However, while many of the approaches perform
reasonably well as we will also illustrate in the subsequent example section, the non-convex nature of the problem does not allow for any guarantees
of finding the global optimizer in general. Here, we employ the so-called alternating
method which alternates between solving two different semi-definite programs obtained
from fixing two subsets of the set of all decision variables in the BMI. In particular,
assuming that $ \Copt, \Dopt, \Dopt^\dagger $ are chosen as in~\eqref{eqDefCoptDdagger},
we alternate between finding the best algorithm $ \Aopt, \Bopt $ in terms of $ \rho $ and $ \gamma $ for fixed $ P, P_\perf $
and solving for $ P, P_\perf $ for this $ \Aopt, \Bopt, \rho, \gamma $.
Building upon an iterative scheme, having a good initial feasible solution is key for successfully applying this procedure. As it turns out
in numerical examples, {certain types of parametrized algorithms are better
suited for initialization} than making use of~\Cref{lemmaConservativeConvexSynthesisNew}.
{In particular, we {let} $ \Aopt, \Bopt $
{be} parametrized by the free {matrices} $ K_i \in \real^{p \times p} $, $ i = 1,2,\dots,n $, as
\ifthenelse{\boolean{singleColumn}}{
\begin{align}
	\Aopt = \Aopt_1 + I_{np} + \Bopt_1 \begin{bmatrix} 0 & K_2 & \dots & K_n \end{bmatrix}, \qquad
	\Bopt = \Bopt_1 K_1, \label{eqStructuredAlgoAB}
\end{align}}{
\begin{subequations}
\begin{align}
	\Aopt &= \Aopt_1 + I_{np} + \Bopt_1 \begin{bmatrix} 0 & K_2 & \dots & K_n \end{bmatrix} \\
	\Bopt &= \Bopt_1 K_1,
\end{align}\label{eqStructuredAlgoAB}
\end{subequations}}
where $ \Aopt_1 \in \real^{np \times np} $, $ \Bopt_1 \in \real^{np \times p} $
are defined as
\begin{align}
	\Aopt_1 = 
	\left[
	\begin{shortArray}{ccccc} 
	0      & I_p     & \\
	0      & 0     & I_p & \\
	\smash[t]{\vdots}  &       &   & \smash[t]{\ddots}  & \\
	\smash[t]{\vdots}  &       &   &        & I_p \\
	0      & \dots &   &        & 0
	\end{shortArray}\right],
	\Bopt_1 = 
	\begin{bmatrix}
		0 \\ 0 \\ \smash[t]{\vdots} \\ \smash[t]{\vdots} \\ I_p
	\end{bmatrix}.
\end{align}
The particular structure is motivated by~\citet{mic2014heavy}, 
\citet{mic2016extremum} and corresponds to  
a Euler discretization of the $n$-th order heavy ball method
presented in the latter references. For the initialization of the
BMI iteration we then fix $ \Copt, \Dopt $ as in~\eqref{eqDefCoptDdagger} 
and utilize~\Cref{lemmaSynthesisParametrizedAlgos} to find suitable
parameters $ K_i $, $ i = 1,2,\dots,n $, see~\Cref{secAppendixDesignStructuredAlgorithms}.
}

\subsubsection{Numerical results}\label{secSynthesisNumericalResults}
In this section we employ the presented results to design novel
optimization algorithms {that guarantee a priori specified convergence rates 
and that are robust against additive noise, i.e., the performance channel 
is defined as explained in~\Cref{secAnalysisExistingAlgorithms}. 
More precisely, for a fixed condition ratio $\kappa = L/m$, we
choose different desired convergence rates
and then design algorithms that ensure this convergence rate
and additionally 
have a minimal bound on the $H_2$-performance level.
The numerical results for two specific condition ratios (left: $ \kappa = 50 $, right: $ \kappa = 100 $) are depicted in~\Cref{figSynthAlgoAdditiveNoise}.
This figure also shows the convergence 
rate and $H_2$-performance levels obtained from analyzing
the Triple Momentum Method, Nesterov's Method and the Gradient
Descent algorithm using~\Cref{lemmaAnalysisConvRate}, \Cref{lemmaAnalysisH2}.
Apparently, we are not only able to recover the convergence rate
and robustness properties of these existing algorithm but our 
design approach enables us to weigh up these two performance specifications
against each other. More specifically, as already observed in our analysis in~\Cref{figH2Analysis},
there is a trade-off between convergence rates and robustness against noise
and the newly generated algorithms lie on the corresponding Pareto frontier. 
The numeric results also show that 
higher order algorithms (i.e., a larger dimension $n$ of the state in~\eqref{eqOptAlgo})
can yield better 
$H_2$-performance levels while providing the same convergence rate
guarantees. Still, we cannot design an algorithm with better
convergence rate guarantees than the Triple Momentum Method.
Interestingly, the convex approach from~\Cref{lemmaConservativeConvexSynthesisWithPerf}
might yield even better performance bounds; however, the approach is limited
to convergence rates worse than the one of the Gradient Descent algorithm.}
\begin{figure}[t!]
	\ifthenelse{\boolean{singleColumn}}{
	\setlength\figurewidth{0.4\textwidth}
	\setlength\figureheight{0.4\textwidth}
	\begin{center}
	\begin{footnotesize}
%
%
\definecolor{mycolor1}{rgb}{0.00000,0.39062,0.00000}%
\definecolor{mycolor2}{rgb}{0.29297,0.00000,0.50781}%
\definecolor{mycolor3}{rgb}{0.85000,0.32500,0.09800}%
\definecolor{mycolor4}{rgb}{0.00000,0.44700,0.74100}%
\definecolor{mycolor5}{rgb}{0.92900,0.69400,0.12500}%
\begin{tikzpicture}

\begin{axis}[%
width=0.951\figurewidth,
height=\figureheight,
at={(0\figurewidth,0\figureheight)},
scale only axis,
unbounded coords=jump,
xmin=0.85,
xmax=1,
xlabel={Convergence rate $\rho$},
ymin=0.005,
ymax=0.4,
ymode=log,
ylabel={$H_2$-performance},
axis background/.style={fill=white},
legend style={at={(0.04,0.04)},anchor=south west,legend cell align=left,align=left,draw=white!15!black,font=\small},
grid=both,
grid style={line width=.1pt, draw=black!30},
major grid style={line width=.2pt,draw=black!80},
ytick={0.01,0.1}
]
\addplot [color=mycolor1,dashed,mark=*,mark options={solid,fill=mycolor1}]
  table[row sep=crcr]{%
0.862871536981504	0.30946152453036\\
0.868543556273941	0.288298636481586\\
0.874215575566379	0.267328992115637\\
0.879887594858816	0.251804914556189\\
0.885559614151253	0.239264453926221\\
0.891231633443691	0.230883435261668\\
0.896903652736128	0.224452274367907\\
0.902575672028565	0.218636272483777\\
0.908247691321003	0.213041078515733\\
0.91391971061344	0.208782422389828\\
0.919591729905877	0.204608588573787\\
0.925263749198315	0.197687811804989\\
0.930935768490752	0.192445233553485\\
0.936607787783189	0.184933986506113\\
0.942279807075627	0.178715952597561\\
0.947951826368064	0.166866285282643\\
0.953623845660501	0.15461276337932\\
0.959295864952939	0.14558710996721\\
0.964967884245376	0.134695296588531\\
0.970639903537813	0.124890226821288\\
0.976311922830251	0.116701440957512\\
0.981983942122688	0.113124350381597\\
0.987655961415125	0.102805955772367\\
0.993327980707563	0.101259898070055\\
0.999	0.101675364929887\\
};
\addlegendentry{BMI synthesis: n = 2};

\addplot [color=mycolor1,dashed,mark=*,mark options={solid,fill=white!50!mycolor1}]
  table[row sep=crcr]{%
0.862871536981504	nan\\
0.868543556273941	0.288481688089103\\
0.874215575566379	0.266947234996696\\
0.879887594858816	0.251914680311602\\
0.885559614151253	0.236312288366904\\
0.891231633443691	0.224751520970588\\
0.896903652736128	0.2148201122362\\
0.902575672028565	0.202914257308799\\
0.908247691321003	0.195601319561357\\
0.91391971061344	0.18705395991594\\
0.919591729905877	0.185284101951933\\
0.925263749198315	0.177876477563977\\
0.930935768490752	0.168862293049972\\
0.936607787783189	0.162758399794006\\
0.942279807075627	0.158858866727713\\
0.947951826368064	0.15561693070124\\
0.953623845660501	0.149971693472256\\
0.959295864952939	0.144353607176755\\
0.964967884245376	0.137964463567028\\
0.970639903537813	0.133212749915167\\
0.976311922830251	0.122083768435539\\
0.981983942122688	0.116669364452387\\
0.987655961415125	0.100943799672138\\
0.993327980707563	0.1010741675774\\
0.999	0.0774816382160537\\
};
\addlegendentry{BMI synthesis: n = 6};

\addplot [color=mycolor2,dashed,mark=*,mark options={solid,fill=white!30!mycolor2}]
  table[row sep=crcr]{%
0.96078431372549	0.141421363792861\\
0.964313725490196	0.134906000978448\\
0.967843137254902	0.128063185732074\\
0.971372549019608	0.120830271910446\\
0.974901960784314	0.113142313612656\\
0.97843137254902	0.104886031800763\\
0.981960784313725	0.0959196089862665\\
0.985490196078431	0.0860280570309588\\
0.989019607843137	0.0748394646391178\\
0.992549019607843	0.061649513407354\\
0.996078431372549	0.0447336084716429\\
0.999607843137255	0.0142635363359695\\
};
\addlegendentry{Convex synthesis: n = 2};

\addplot [color=mycolor3,only marks,mark=*,mark options={solid,fill=mycolor3},mark size=3.5pt]
  table[row sep=crcr]{%
0.858578643762691	0.350239408530843\\
};
\addlegendentry{Triple Momentum};

\addplot [color=mycolor4,only marks,mark=*,mark options={solid,fill=mycolor4},mark size=3.5pt]
  table[row sep=crcr]{%
0.89498743710662	0.210519703384064\\
};
\addlegendentry{Nesterov's Method};

\addplot [color=mycolor5,only marks,mark=*,mark options={solid,fill=mycolor5},mark size=3.5pt]
  table[row sep=crcr]{%
0.96078431372549	0.141419971567621\\
};
\addlegendentry{Gradient Descent};

\end{axis}
\end{tikzpicture}%
%
%
\definecolor{mycolor1}{rgb}{0.00000,0.39062,0.00000}%
\definecolor{mycolor2}{rgb}{0.29297,0.00000,0.50781}%
\definecolor{mycolor3}{rgb}{0.85000,0.32500,0.09800}%
\definecolor{mycolor4}{rgb}{0.00000,0.44700,0.74100}%
\definecolor{mycolor5}{rgb}{0.92900,0.69400,0.12500}%
\begin{tikzpicture}

\begin{axis}[%
width=0.951\figurewidth,
height=\figureheight,
at={(0\figurewidth,0\figureheight)},
scale only axis,
xmin=0.85,
xmax=1,
xlabel={Convergence rate $\rho$},
ymin=0.005,
ymax=0.4,
ymode=log,
ylabel={$H_2$-performance},
axis background/.style={fill=white},
legend style={at={(0.04,0.04)},anchor=south west,legend cell align=left,align=left,draw=white!15!black,font=\small},
grid=both,
grid style={line width=.1pt, draw=black!30},
major grid style={line width=.2pt,draw=black!80}
]
\addplot [color=mycolor1,dashed,mark=*,mark options={solid,fill=mycolor1}]
  table[row sep=crcr]{%
0.9045	0.261567264083616\\
0.9084375	0.243555153635238\\
0.912375	0.231129695842992\\
0.9163125	0.215287805386304\\
0.92025	0.204588889255119\\
0.9241875	0.198329490965216\\
0.928125	0.191266391487824\\
0.9320625	0.185835050658311\\
0.936	0.181813050688669\\
0.9399375	0.177015932225271\\
0.943875	0.169426947416612\\
0.9478125	0.166620505054766\\
0.95175	0.16362239046844\\
0.9556875	0.159687888968721\\
0.959625	0.151980294226251\\
0.9635625	0.143500869773141\\
0.9675	0.135607304479228\\
0.9714375	0.12398282818492\\
0.975375	0.111812537181058\\
0.9793125	0.101853499719628\\
0.98325	0.0927881858801923\\
0.9871875	0.0833378666209723\\
0.991125	0.0751994950863044\\
0.9950625	0.070845975044466\\
0.999	0.0701835359047382\\
};
\addlegendentry{BMI synthesis: n = 2};

\addplot [color=mycolor1,dashed,mark=*,mark options={solid,fill=white!50!mycolor1}]
  table[row sep=crcr]{%
0.9045	0.272405601077845\\
0.9084375	0.248368929085234\\
0.912375	0.228724028675297\\
0.9163125	0.215874028821559\\
0.92025	0.202784244802335\\
0.9241875	0.191177546250203\\
0.928125	0.183860790345129\\
0.9320625	0.17255063338454\\
0.936	0.165202650504344\\
0.9399375	0.160965808779423\\
0.943875	0.151504150535386\\
0.9478125	0.147702225697656\\
0.95175	0.144422527223611\\
0.9556875	0.139941439610919\\
0.959625	0.13479136467685\\
0.9635625	0.129441724767907\\
0.9675	0.125878477270222\\
0.9714375	0.118994571725512\\
0.975375	0.113565809605556\\
0.9793125	0.10844956157998\\
0.98325	0.100089219577991\\
0.9871875	0.0907196323107382\\
0.991125	0.0852976125711656\\
0.9950625	0.070742114149907\\
0.999	0.0529807306114738\\
};
\addlegendentry{BMI synthesis: n = 6};

\addplot [color=mycolor2,dashed,mark=*,mark options={solid,fill=white!30!mycolor2}]
  table[row sep=crcr]{%
0.98019801980198	0.100000019042179\\
0.981980198019802	0.0953754951992458\\
0.983762376237624	0.0905416855625699\\
0.985544554455446	0.085427655425053\\
0.987326732673267	0.0799903695659748\\
0.989108910891089	0.0741565418750907\\
0.990891089108911	0.0678132192084071\\
0.992673267326733	0.0608288887500978\\
0.994455445544554	0.0529196855675857\\
0.996237623762376	0.04362275140369\\
0.998019801980198	0.0316502773980378\\
0.99980198019802	0.0102934154763341\\
};
\addlegendentry{Convex synthesis: n = 2};

\addplot [color=mycolor3,only marks,mark=*,mark options={solid,fill=mycolor3},mark size=3.5pt]
  table[row sep=crcr]{%
0.9	0.336189617584283\\
};
\addlegendentry{Triple Momentum};

\addplot [color=mycolor4,only marks,mark=*,mark options={solid,fill=mycolor4},mark size=3.5pt]
  table[row sep=crcr]{%
0.926786189044345	0.181669447045039\\
};
\addlegendentry{Nesterov's Method};

\addplot [color=mycolor5,only marks,mark=*,mark options={solid,fill=mycolor5},mark size=3.5pt]
  table[row sep=crcr]{%
0.98019801980198	0.0999998445393456\\
};
\addlegendentry{Gradient Descent};

\end{axis}
\end{tikzpicture}%
	\end{footnotesize}
	\end{center}
	\caption{Guaranteed convergence rates and corresponding $H_2$-performance levels for $m=1$ and $L=50$ (left) as well as $L=100$ (right).
	For comparison, we also plot the $H_2$-performance levels {for the Triple Momentum Method, Nesterov's Method and the Gradient Descent algorithm obtained from analysis using~\Cref{lemmaAnalysisConvRate}, \Cref{lemmaAnalysisH2}.}}
	}{
	\setlength\figurewidth{0.39\textwidth}
	\setlength\figureheight{0.35\textwidth}
	\begin{center}
%
%
\definecolor{mycolor1}{rgb}{0.00000,0.39062,0.00000}%
\definecolor{mycolor2}{rgb}{0.29297,0.00000,0.50781}%
\definecolor{mycolor3}{rgb}{0.85000,0.32500,0.09800}%
\definecolor{mycolor4}{rgb}{0.00000,0.44700,0.74100}%
\definecolor{mycolor5}{rgb}{0.92900,0.69400,0.12500}%
\begin{tikzpicture}

\begin{axis}[%
width=0.951\figurewidth,
height=\figureheight,
at={(0\figurewidth,0\figureheight)},
scale only axis,
unbounded coords=jump,
xmin=0.85,
xmax=1,
xlabel={Convergence rate $\rho$},
ymin=0.005,
ymax=0.4,
ymode=log,
ylabel={$H_2$-performance},
axis background/.style={fill=white},
legend style={at={(0.04,0.04)},anchor=south west,legend cell align=left,align=left,draw=white!15!black,font=\small},
grid=both,
grid style={line width=.1pt, draw=black!30},
major grid style={line width=.2pt,draw=black!80},
ytick={0.01,0.1}
]
\addplot [color=mycolor1,dashed,mark=*,mark options={solid,fill=mycolor1}]
  table[row sep=crcr]{%
0.862871536981504	0.30946152453036\\
0.868543556273941	0.288298636481586\\
0.874215575566379	0.267328992115637\\
0.879887594858816	0.251804914556189\\
0.885559614151253	0.239264453926221\\
0.891231633443691	0.230883435261668\\
0.896903652736128	0.224452274367907\\
0.902575672028565	0.218636272483777\\
0.908247691321003	0.213041078515733\\
0.91391971061344	0.208782422389828\\
0.919591729905877	0.204608588573787\\
0.925263749198315	0.197687811804989\\
0.930935768490752	0.192445233553485\\
0.936607787783189	0.184933986506113\\
0.942279807075627	0.178715952597561\\
0.947951826368064	0.166866285282643\\
0.953623845660501	0.15461276337932\\
0.959295864952939	0.14558710996721\\
0.964967884245376	0.134695296588531\\
0.970639903537813	0.124890226821288\\
0.976311922830251	0.116701440957512\\
0.981983942122688	0.113124350381597\\
0.987655961415125	0.102805955772367\\
0.993327980707563	0.101259898070055\\
0.999	0.101675364929887\\
};
\addlegendentry{BMI synthesis: n = 2};

\addplot [color=mycolor1,dashed,mark=*,mark options={solid,fill=white!50!mycolor1}]
  table[row sep=crcr]{%
0.862871536981504	nan\\
0.868543556273941	0.288481688089103\\
0.874215575566379	0.266947234996696\\
0.879887594858816	0.251914680311602\\
0.885559614151253	0.236312288366904\\
0.891231633443691	0.224751520970588\\
0.896903652736128	0.2148201122362\\
0.902575672028565	0.202914257308799\\
0.908247691321003	0.195601319561357\\
0.91391971061344	0.18705395991594\\
0.919591729905877	0.185284101951933\\
0.925263749198315	0.177876477563977\\
0.930935768490752	0.168862293049972\\
0.936607787783189	0.162758399794006\\
0.942279807075627	0.158858866727713\\
0.947951826368064	0.15561693070124\\
0.953623845660501	0.149971693472256\\
0.959295864952939	0.144353607176755\\
0.964967884245376	0.137964463567028\\
0.970639903537813	0.133212749915167\\
0.976311922830251	0.122083768435539\\
0.981983942122688	0.116669364452387\\
0.987655961415125	0.100943799672138\\
0.993327980707563	0.1010741675774\\
0.999	0.0774816382160537\\
};
\addlegendentry{BMI synthesis: n = 6};

\addplot [color=mycolor2,dashed,mark=*,mark options={solid,fill=white!30!mycolor2}]
  table[row sep=crcr]{%
0.96078431372549	0.141421363792861\\
0.964313725490196	0.134906000978448\\
0.967843137254902	0.128063185732074\\
0.971372549019608	0.120830271910446\\
0.974901960784314	0.113142313612656\\
0.97843137254902	0.104886031800763\\
0.981960784313725	0.0959196089862665\\
0.985490196078431	0.0860280570309588\\
0.989019607843137	0.0748394646391178\\
0.992549019607843	0.061649513407354\\
0.996078431372549	0.0447336084716429\\
0.999607843137255	0.0142635363359695\\
};
\addlegendentry{Convex synthesis: n = 2};

\addplot [color=mycolor3,only marks,mark=*,mark options={solid,fill=mycolor3},mark size=3.5pt]
  table[row sep=crcr]{%
0.858578643762691	0.350239408530843\\
};
\addlegendentry{Triple Momentum};

\addplot [color=mycolor4,only marks,mark=*,mark options={solid,fill=mycolor4},mark size=3.5pt]
  table[row sep=crcr]{%
0.89498743710662	0.210519703384064\\
};
\addlegendentry{Nesterov's Method};

\addplot [color=mycolor5,only marks,mark=*,mark options={solid,fill=mycolor5},mark size=3.5pt]
  table[row sep=crcr]{%
0.96078431372549	0.141419971567621\\
};
\addlegendentry{Gradient Descent};

\end{axis}
\end{tikzpicture}%
%
%
\definecolor{mycolor1}{rgb}{0.00000,0.39062,0.00000}%
\definecolor{mycolor2}{rgb}{0.29297,0.00000,0.50781}%
\definecolor{mycolor3}{rgb}{0.85000,0.32500,0.09800}%
\definecolor{mycolor4}{rgb}{0.00000,0.44700,0.74100}%
\definecolor{mycolor5}{rgb}{0.92900,0.69400,0.12500}%
\begin{tikzpicture}

\begin{axis}[%
width=0.951\figurewidth,
height=\figureheight,
at={(0\figurewidth,0\figureheight)},
scale only axis,
xmin=0.85,
xmax=1,
xlabel={Convergence rate $\rho$},
ymin=0.005,
ymax=0.4,
ymode=log,
ylabel={$H_2$-performance},
axis background/.style={fill=white},
legend style={at={(0.04,0.04)},anchor=south west,legend cell align=left,align=left,draw=white!15!black,font=\small},
grid=both,
grid style={line width=.1pt, draw=black!30},
major grid style={line width=.2pt,draw=black!80}
]
\addplot [color=mycolor1,dashed,mark=*,mark options={solid,fill=mycolor1}]
  table[row sep=crcr]{%
0.9045	0.261567264083616\\
0.9084375	0.243555153635238\\
0.912375	0.231129695842992\\
0.9163125	0.215287805386304\\
0.92025	0.204588889255119\\
0.9241875	0.198329490965216\\
0.928125	0.191266391487824\\
0.9320625	0.185835050658311\\
0.936	0.181813050688669\\
0.9399375	0.177015932225271\\
0.943875	0.169426947416612\\
0.9478125	0.166620505054766\\
0.95175	0.16362239046844\\
0.9556875	0.159687888968721\\
0.959625	0.151980294226251\\
0.9635625	0.143500869773141\\
0.9675	0.135607304479228\\
0.9714375	0.12398282818492\\
0.975375	0.111812537181058\\
0.9793125	0.101853499719628\\
0.98325	0.0927881858801923\\
0.9871875	0.0833378666209723\\
0.991125	0.0751994950863044\\
0.9950625	0.070845975044466\\
0.999	0.0701835359047382\\
};
\addlegendentry{BMI synthesis: n = 2};

\addplot [color=mycolor1,dashed,mark=*,mark options={solid,fill=white!50!mycolor1}]
  table[row sep=crcr]{%
0.9045	0.272405601077845\\
0.9084375	0.248368929085234\\
0.912375	0.228724028675297\\
0.9163125	0.215874028821559\\
0.92025	0.202784244802335\\
0.9241875	0.191177546250203\\
0.928125	0.183860790345129\\
0.9320625	0.17255063338454\\
0.936	0.165202650504344\\
0.9399375	0.160965808779423\\
0.943875	0.151504150535386\\
0.9478125	0.147702225697656\\
0.95175	0.144422527223611\\
0.9556875	0.139941439610919\\
0.959625	0.13479136467685\\
0.9635625	0.129441724767907\\
0.9675	0.125878477270222\\
0.9714375	0.118994571725512\\
0.975375	0.113565809605556\\
0.9793125	0.10844956157998\\
0.98325	0.100089219577991\\
0.9871875	0.0907196323107382\\
0.991125	0.0852976125711656\\
0.9950625	0.070742114149907\\
0.999	0.0529807306114738\\
};
\addlegendentry{BMI synthesis: n = 6};

\addplot [color=mycolor2,dashed,mark=*,mark options={solid,fill=white!30!mycolor2}]
  table[row sep=crcr]{%
0.98019801980198	0.100000019042179\\
0.981980198019802	0.0953754951992458\\
0.983762376237624	0.0905416855625699\\
0.985544554455446	0.085427655425053\\
0.987326732673267	0.0799903695659748\\
0.989108910891089	0.0741565418750907\\
0.990891089108911	0.0678132192084071\\
0.992673267326733	0.0608288887500978\\
0.994455445544554	0.0529196855675857\\
0.996237623762376	0.04362275140369\\
0.998019801980198	0.0316502773980378\\
0.99980198019802	0.0102934154763341\\
};
\addlegendentry{Convex synthesis: n = 2};

\addplot [color=mycolor3,only marks,mark=*,mark options={solid,fill=mycolor3},mark size=3.5pt]
  table[row sep=crcr]{%
0.9	0.336189617584283\\
};
\addlegendentry{Triple Momentum};

\addplot [color=mycolor4,only marks,mark=*,mark options={solid,fill=mycolor4},mark size=3.5pt]
  table[row sep=crcr]{%
0.926786189044345	0.181669447045039\\
};
\addlegendentry{Nesterov's Method};

\addplot [color=mycolor5,only marks,mark=*,mark options={solid,fill=mycolor5},mark size=3.5pt]
  table[row sep=crcr]{%
0.98019801980198	0.0999998445393456\\
};
\addlegendentry{Gradient Descent};

\end{axis}
\end{tikzpicture}%
	\end{center}
	\caption{Guaranteed convergence rates and corresponding $H_2$-performance levels for $m=1$ and $L=50$ (top) as well as $L=100$ (bottom).
	For comparison, we also plot the $H_2$-performance levels {for the Triple Momentum Method, Nesterov's Method and the Gradient Descent algorithm obtained from analysis using~\Cref{lemmaAnalysisConvRate}, \Cref{lemmaAnalysisH2}.}}}
	\label{figSynthAlgoAdditiveNoise}
\end{figure}

\section{Structure exploiting algorithms}\label{secStructureExploitingAlgorithms}
Up to now the only assumption on the objective function was that $ \fObj \in {\classObj{m}{L}} $.
However, in many situations it is possible to further characterize the objective
function in terms of its structure. In such cases, it is to be expected that a tailored algorithm that 
exploits these properties has much better guaranteed convergence rates
{than} a standard algorithm for the class $ {\classObj{m}{L}} $. 
In this section we elaborate on how additional structural
knowledge can be incorporated in the presented framework.
To this end, for any{\sout{$ q \in \nat $ and}} $ L \geq m > 0 $,
we define the following two sets of functions: 
\ifthenelse{\boolean{singleColumn}}{
\begin{align}
	\begin{split}
	{{\classObj{m}{L}^{\nonrep}}} = \big\lbrace {\fObj \in {\classObj{m}{L}} }  \sepSet~& \exists \phi_i : \real \to \real, 
	\phi_i \in \mathcal{C}^0, i=1,2,\dots,p, \text{ such that } \\
	&\nabla \fObj(x) = \begin{bmatrix} \phi_1(x_1) & \dots & \phi_p(x_p) \end{bmatrix}^\top \\
	&\text{for all } x = \begin{bmatrix} x_1 & x_2 & \dots & x_p \end{bmatrix}^\top \in \real^p \big\rbrace
	\end{split} \label{eqDefFunctionDiagonalHessian} \\
	\begin{split}
	{{\classObj{m}{L}^{\rep}}} = \big\lbrace {\fObj \in {\classObj{m}{L}}} \sepSet~& \exists \phi : \real \to \real, 
	\phi \in \mathcal{C}^0, \text{ such that } \\
	& \nabla \fObj(x) = \begin{bmatrix} \phi(x_1) & \dots & \phi(x_p) \end{bmatrix}^\top \\
	& \text{for all } x = \begin{bmatrix} x_1 & x_2 & \dots & x_p \end{bmatrix}^\top \in \real^p \big\rbrace.
	\end{split} \label{eqDefFunctionDiagonalRepeatedHessian}
\end{align}}{
\begin{align}
	\begin{split}
	{{\classObj{m}{L}^{\nonrep}}} &= \big\lbrace {\fObj \in {\classObj{m}{L}}}  \sepSet \exists \phi_i : \real \to \real, \\
	&\hphantom{= \big\lbrace}\phi_i \in \mathcal{C}^0, i=1,2,\dots,p, \text{ such that } \\
	&\hphantom{= \big\lbrace}\nabla \fObj(x) = \begin{bmatrix} \phi_1(x_1) & \dots & \phi_p(x_p) \end{bmatrix}^\top \\
	&\hphantom{= \big\lbrace}\text{for all } x = \begin{bmatrix} x_1 & x_2 & \dots & x_p \end{bmatrix}^\top \in \real^p \big\rbrace
	\end{split} \label{eqDefFunctionDiagonalHessian} \\
	\begin{split}
	{{\classObj{m}{L}^{\rep}}} &= \big\lbrace {\fObj \in {\classObj{m}{L}}} \sepSet \exists \phi : \real \to \real, \\
	&\hphantom{= \big\lbrace}\phi \in \mathcal{C}^0, \text{ such that } \\
	&\hphantom{= \big\lbrace} \nabla \fObj(x) = \begin{bmatrix} \phi(x_1) & \dots & \phi(x_p) \end{bmatrix}^\top \\
	&\hphantom{= \big\lbrace}\text{for all } x = \begin{bmatrix} x_1 & x_2 & \dots & x_p \end{bmatrix}^\top \in \real^p \big\rbrace.
	\end{split} \label{eqDefFunctionDiagonalRepeatedHessian}
\end{align}}
Note that {$ {\classObj{m}{L}^{\nonrep}} $}
is the set of all functions with a diagonal Hessian and {$ {\classObj{m}{L}^{\rep}} $} is the set 
of all functions with a diagonally repeated Hessian. Objective functions
of this form appear, e.g., in distributed optimization problems where
the objective function is a sum of the agents' individual objective functions.
We discuss how this information
can be translated to IQCs in~\Cref{subsecIQCsForStructuredObjectiveFunctions}. 
{It is straightforward to apply the same procedure to the 
case where each $ \phi_i $, $ i = 1,2,\dots,p $, in~\eqref{eqDefFunctionDiagonalHessian}
is a map from $ \real^{q_i} $ to $ \real^{q_i} $, $ q_i \geq 1 $ (or
$ \phi: \real^q \to \real^q $ in~\eqref{eqDefFunctionDiagonalRepeatedHessian}, respectively).
However, this will yield the same matrix inequalities for 
algorithm analysis and synthesis after a dimension reduction;
thus we directly consider the reduced case here. 
In other words, the resulting matrix inequalities are scalable
with respect to the dimensions $q_i$ (or $q$, respectively).} 
Additionally, in~\Cref{subsecParametrizedObjectiveFunctions}, we consider 
{parametrized objective functions consisting} of a known quadratic and a specifically structured
unknown part.

\subsection{IQCs for {the classes $ \unboldmath{\classObj{m}{L}^{\normalfont\nonrep}} $, $ \unboldmath{\classObj{m}{L}^{\normalfont\rep}} $}}\label{subsecIQCsForStructuredObjectiveFunctions}
In the following we want to derive IQCs for uncertainties of the form~\eqref{eqDefDelta}
and the transformed version~\eqref{eqDefDeltaRho} thereof under the additional assumption
that the objective function $H$ is in {$ {\classObj{m}{L}^{\nonrep}} $ or $ {\classObj{m}{L}^{\rep}} $}, respectively.
{At this point, we benefit from embedding our approach in the standard IQC framework
rendering the following derivation straightforward.}
{In the spirit of our previous discussions,} we define the {following} sets of uncertainties 
\begin{subequations}
\begin{align}
	\deltaNonrep(m,L) &= \big\lbrace \Delta_\fObj \sepSet \fObj \in {{\classObj{m}{L}^{\nonrep}}} \big\rbrace, \label{eqDefSetOfDeltaNonRepeated} \\
	\deltaRep(m,L) &= \big\lbrace \Delta_\fObj \sepSet \fObj \in {{\classObj{m}{L}^{\rep}}} \big\rbrace, \label{eqSetOfDeltaRepeated}
\end{align}
\end{subequations}
with $ \Delta_\fObj(y)\k := \nabla \fObj ( y_k + \optVar^\optSign ) - m y_k $, see~\eqref{eqDefDelta}.
Observe that $ \deltaRep(m,L) \subseteq \deltaNonrep(m,L) \subseteq \deltaUnstr(m,L) $.
Operators of this type have as well been studied in the literature.
In fact, the following modification of~\Cref{lemmaZFGeneralTimeDomain} holds, see~\citet{fetzer2017absolute}.
{
\begin{lemma}\label{lemmaZFstructuredTimeDomain}
Let $ L \geq m > 0 $ and $ \dimCausal, \dimAnticausal \in \nat $ be given.
	Let $ \deltaNonrep(m,L), \deltaRep(m,L) $ be defined according to~\eqref{eqDefSetOfDeltaNonRepeated}, \eqref{eqSetOfDeltaRepeated} and 
	let $ {M}: \ltwo^p \to \ltwo^p $ be a bounded linear operator defined as
	\begin{align}
	{M} = \toep\big(\begin{bmatrix} M_{-\dimCausal} & M_{-\dimCausal+1} & \dots & M_{\dimAnticausal} \end{bmatrix} \big), 	
	\end{align}
	where $ M_i \in \real^{p \times p} $.
	Then~\eqref{eqIneqZamesFalbTimeDomain} persists to hold 
	\begin{enumerate}[leftmargin=*,label={(\alph*)}]
		\item for all $ \Delta \in \deltaNonrep(m,L) $ if, for $ i \in \lbrace -\dimCausal, \dots, \dimAnticausal \rbrace $,
			\begin{align}
				M_i = \textup{blkdiag}( m_{i,1}, m_{i,2}, \dots, m_{i,p} )
			\end{align}
			and, for each $ k = 1,2,\dots,p $, the operator
			\begin{align}
				\toep( \left[\begin{shortArray}{ccccc} m_{-\dimCausal,k} & \dots & m_{0,k} & \dots & m_{\dimAnticausal,k} \end{shortArray}\right] )
			\end{align}
			is doubly hyperdominant.
		\item for all $ \Delta \in \deltaRep(m,L) $ if $M$ is doubly hyperdominant. \qedhere
	\end{enumerate}
\end{lemma}}
By our discussions {in} \Cref{subsecExponentialConvergenceViaIQCs}, it is {not difficult}
to derive IQCs for exponential stability analysis. We first define the 
corresponding sets of transformed uncertainties as
\begin{subequations}
\begin{align}
	\mathbf{\Delta}_{\nonrep,\rho}(m,L) &= \lbrace \rho_- \circ \Delta \circ \rho_+ \sepSet \Delta \in \mathbf{\Delta}_\nonrep(m,L) \rbrace \label{eqDefSetOfDeltaRhoNonrepeated} \\
	\mathbf{\Delta}_{\rep,\rho}(m,L) &= \lbrace \rho_- \circ \Delta \circ \rho_+ \sepSet \Delta \in \mathbf{\Delta}_\rep(m,L) \rbrace. \label{eqDefSetOfDeltaRhoRepeated}
\end{align}
\end{subequations}
Following the same steps, again we demand $ M_\rho = \rho_- \circ M \circ \rho_- $ 
to adhere to {condition~(a)} in~\Cref{lemmaZFstructuredTimeDomain} if $ \Delta_\rho \in \mathbf{\Delta}_{\nonrep,\rho}(m,L) $
or {to~(b)} if $ \Delta_\rho \in \mathbf{\Delta}_{\rep,\rho}(m,L) $.
The sets of admissible $ [ M_i ]_{i \in \lbrace-\dimCausal,\dots,\dimAnticausal\rbrace} $ are then defined as
\ifthenelse{\boolean{singleColumn}}{
\begin{align}
	\mathbb{M}_\rep(\rho,\dimAnticausal,\dimCausal,p) = \big\lbrace \left[\begin{shortArray}{cccc} M_{-\dimCausal} & M_{-\dimCausal +1} & \dots & M_{\dimAnticausal} \end{shortArray}\right] \sepSet 
	~& M_i \in \real^{p \times p},  M_i \leq 0 \text{ for all } i \neq 0, \nonumber \\ 
	& ( \textstyle{\sum_{i=-\dimCausal}^{\dimAnticausal} M_i \rho^{-i}} ) \mathbf{1} \geq 0, \nonumber \\
	& \mathbf{1}^\top ( \textstyle{\sum_{i=-\dimCausal}^{\dimAnticausal} M_i \rho^{i}} ) \geq 0\big\rbrace 
	\label{eqDefSetMrepeated} 
\end{align}
}{
\begin{align}
	\mathbb{M}_\rep(\rho,\dimAnticausal,\dimCausal,p) =~& \big\lbrace \left[\begin{shortArray}{cccc} M_{-\dimCausal} & M_{-\dimCausal +1} & \dots & M_{\dimAnticausal} \end{shortArray}\right] \sepSet \nonumber \\
	& M_i \in \real^{p \times p},  M_i \leq 0 \text{ for all } i \neq 0, \nonumber \\ 
	& ( \textstyle{\sum_{i=-\dimCausal}^{\dimAnticausal} M_i \rho^{-i}} ) \mathbf{1} \geq 0, \nonumber \\
	& \mathbf{1}^\top ( \textstyle{\sum_{i=-\dimCausal}^{\dimAnticausal} M_i \rho^{i}} ) \geq 0\big\rbrace 
	\label{eqDefSetMrepeated} 
\end{align}}
in the repeated case and as
\ifthenelse{\boolean{singleColumn}}{
\begin{align}
 	\mathbb{M}_\nonrep(\rho,\dimAnticausal,\dimCausal,p) =
 	\big\lbrace \left[\begin{shortArray}{cccc} M_{-\dimCausal} & M_{-\dimCausal +1} & \dots & M_{\dimAnticausal} \end{shortArray}\right] \sepSet
	~&M_i = \textup{blkdiag}( m_{i,1}, m_{i,2}, \dots, m_{i,p} ), \nonumber \\
	&M_i \leq 0 \text{ for all } i \neq 0, \nonumber \\
	& ( \textstyle{\sum_{i=-\dimCausal}^{\dimAnticausal} M_i \rho^{-i}} ) \mathbf{1} \geq 0, \nonumber \\
	& \mathbf{1}^\top ( \textstyle{\sum_{i=-\dimCausal}^{\dimAnticausal} M_i \rho^{i}} ) \geq 0\big\rbrace \label{eqDefSetMnonRepeated}
\end{align}}{
\begin{align}
 	\mathbb{M}_\nonrep(\rho,\dimAnticausal,\dimCausal,p) =~& 
 	\big\lbrace \left[\begin{shortArray}{cccc} M_{-\dimCausal} & M_{-\dimCausal +1} & \dots & M_{\dimAnticausal} \end{shortArray}\right] \sepSet \nonumber \\
	&M_i = \textup{blkdiag}( m_{i,1}, m_{i,2}, \dots, m_{i,p} ), \nonumber \\
	&M_i \leq 0 \text{ for all } i \neq 0, \nonumber \\
	& ( \textstyle{\sum_{i=-\dimCausal}^{\dimAnticausal} M_i \rho^{-i}} ) \mathbf{1} \geq 0, \nonumber \\
	& \mathbf{1}^\top ( \textstyle{\sum_{i=-\dimCausal}^{\dimAnticausal} M_i \rho^{i}} ) \geq 0\big\rbrace \label{eqDefSetMnonRepeated}
\end{align}}
for the non-repeated case. The corresponding Zames-Falb multipliers $ \mathbf{\Pi}_{\Delta,\rho,\nonrep}^p(m,L) $,
$ \mathbf{\Pi}_{\Delta,\rho,\rep}^p(m,L) $ are obtained
from~\eqref{eqDefZamesFalbMultipliers} simply by replacing $ \mathbb{M} $
by $ \mathbb{M}_\nonrep(\rho,\dimAnticausal,\dimCausal,p) $ and $ \mathbb{M}_\rep(\rho,\dimAnticausal,\dimCausal,p) $,
respectively.
Analogously to~\Cref{lemmaZFIQCrho}, we then have the following result.
\begin{theorem}\label{lemmaZFIQCrhoStructured}
Let $ L \geq m > 0 $ and let $ \mathbf{\Delta}_{\nonrep,\rho}(m,L) $, $\mathbf{\Delta}_{\rep,\rho}(m,L) $
	be defined as in~\eqref{eqDefSetOfDeltaRhoNonrepeated}, \eqref{eqDefSetOfDeltaRhoRepeated}.
	Then, for each $ \rho \in (0,1] $, $ \Delta_\rho $ satisfies the IQC defined by $ \Pi $
	\begin{enumerate}[leftmargin=*]
		\item for each $ \Delta_\rho \in \mathbf{\Delta}_{\nonrep,\rho}(m,L) $ and each $ \Pi \in \mathbf{\Pi}_{\Delta,\rho,\nonrep}^p(m,L) $;
		\item for each $ \Delta_\rho \in \mathbf{\Delta}_{\rep,\rho}(m,L) $ and each $ \Pi \in \mathbf{\Pi}_{\Delta,\rho,\rep}^p(m,L) $. \qedhere
	\end{enumerate}
\end{theorem}
It is then straightforward to include additional structural assumptions
in the results from~\Cref{secRobustOptimizationAlgorithms}. In a nutshell,
we only need to replace the constraint on $ \begin{bmatrix} M_- & M_0 & M_+ \end{bmatrix} $
by the respective structured counterpart.

\subsection{{Parametrized objective functions}}\label{subsecParametrizedObjectiveFunctions}
In the following we assume that the gradient of the objective function $ \fObj : \real^p \to \real $
admits the form
\begin{align}
	\nabla \fObj(\optVar) = \fObj_1 \optVar +  \trafoStructured^\top \nabla \fObj_2 (\trafoStructured \optVar), \label{eqExampleStructureGradient}
\end{align}
where $ \fObj_1 \in \real^{p \times p} $, $ \fObj_1 \succ 0 $, $ \trafoStructured \in \real^{q \times p } $
are known whereas $ \fObj_2 : \real^q \to \real $ is unknown but 
fulfills $ \fObj_2 \in {\classObj{m_2}{L_2}} $ for some known constants
$ L_2 \geq m_2 > 0 $. 
Objective functions of this specific form arise in the context of linear
relaxed logarithmic barrier function based model predictive control~\citep{feller2017relaxed},
where fast converging optimization algorithms are crucial for the practical
applicability of the control scheme.
It is clear that $ \fObj \in {\classObj{m}{L}} $
with $ m = \lambda_{\textup{min}}( \fObj_1 ) $, $ L = \lambda_{\textup{max}}( \fObj_1 + \trafoStructured^\top \trafoStructured L_2 ) $,
where $ \lambda_{\textup{min}} $, $ \lambda_{\textup{max}} $ denote the minimal
and maximal eigenvalue, respectively. {Again, we denote by $ \optVar^\optSign $
the minimizer of $ \fObj $.}
Note that $ m, L $ can both be computed under the assumption that $ \fObj_1 $
and $\trafoStructured$ are known. Hence, we can make use of any algorithm for the class
$ {\classObj{m}{L}} $; in particular, we can make use of the fastest known
algorithm in the class of considered algorithms, i.e., the \tripleMomentum~\citep{vanScoy2018fastest}.
However, with the structure of $ \fObj $ being
known and taking the form~\eqref{eqExampleStructureGradient}, 
it is possible to obtain algorithms with improved convergence rate
guarantees in many cases using the presented framework.

To this end, consider an algorithm of the form~\eqref{eqOptAlgo} under the 
structural assumption~\eqref{eqExampleStructureGradient}
which then takes the form
\begin{subequations}
\begin{align}
	\xOpt\kk &= ( \Aopt + \Bopt \fObj_1 \Copt ) \xOpt\k + \Bopt \trafoStructured^\top \nabla \fObj_2 ( \trafoStructured \Copt \xOpt\k ) \\
	{\optVar\k} &= {\Dopt \xOpt\k}.
\end{align}
\end{subequations}
By the same state transformation $ \xTrafo\k = \xOpt\k - \Dopt^\dagger \optVar^\optSign $ and
under the assumptions~\eqref{eqConditionABCD} we then 
obtain the transformed dynamics
\begin{subequations}
\ifthenelse{\boolean{singleColumn}}{
\begin{align}
	\xTrafo\kk =~& ( \Aopt + \Bopt \fObj_1 \Copt ) \xTrafo\k + \Bopt \trafoStructured^\top \nabla \fObj_2( \trafoStructured \Copt \xTrafo\k + \trafoStructured \optVar^\optSign) + \Bopt \fObj_1 \optVar^\optSign \\
	{\optVar\k} =~& {\Dopt \xTrafo\k + \optVar^\optSign}.
\end{align}}{
\begin{align}
	\xTrafo\kk =~& ( \Aopt + \Bopt \fObj_1 \Copt ) \xTrafo\k + \Bopt \trafoStructured^\top \nabla \fObj_2( \trafoStructured \Copt \xTrafo\k + \trafoStructured \optVar^\optSign) \\
	&+ \Bopt \fObj_1 \optVar^\optSign \nonumber \\
	{\optVar\k} =~& {\Dopt \xTrafo\k + \optVar^\optSign}.
\end{align}}\label{eqTransformedDynamicsStructured}
\end{subequations}
{Since $ \optVar^\optSign $ is the minimizer of $ \fObj $}, we
have $ \nabla \fObj(\optVar^\optSign) = \fObj_1 \optVar^\optSign +  \trafoStructured^\top \nabla \fObj_2 (\trafoStructured \optVar^\optSign) = 0 $.
Since $ \fObj_1 \succeq {m} I \succ 0 $, $ \fObj_1 $ is non-singular and hence
the previous {equality} implies that $ \optVar^\optSign = - \fObj_1^{-1} \trafoStructured^\top \nabla \fObj_2 (\trafoStructured \optVar^\optSign) $.
Using this in~\eqref{eqTransformedDynamicsStructured}, we obtain
\begin{subequations}
\ifthenelse{\boolean{singleColumn}}{
\begin{align}
	\xTrafo\kk =~& ( \Aopt + \Bopt \fObj_1 \Copt ) \xTrafo\k 
	+ \Bopt \trafoStructured^\top \big( \nabla \fObj_2( \trafoStructured \Copt \xTrafo\k + \trafoStructured \optVar^\optSign) - \nabla \fObj_2 (\trafoStructured \optVar^\optSign) \big)  \\
	{\optVar\k =}~&{\Dopt \xTrafo\k + \optVar^\optSign}.
\end{align}}{
\begin{align}
	\xTrafo\kk =~& ( \Aopt + \Bopt \fObj_1 \Copt ) \xTrafo\k \\
	&+ \Bopt \trafoStructured^\top \big( \nabla \fObj_2( \trafoStructured \Copt \xTrafo\k + \trafoStructured \optVar^\optSign) - \nabla \fObj_2 (\trafoStructured \optVar^\optSign) \big) \nonumber  \\
	{\optVar\k =}~&{\Dopt \xTrafo\k + \optVar^\optSign}.
\end{align}}\label{eqTransformedDynamicsStructured2}
\end{subequations}
In the same spirit as in the unstructured case, we hence write 
the feedback interconnection in the standard form as
\begin{subequations}
\begin{align}
	\xTrafo\kk &= ( \Aopt + \Bopt \fObj_1 \Copt + m_2 \Bopt \trafoStructured^\top \trafoStructured \Copt ) \xTrafo\k + w\k \\
	y\k &= \trafoStructured \Copt \xTrafo\k \\
	w\k &= \Delta_{\fObj_2}(y)\k
\end{align}\label{eqOptAlgoStandardFormStructured}
\end{subequations}
with the new uncertainty $ \Delta_{\fObj_2}: \ltwo^q \to \ltwo^q $, defined as
\begin{align}
	\Delta_{\fObj_2}(y)\k = \nabla\fObj_2 ( y\k + \trafoStructured \optVar^\optSign ) - \nabla\fObj_2 ( \trafoStructured \optVar^\optSign ) - m_2 y_k
\end{align}
for any $ y = [ \dots, y_0, y_1, y_2, \dots ] \in \ltwo^q $ and the 
corresponding class of uncertainties 
\begin{align}
	\mathbf{\Delta} = \lbrace \Delta_{\fObj_2} : \fObj_2 \in {\classObj{0}{L_2-m_2}}  \rbrace.
	\label{eqDefClassOfUncertaintiesStructuredObj}
\end{align}
Note that, again, the so-defined $ \Delta \in \mathbf{\Delta} $ is a slope-restricted
operator in the sector $ [ 0,L_2-m_2 ] ${. Hence} we can make use of all of
the previously derived IQCs for the new uncertainty. In a nutshell, this means
that we can employ the very same techniques with the substitutions
$ \Aopt + m \Bopt \Copt \mapsto \Aopt + \Bopt \fObj_1 \Copt + m_2 \Bopt \trafoStructured^\top \trafoStructured \Copt $, $ \Bopt \mapsto \Bopt \trafoStructured^\top $,
$ \Copt \mapsto \trafoStructured \Copt $, $ L-m \mapsto L_2 - m_2 $.
We emphasize that this also applies to the convex synthesis procedures
presented in~\Cref{subsecConservativeConvexSynthesis} using corresponding
substitutions for $ Q_\Aopt $, $ Q_\Bopt $.

In the following numerical example we show that exploiting structural
knowledge about the objective function can lead to a significant
improvement in convergence rate guarantees. We choose $ \fObj_1 = \textup{blkdiag}(1,2,10,4) $,
i.e., $ {m} = 1 $, and further let $\trafoStructured$ be given by
\begin{align}
	\trafoStructured =
	\left[
	\begin{array}{rrrr}
		2  & -7 &  0 & 5 \\
		-1 &  4 & -3 & 2 \\
		0  & -2 &  1 & 0 
	\end{array}\right]
\end{align}
$ m_2 = 1 $, and let $ L_2 $ vary from $ 1 $ to $ 20 $,
resulting in $ L $ ranging from $ 89.677 $ to $ 1774.5 $.
For the probelm at hand, we design tailored algorithms
based on the convex synthesis procedure from~\Cref{subsecConservativeConvexSynthesis}
and the BMI optimization approach described in~\Cref{secBMIoptimization}.
The convergence rates of the respective algorithms are depicted in~\Cref{figSynthStructured}.
In the considered example, structure exploiting algorithms provide much better
guaranteed convergence rates compared to the fastest known standard algorithm,
the Triple Momentum Method, using only $ m $ and $ L $ as known
parameters but neglecting any structural properties of $ \fObj $.
This effect is most considerable for small values of $ L_2 $
resulting in a small class of uncertainties $ \mathbf{\Delta} $ in~\eqref{eqDefClassOfUncertaintiesStructuredObj}.
{In particular, for $ L_2 = m_2 = 1 $, there is no uncertainty
and \eqref{eqOptAlgoStandardFormStructured} is linear, hence resulting in arbitrarily small convergence rates.}
Still, we emphasize that the achievable improvements heavily depend on the 
specific problem. The numerical results also show that additionally
assuming that {$ \fObj_2 \in {\classObj{m}{L}^{\rep}} $} and employing the 
results from~\Cref{subsecIQCsForStructuredObjectiveFunctions}
can further improve the guaranteed convergence rates.
\ifthenelse{\boolean{singleColumn}}{
\setlength\figureheight{0.5\textwidth}
\setlength\figurewidth{0.5\textwidth}
}{
\setlength\figureheight{0.85\columnwidth}
\setlength\figurewidth{0.85\columnwidth}
}
\begin{figure}[t]
	\centering
	\input{../figs/exampleStructure_paper.tex}
	\caption{Convergence rate guarantees provided by structure exploiting algorithms
	designed for the example from~\Cref{subsecParametrizedObjectiveFunctions}
	using the conservative convex approach ({\protect\tikz\protect\node[rectangle,fill=mycolor2,inner sep=0pt,minimum width=4pt,minimum height=4pt] {};}),
	the BMI optimization approach ({\protect\tikz\protect\node[regular polygon, regular polygon sides=3,fill=mycolor3,inner sep=0pt,minimum width=6pt,minimum height=6pt] {};}), 
	as well as the BMI optimization approach under the additional assumption that {$ \fObj_2 \in {\classObj{m}{L}^{\rep}} $}
	({\protect\tikz\protect\node[regular polygon, regular polygon sides=5,fill=mycolor4,inner sep=0pt,minimum width=6pt] {};}).
	For comparison, we plot the convergence rates provided by the Triple Momentum Method 
	({\protect\tikz\protect\node[circle,fill=mycolor1,inner sep=0pt,minimum width=4pt] {};})
	and the lower bound for any first order algorithm when the structure is neglected
	({\protect\tikz\protect\node[diamond,fill=black,inner sep=0pt,minimum width=4pt,minimum height=6pt] {};}).
	}
	\label{figSynthStructured}
\end{figure}

\section{Conclusions and Outlook}
We presented a novel and general framework for analyzing and designing 
robust {and structure exploiting} optimization algorithms 
suitable for solving a class of unconstrained optimization
problems with strongly convex cost and possible additional structure.
Building upon the well-studied field
of robust control theory based on integral quadratic constraints
and adapting it to our needs, we provide an approach to the design
of {robust and structure exploiting} optimization algorithms 
specifically tailored to the class of optimization problems at
hand. 
Several numerical examples illustrate that tailored algorithms
designed following the presented methodology can outperform
standard algorithms in terms of robustness {against noise} and guaranteed convergence
rates in the considered scenarios.

One key advantage of the approach is that it allows for systematic
extensions by suitable adaptations of existing results from robust
control theory. In particular, it is to be expected that further
characterizations of the class of objective functions as well as
performance characterizations, both stated in terms of suitable 
integral quadratic constraints, can be embedded in the presented
framework. Further, extending the class of optimization algorithms
to primal dual dynamics or utilizing a {(relaxed)} barrier function based approach,
it might be possible to establish a similar
framework applicable to constrained optimization problems.

\section{Acknowledgements}
The second author would like to cordially thank Pete Seiler 
and Joaquin Carrasco for bringing references~\citet{freeman2018noncausal} and~\citet{zhang2019noncausal}
to our attention and for the continued interesting discussions 
on IQC theory and Zames-Falb multipliers. 
The second author gratefully acknowledges funding for this project by Deutsche
Forschungsgemeinschaft (DFG, German Research Foundation) under Germany's Excellence
Strategy - EXC 2075 - 390740016.

\bibliographystyle{abbrv}        
\bibliography{subfiles/bibfile,subfiles/ownWork}        

\begin{thebibliography}{10}

\bibitem{antipin1994minimization}
A.~S. Antipin.
\newblock Minimization of convex functions on convex sets by means of
  differential equations.
\newblock {\em Differential equations}, 30(9):1365--1375, 1994.

\bibitem{aybat2019robust}
N.~S. Aybat, A.~Fallah, M.~G\"{u}rb\"{u}zbalaban, and A.~Ozdaglar.
\newblock Robust accelerated gradient methods for smooth strongly convex
  functions.
\newblock {\em ArXiv e-prints arXiv:1805.10579}, 2019.

\bibitem{bhaya2006control}
A.~Bhaya and E.~Kaszkurewicz.
\newblock {\em Control perspectives on numerical algorithms and matrix
  problems}, volume~10 of {\em {Advances in Design and Control}}.
\newblock {SIAM}, 2006.

\bibitem{boczar2017exponential}
R.~Boczar, L.~Lessard, A.~Packard, and B.~Recht.
\newblock Exponential stability analysis via {I}ntegral {Q}uadratic
  {C}onstraints.
\newblock {\em {ArXiv e-prints arXiv:1706.01337}}, 2017.

\bibitem{boczar2015exponential}
R.~Boczar, L.~Lessard, and B.~Recht.
\newblock Exponential convergence bounds using integral quadratic constraints.
\newblock In {\em {Proc.\ 54th IEEE Conf.\ Decision and Control (CDC)}}, pages
  7516--7521, 2015.

\bibitem{amato2001repeated}
F.~J. D'{A}mato, M.~A. Rotea, A.~Megretski, and U.~T. J{\"o}nsson.
\newblock New results for analysis of systems with repeated nonlinearities.
\newblock {\em Automatica}, 37(5):739--747, 2001.

\bibitem{desoer1975feedback}
C.~A. Desoer and M.~Vidyasagar.
\newblock {\em Feedback systems: {I}nput-output properties}, volume~55.
\newblock {SIAM}, 1975.

\bibitem{drori2018efficient}
Y.~Drori and A.~B. Taylor.
\newblock Efficient first-order methods for convex minimization: {A}
  constructive approach.
\newblock {\em {ArXiv e-prints arXiv:1803.05676}}, 2018.

\bibitem{drori2014firstorder}
Y.~Drori and M.~Teboulle.
\newblock Performance of first-order methods for smooth convex minimization: A
  novel approach.
\newblock {\em {Mathematical Programming}}, 145(1):451--482, 2014.

\bibitem{duerr2012algorithms}
H.-B. D\"{u}rr and C.~Ebenbauer.
\newblock On a class of smooth optimization algorithms with applications in
  control.
\newblock In {\em Proc.\ 4th {IFAC} {Conference} on {Nonlinear} {Model}
  {Predictive} {Control}}, pages 291--298, 2012.

\bibitem{fazlyab2018design}
M.~Fazlyab, M.~Morari, and V.~M. Preciado.
\newblock Design of first-order optimization algorithms via sum-of-squares
  programming.
\newblock In {\em {Proc.\ 57th IEEE Conf.\ Decision and Control (CDC)}}, pages
  4445--4452, 2018.

\bibitem{fazlyab2018iqc}
M.~Fazlyab, A.~Ribeiro, M.~Morari, and V.~Preciado.
\newblock Analysis of optimization algorithms via {Integral Quadratic
  Constraints}: {N}onstrongly convex problems.
\newblock {\em {{SIAM} {J}ournal on {O}ptimization}}, 28(3):2654--2689, 2018.

\bibitem{feller2017relaxed}
C.~Feller and C.~Ebenbauer.
\newblock Relaxed logarithmic barrier function based model predictive control
  of linear systems.
\newblock {\em {IEEE Transactions on Automatic Control}}, 62(3):1223--1238,
  2017.

\bibitem{fetzer2017absolute}
M.~Fetzer and C.~W. Scherer.
\newblock Absolute stability analysis of discrete time feedback
  interconnections.
\newblock In {\em Proc.\ 20th {IFAC} {W}orld {C}ongress}, pages 8447--8453,
  2017.

\bibitem{freeman2018noncausal}
R.~A. Freeman.
\newblock Noncausal {Z}ames-{F}alb multipliers for tighter estimates of
  exponential convergence rates.
\newblock In {\em {Proc.\ American Control Conf.\ (ACC)}}, pages 2984--2989,
  2018.

\bibitem{kao2012discrete}
C.~Kao.
\newblock On stability of discrete-time {LTI} systems with varying time delays.
\newblock {\em {IEEE Transactions on Automatic Control}}, 57(5):1243--1248,
  2012.

\bibitem{lessard2016analysis}
L.~Lessard, B.~Recht, and A.~Packard.
\newblock Analysis and design of optimization algorithms via integral quadratic
  constraints.
\newblock {\em {SIAM Journal on Optimization}}, 26(1):57--95, 2016.

\bibitem{lessard2019synthesis}
L.~Lessard and P.~Seiler.
\newblock Direct synthesis of iterative algorithms with bounds on achievable
  worst-case convergence rate.
\newblock {\em ArXiv e-prints arXiv:1904.09046}, 2019.

\bibitem{lofberg2004yalmip}
J.~L{\"{o}}fberg.
\newblock {YALMIP : A Toolbox for Modeling and Optimization in MATLAB}.
\newblock In {\em {Proc.\ CACSD Conference}}, 2004.

\bibitem{lure1944theory}
A.~I. Lur'e and V.~N. Postnikov.
\newblock On the theory of stability of control systems.
\newblock {\em Applied mathematics and mechanics}, 8(3):246--248, 1944.

\bibitem{mancera2005multipliers}
R.~Mancera and M.~G. Safonov.
\newblock All stability multipliers for repeated {MIMO} nonlinearities.
\newblock {\em Systems \& {C}ontrol {L}etters}, 54(4):389--397, 2005.

\bibitem{megretski1997iqcs}
A.~Megretski and A.~Rantzer.
\newblock System analysis via integral quadratic constraints.
\newblock {\em {IEEE Transactions on Automatic Control}}, 42(6):819--830, 1997.

\bibitem{mic2014heavy}
S.~Michalowsky and C.~Ebenbauer.
\newblock The multidimensional $n$-th order heavy ball method and its
  application to extremum seeking.
\newblock In {\em {Proc.\ 53rd IEEE Conf.\ Decision and Control (CDC)}}, pages
  2660--2666, 2014.

\bibitem{mic2016extremum}
S.~Michalowsky and C.~Ebenbauer.
\newblock Extremum control of linear systems based on output feedback.
\newblock In {\em {Proc.\ 55th IEEE Conf.\ Decision and Control (CDC)}}, pages
  2963--2968, 2016.

\bibitem{mic2019algoDesignArxiv}
S.~Michalowsky, C.~Scherer, and C.~Ebenbauer.
\newblock Robust and structure exploiting algorithms: {A}n integral quadratic
  constraint approach.
\newblock {\em ArXiv e-prints arXiv:1905.00279}, 2019.

\bibitem{mohammadi2018variance}
H.~Mohammadi, M.~Razaviyayn, and M.~R. Jovanovi{\'c}.
\newblock Variance amplification of accelerated first-order algorithms for
  strongly convex quadratic optimization problems.
\newblock In {\em {Proc.\ 57th IEEE Conf.\ Decision and Control (CDC)}}, pages
  5753--5758, 2018.

\bibitem{murphy2012machine}
K.~P. Murphy and F.~Bach.
\newblock {\em Machine {L}earning: {A} Probabilistic Perspective}.
\newblock {Adaptive Computation and Machine Learning}. {MIT Press}, 2012.

\bibitem{nesterov2004introductory}
Y.~Nesterov.
\newblock {\em Introductory lectures on convex optimization}, volume~87.
\newblock {Springer Science \& Business Media}, 2004.

\bibitem{paganini2000h2}
F.~Paganini and E.~Feron.
\newblock Linear matrix inequality methods for robust {$H_2$} analysis: {A}
  survey with comparisons.
\newblock In {\em Advances in Linear Matrix Inequality Methods in Control},
  chapter~7, pages 129--151. 2000.

\bibitem{polyak1987introduction}
B.~T. Polyak.
\newblock {\em Introduction to optimization}.
\newblock {Optimization Software Inc.}, 1987.

\bibitem{rockafellar1970convex}
R.~T. Rockafellar.
\newblock {\em Convex {A}nalysis}.
\newblock Princeton {U}niversity {P}ress, 1970.

\bibitem{romer2017sampling}
A.~Romer, J.~M. Montenbruck, and F.~Allg{\"o}wer.
\newblock Sampling strategies for data-driven inference of passivity
  properties.
\newblock In {\em {Proc.\ 56th IEEE Conf.\ Decision and Control (CDC) Decision
  and Control (CDC)}}, pages 6389--6394, 2017.

\bibitem{safavi2018explicit}
S.~Safavi, B.~Joshi, G.~Fran\c{c}a, and J.~Bento.
\newblock An explicit convergence rate for {N}esterov's method from {SDP}.
\newblock In {\em {2018 IEEE International Symposium on Information Theory
  (ISIT)}}, pages 1560--1564, 2018.

\bibitem{taylor2017pesto}
A.~B. Taylor, J.~M. Hendrickx, and F.~Glineur.
\newblock Performance estimation toolbox ({PESTO}): {A}utomated worst-case
  analysis of first-order optimization methods.
\newblock In {\em {Proc.\ 56th IEEE Conf.\ Decision and Control (CDC)}}, pages
  1278--1283, 2017.

\bibitem{vanScoy2018fastest}
B.~Van~Scoy, R.~A. Freeman, and K.~M. Lynch.
\newblock The fastest known globally convergent first-order method for
  minimizing strongly convex functions.
\newblock {\em {{IEEE} {C}ontrol {S}ystems {L}etters}}, 2(1):49--54, 2018.

\bibitem{veenman2016robust}
J.~Veenman, C.~W. Scherer, and H.~K{\"o}ro{\u g}lu.
\newblock Robust stability and performance analysis based on integral quadratic
  constraints.
\newblock {\em European J.\ Control}, 31:1--32, 2016.

\bibitem{willems1971analysis}
J.~Willems.
\newblock {\em The analysis of feedback systems}.
\newblock {MIT} {P}ress, 1971.

\bibitem{wilson2016lyapunov}
A.~C. Wilson, B.~Recht, and M.~I. Jordan.
\newblock A {L}yapunov analysis of momentum methods in optimization.
\newblock {\em {ArXiv e-prints arXiv:1611.02635}}, 2016.

\bibitem{zames1969stability}
G.~Zames and P.~Falb.
\newblock Stability conditions for systems with monotone and slope-restricted
  nonlinearities.
\newblock {\em {SIAM Journal on Control}}, 6(1):89--108, 1968.

\bibitem{zhang2019noncausal}
J.~Zhang, P.~Seiler, and J.~Carrasco.
\newblock Noncausal {FIR} {Z}ames-{F}alb multiplier search for exponential
  convergence rate.
\newblock {\em ArXiv e-prints arXiv:1902.09473}, 2019.

\end{thebibliography}

\appendix
\section{Appendix}
\ifthenelse{\boolean{longVersion}}{
\subsection{Some standard results}\label{subsecAppendixStandardResults}
For the sake of completeness, we state some well-known
results in the following that can be found in standard textbooks.
\begin{lemma}[Discrete-time KYP-Lemma]\label{lemmaKYP}
	Let $G$ be a real rational and proper transfer matrix and let
	$ (A,B,C,D) $ denote a state-space representation of $G$. 
	Suppose that $ A $ has no eigenvalues on the unit circle 
	and let $ M $ be a real symmetric matrix. Then the following frequency domain inequality 
	\begin{align}
		G(z)^\complConj M G(z) {\stackrel{\scriptscriptstyle\unitCircle}{\prec} 0} 
	\end{align}
	holds if and only if there exists a $ P = P^\top $ such that
	\begin{align}
		\left[
		\begin{array}{cc}
			A & B \\
			I & 0 \\
			C & D 
		\end{array}
		\right]^\top 
		\left[
		\begin{array}{ccc}
			P & 0  & 0 \\
			0 & -P & 0 \\
			0 & 0  & M
		\end{array}
		\right]
		\left[
		\begin{array}{cc}
			A & B \\
			I & 0 \\
			C & D 
		\end{array}
		\right]
		\prec 0 .
	\end{align}
\end{lemma}
}{}

\ifthenelse{\boolean{proofsAppendix}}{
	\subsection{Proofs}
	We collect all proofs of the presented results in the following.
	\subsubsection{Proof of~\Cref{lemmaIQC}}\label{secProofLemmaIQC}
	
	\subsubsection{Proof of~\Cref{lemmaConditionsAlgorithm}}\label{secProofConditionsAlgorithm}
	
	\subsubsection{Proof of~\Cref{lemmaIQCexponentialStability}}\label{secProofIQCexponentialStability}
	
	\subsubsection{Proof of~\Cref{lemmaZFGeneralTimeDomain}}\label{secProofZFGeneralTimeDomain}
	
	\subsubsection{Proof of~\Cref{lemmaMultipliersTransformedLoop}}\label{secProofLemmaMultipliersTransformedLoop}
	 {Let $ \opTimeDomain_T = \trunc{T} \opTimeDomain \trunc{T} $, $ T \in \nat $.} {
Since $ \rho_+(\ltwo^p) = \ltworho^p $, we infer from~\eqref{eqIQClemmaMultipliersTransformedLoop} {that}
\begin{align}
	\left\langle \begin{bmatrix} \rho_+(y) \\ \Delta(\rho_+(y)) \end{bmatrix}, {\opTimeDomain_T} \begin{bmatrix} \rho_+(y) \\ \Delta(\rho_+(y)) \end{bmatrix} \right\rangle \geq 0 
\end{align}
for all $ y \in \ltwo^p $, $ \Delta \in \mathbf{\Delta} $, 
{$T \in \nat$}.
With $ \rho_+ \circ \rho_- = \text{id}_{\ltworho} $ this is 
equivalently formulated as
\begin{align}
	&\left\langle \begin{bmatrix} \rho_+(y) \\ \rho_+\Big( \rho_- \big(\Delta(\rho_+(y)) \big)\Big) \end{bmatrix}, {\opTimeDomain_T} \begin{bmatrix} \rho_+(y) \\ \rho_+\Big( \rho_- \big(\Delta(\rho_+(y)) \big)\Big) \end{bmatrix} \right\rangle \nonumber \\
	=~& \left\langle  \begin{bmatrix} y \\ \Delta_\rho(y) \end{bmatrix}, \rho_+\Big( {\opTimeDomain_T} \big( \rho_+ (  \begin{bmatrix} y \\ \Delta_\rho(y) \end{bmatrix} ) \big) \Big) \right\rangle \nonumber \\
	=~& \left\langle  \begin{bmatrix} y \\ \Delta_\rho(y) \end{bmatrix}, {\trunc{T} \Big( \tilde{\opTimeDomain} \big( \trunc{T} (}  \begin{bmatrix} y \\ \Delta_\rho(y) \end{bmatrix} ) \big) \Big) \right\rangle
	\geq 0 
	\label{eqInequalityProofMultipliersTransformedLoop}
\end{align}
for all $ y \in \ltwo^p $, $ \Delta \in \mathbf{\Delta} $, {$T \in \nat$}.
{Since $ \Delta (\ltworho^p) \subset \ltworho^p $, we have $ \Delta_\rho: \ltwo^p \to \ltwo^p $;
hence, by the assumption that $ \tilde{\opTimeDomain} $ is bounded on $ \ltwo^p $,
{we can take the limit $ T \to \infty $ in~\eqref{eqInequalityProofMultipliersTransformedLoop}
to obtain that}~\eqref{eqIQCdeltaRho} holds.}}

	\subsubsection{Proof of~\Cref{lemmaZFRhoTimeDomain}}\label{secProofZFRhoTimeDomain}
	
	\subsubsection{Proof of~\Cref{lemmaAnalysisConvRate}}\label{secProofLemmaAnalysisConvRate}
	
	\ifthenelse{\boolean{longVersion}}{\subsubsection{Proof of dimensionality reduction in~\Cref{lemmaAnalysisConvRate}}\label{secAppendixProofDimensionalityReduction}
	We follow the {arguments} of~\citet{lessard2016analysis} to prove this statement. It is clear that
if~\eqref{eqLMIConvAnalysisReduced} has a solution $ \bar{P}, m_+, ,m_-, m_0 $, then 
$ P = \overline{P} \otimes I_p, M_+ = m_+ \otimes I_p, m_- \otimes I_p, m_0 \otimes I_p $. Now suppose that~\eqref{eqLMIConvAnalysis}
has a solution $ P, M_+, M_-, M_0 $. Multiply~\eqref{eqLMIConvAnalysis}
from right and left by $ \textup{blkdiag}( I_{n_\complete} \otimes e_1, I_{n_\complete} \otimes e_1 ) $ and its transpose, where
$ e_1 \in \real^{p \times 1} $ is the first unit vector. Observing 
that for any $ U =\bar{U} \otimes I_p $, $ \bar{U} \in \real^{n_\complete \times r} $, we have
$ U (I_n \otimes e_1) = ( \bar{U} \otimes I_p ) ( I_{n_\complete} \otimes e_1 ) = \bar{U} \otimes e_1 = ( I_{r} \otimes e_1 ) ( \bar{U} \otimes 1 ) $,
we obtain 
\ifthenelse{\boolean{singleColumn}}{
\begin{align}
	\begin{bmatrix} \overline{\AcompleteNoPerf} & \overline{\BcompleteNoPerf} \\ I & 0 \end{bmatrix}^\top 
	\begin{bmatrix} \overline{P} & 0 \\ 0 & -\overline{P} \end{bmatrix}
	\begin{bmatrix} \overline{\AcompleteNoPerf} & \overline{\BcompleteNoPerf} \\ I & 0 \end{bmatrix} 
	+
	\begin{bmatrix} \overline{\CcompleteNoPerf} & \overline{\DcompleteNoPerf} \end{bmatrix}^\top ( I_{3p} \otimes e_1 )^\top M_\Delta ( I_{3p} \otimes e_1 ) \begin{bmatrix} \overline{\CcompleteNoPerf} & \overline{\DcompleteNoPerf} \end{bmatrix}
	\prec 0 
\end{align}}{
\begin{align}
	&\begin{bmatrix} \overline{\AcompleteNoPerf} & \overline{\BcompleteNoPerf} \\ I & 0 \end{bmatrix}^\top 
	\begin{bmatrix} \overline{P} & 0 \\ 0 & -\overline{P} \end{bmatrix}
	\begin{bmatrix} \overline{\AcompleteNoPerf} & \overline{\BcompleteNoPerf} \\ I & 0 \end{bmatrix} \\
	&+
	\begin{bmatrix} \overline{\CcompleteNoPerf} & \overline{\DcompleteNoPerf} \end{bmatrix}^\top ( I_{3p} \otimes e_1 )^\top M_\Delta ( I_{3p} \otimes e_1 ) \begin{bmatrix} \overline{\CcompleteNoPerf} & \overline{\DcompleteNoPerf} \end{bmatrix}
	\prec 0 \nonumber
\end{align}}
with $ \bar{P} = (I_n \otimes e_1)^\top P (I_n \otimes e_1) $. 
Let $ M_\Delta(M_+,M_-,M_0) $ be defined as in~\eqref{eqDefMDelta}
explicitly including the dependency on $M_+,M_-,M_0$,
and note that
\begin{align}
	&( I_{3p} \otimes e_1 )^\top M_\Delta(M_+,M_-,M_0) ( I_{3p} \otimes e_1 ) \nonumber \\
	=~& ( I_{3p} \otimes e_1 )^\top ( M_\Delta(m_+,m_-,m_0) \otimes I_p ) ( I_{3p} \otimes e_1 ) \nonumber \\
	=~& ( I_{3p} \otimes e_1 )^\top ( M_\Delta(m_+,m_-,m_0) \otimes e_1 ) \nonumber \\
	=~& M_\Delta(m_+,m_-,m_0).
\end{align}
Hence, $ \bar{P} = (I_n \otimes e_1)^\top P (I_n \otimes e_1), m_+, m_-, m_0 $
is a solution to~\eqref{eqLMIConvAnalysisReduced} which concludes the proof.
}{}
	\subsubsection{Proof of~\Cref{lemmaAnalysisH2}}\label{secProofLemmaAnalysisH2}
	{We first observe that robust stability follows directly from~\Cref{lemmaAnalysisConvRate}
noting that $ \Acomplete = \AcompleteNoPerf(1) $, $ \Bcomplete{}_{,1} = \BcompleteNoPerf $,
$ \Ccomplete{}_{,1} = \CcompleteNoPerf(1) $, $ \Dcomplete{}_{,11} = \DcompleteNoPerf $ 
by~\eqref{eqSSrealizationCompleteSystem} and $ \Ccomplete{}_{,1}^\top \Ccomplete{}_{,1} \succeq 0 $;
hence~\eqref{eqLMIH2Analysis} implies that~\eqref{eqAllConstraintsConvAnalysis} holds with $ \rho=1 $.}
{Consider the complete transfer matrix $ \tfcomplete $ as defined
in~\eqref{eqDefCompleteTransferFunctionFDI} and let a state-space realization
be given as in~\eqref{eqSSCompleteSystemWithPerformance} specifically stated
in~\eqref{eqSSrealizationCompleteSystem} in~\Cref{secAppendixStateSpaceRealizations}.}
Let $ x_\complete $ denote the solution of
\begin{subequations}
\begin{align}	
	x_{\complete,k+1} &= \Acomplete x_{\complete,k} + \Bcomplete{}_{,1} w_k + \Bcomplete{}_{,2} w_{\perf,k}  \\
	w_k &= \Delta( { \Copt N^\top {x}_{\complete,k} } ) 
\end{align}
\end{subequations}
with initial condition $ x_{\complete,0} = 0 $. 
{Note that with $ \bar{x}_{\complete,k} = N^\top {x}_{\complete,k} $
we have 
\begin{align}
	\bar{x}_{\complete,k+1} = ( \Aopt + m \Bopt \Copt ) \bar{x}_{\complete,k} + \Bopt \Delta( \Copt \bar{x}_{\complete,k} ) + \Bperf w_{\perf,k},
\end{align}
i.e., $ \bar{x}_{\complete,k} $ follows the same dynamics as~\eqref{eqRobustStandardForm} and
the performance output is given by $ y_{\perf,k} = \Cperf \bar{x}_{\complete,k} = \Ccomplete{}_{,2} x_{\complete,k} $.}
Consequently, multiplying~\eqref{eqLMIH2Analysis} from left by $ \begin{bmatrix} x_{\complete,k}^\top & w_k^\top \end{bmatrix} $
and from right by its transpose, we infer
\ifthenelse{\boolean{singleColumn}}{
\begin{align}
	( \star )^\top P_\perf ( \Acomplete x_{\complete,k} + \Bcomplete{}_{,1} w_k ) - x_{\complete,k}^\top P_\perf x_{\complete,k} 
	{\leq} - y_{\perf,k}^\top y_{\perf,k} - (\star)^\top M_\Delta ( \Ccomplete{}_{,1} x_{\complete,k} + \Dcomplete{}_{,11} w_k ).
	\label{eqIneqProofH2_1}
\end{align}
}{
\begin{align}
	( \star )^\top P_\perf ( \Acomplete x_{\complete,k} + \Bcomplete{}_{,1} w_k ) - x_{\complete,k}^\top P_\perf x_{\complete,k} \nonumber \\
	{\leq} - y_{\perf,k}^\top y_{\perf,k} - (\star)^\top M_\Delta ( \Ccomplete{}_{,1} x_{\complete,k} + \Dcomplete{}_{,11} w_k ).
	\label{eqIneqProofH2_1}
\end{align}}
Let 
{
\begin{align}
	P_\perf' = N^\top P_\perf N, \qquad
	P_\perf'' = P_\perf - NN^\top P_\perf N N^\top,
	\label{eqDefPperfPrime}
\end{align}}
i.e., $ P_\perf = P_\perf'' + N P_\perf' N^\top $ {and $ N^\top P_\perf'' N = 0 $}.
{By~\eqref{eqSSrealizationCompleteSystem}, we have $ \Bcomplete{}_{,2} = N \Bperf $;
hence we infer that $ \Bcomplete{}_{,2}{}^\top P_\perf \Bcomplete{}_{,2} = \Bcomplete{}_{,2}^\top N P_\perf' N^\top \Bcomplete{}_{,2} $.}
We then calculate 
\ifthenelse{\boolean{singleColumn}}{
\begin{align}
	&x_{\complete,k+1}^\top P_\perf x_{\complete,k+1} - x_{\complete,k}^\top P_\perf x_{\complete,k} \nonumber \\
	=~& ( \star )^\top P_\perf ( \Acomplete x_{\complete,k} + \Bcomplete{}_{,1} w_k ) - x_{\complete,k}^\top P_\perf x_{\complete,k} 
	+ 2 { x_{\complete,k+1} ^\top} {P_\perf} \Bcomplete{}_{,2} w_{\perf,k} \nonumber \\
	&+ w_{\perf,k}^\top \Bcomplete{}_{,2}^\top N P_\perf' N^\top \Bcomplete{}_{,2} w_{\perf,k}.
	\label{eqIneqProofH2_2}
\end{align}}{
\begin{align}
	&x_{\complete,k+1}^\top P_\perf x_{\complete,k+1} - x_{\complete,k}^\top P_\perf x_{\complete,k} \nonumber \\
	=~& ( \star )^\top P_\perf ( \Acomplete x_{\complete,k} + \Bcomplete{}_{,1} w_k ) - x_{\complete,k}^\top P_\perf x_{\complete,k} 	\label{eqIneqProofH2_2} \\
	& + 2 {x_{\complete,k+1}^\top} {P_\perf} \Bcomplete{}_{,2} w_{\perf,k} 
	+ w_{\perf,k}^\top \Bcomplete{}_{,2}^\top N P_\perf' N^\top \Bcomplete{}_{,2} w_{\perf,k}. \nonumber
\end{align}}
Using~\eqref{eqIneqProofH2_1}, \eqref{eqIneqProofH2_2} and summing
from $ k = 0 $ to $ k_{\textup{max}} \in \nat $, we obtain
\ifthenelse{\boolean{singleColumn}}{
\begin{align}
	\sum\limits_{k=0}^{k_{\textup{max}}} x_{\complete,k+1}^\top P_\perf x_{\complete,k+1} - x_{\complete,k}^\top P_\perf x_{\complete,k} 
	{\leq}~& \sum\limits_{k=0}^{k_{\textup{max}}} \Big( - y_{\perf,k}^\top y_{\perf,k} - (\star)^\top M_\Delta ( \Ccomplete{}_{,1} x_{\complete,k} + \Dcomplete{}_{,11} w_k ) \nonumber \\
	& + 2 {x_{\complete,k+1}^\top} {P_\perf} \Bcomplete{}_{,2} w_{\perf,k} \nonumber \\
	&+ w_{\perf,k}^\top \Bcomplete{}_{,2}^\top N P_\perf' N^\top \Bcomplete{}_{,2} w_{\perf,k} \Big). \label{eqIneqProofH2_3}
\end{align}}{
\begin{align}
	&\sum\limits_{k=0}^{k_{\textup{max}}} x_{\complete,k+1}^\top P_\perf x_{\complete,k+1} - x_{\complete,k}^\top P_\perf x_{\complete,k} \nonumber \\
	{\leq}~& \sum\limits_{k=0}^{k_{\textup{max}}} \Big( - y_{\perf,k}^\top y_{\perf,k} - (\star)^\top M_\Delta ( \Ccomplete{}_{,1} x_{\complete,k} + \Dcomplete{}_{,11} w_k ) \nonumber \\
	& + 2 {x_{\complete,k+1}^\top} {P_\perf} \Bcomplete{}_{,2} w_{\perf,k} + w_{\perf,k}^\top \Bcomplete{}_{,2}^\top N P_\perf' N^\top \Bcomplete{}_{,2} w_{\perf,k} \Big). \label{eqIneqProofH2_3}
\end{align}}
{We note that $ w_{\perf,k} $ is a discrete-time white noise process}
with independent components; hence, {$ w_{\perf,k} $ and $ x_{\complete,k+1} $
are independent for any $ k \in \nat $}, $ \expec(w_{\perf,k}) = 0 $ 
and $ \expec( w_{\perf,k}^\top X w_{\perf,k} ) = \textup{tr}(X) $ for any $X \in \real^{n_{w_\perf} \times n_{w_\perf}} $.
{Taking expectations, we hence infer}
\ifthenelse{\boolean{singleColumn}}{
\begin{align}
	\expec\big( \sum\limits_{k=0}^{k_{\textup{max}}} x_{\complete,k+1}^\top P_\perf x_{\complete,k+1} - x_{\complete,k}^\top P_\perf x_{\complete,k} \big)
	{\leq}~& \expec\big( \sum\limits_{k=0}^{k_{\textup{max}}} -  y_{\perf,k}^\top y_{\perf,k}  - (\star)^\top M_\Delta ( \Ccomplete{}_{,1} x_{\complete,k} + \Dcomplete{}_{,11} w_k ) \big) \nonumber \\
	&+ k_{\textup{max}} \textup{tr}( \Bcomplete{}_{,2}^\top N P_\perf' N^\top \Bcomplete{}_{,2} ). \label{eqIneqProofH2_4}
\end{align}}{
\begin{align}
	&\expec\big( \sum\limits_{k=0}^{k_{\textup{max}}} x_{\complete,k+1}^\top P_\perf x_{\complete,k+1} - x_{\complete,k}^\top P_\perf x_{\complete,k} \big) \nonumber \\
	{\leq}~& \expec\big( \sum\limits_{k=0}^{k_{\textup{max}}} -  y_{\perf,k}^\top y_{\perf,k}  - (\star)^\top M_\Delta ( \Ccomplete{}_{,1} x_{\complete,k} + \Dcomplete{}_{,11} w_k ) \big) \nonumber \\
	&+ k_{\textup{max}} \textup{tr}( \Bcomplete{}_{,2}^\top N P_\perf' N^\top \Bcomplete{}_{,2} ). \label{eqIneqProofH2_4}
\end{align}}
{With $ x_{\complete,0} = 0 $ we further note that 
for any $ k_{\textup{max}} \in \nat $ we have
\begin{align}
	\sum\limits_{k=0}^{k_{\textup{max}}} x_{\complete,k+1}^\top P_\perf x_{\complete,k+1} - x_{\complete,k}^\top P_\perf x_{\complete,k} = x_{\complete,k_{\textup{max}}+1}^\top P_\perf x_{\complete,k_{\textup{max}}+1}.
	\label{eqEqProofH2_5}
\end{align}}
Combining~\eqref{eqIneqProofH2_4}
with~\eqref{eqEqProofH2_5} and using~\eqref{eqTraceH2Analysis} we then obtain
\begin{align}
	&\tfrac{1}{k_{\textup{max}}} \sum\limits_{k=0}^{k_{\textup{max}}} \expec( y_{\perf,k}^\top y_{\perf,k} ) {+ \tfrac{1}{k_{\textup{max}}} \expec(x_{\complete,k_{\textup{max}}+1}^\top P_\perf x_{\complete,k_{\textup{max}}+1})} \label{eqEqProofH2_6} \\
	& \quad {\leq} - \tfrac{1}{k_{\textup{max}}}  \expec\big( \sum\limits_{k=0}^{k_{\textup{max}}} (\star)^\top M_\Delta ( \Ccomplete{}_{,1} x_{\complete,k} + \Dcomplete{}_{,11} w_k ) \big) + \gamma^2 \nonumber 
\end{align}
Now note that, {for a fixed realization $ w_{\perf,k} $}, $ \sum_{k=0}^{\infty} (\star)^\top M_\Delta ( \Ccomplete{}_{,1} x_{\complete,k} + \Dcomplete{}_{,11} w_k ) $
is a time-domain representation of $ \IQC\big( \psi_\Delta M_\Delta \psi_\Delta, y, \Delta_\fObj(y) \big) $.
{Since~\eqref{eqIneqZamesFalbTimeDomain} is a hard IQC, we infer that
\begin{align}
	\sum\limits_{k=0}^{k_{\textup{max}}} (\star)^\top M_\Delta ( \Ccomplete{}_{,1} x_{\complete,k} + \Dcomplete{}_{,11} w_k ) \geq 0
\end{align}
for any realization $ w_{\perf,k} $ and any $ k_{\textup{max}} $, and the
latter inequality persists to hold after taking expectations.}
{Additionally, with $P_\perf$ being positive definite, 
$ \expec(x_{\complete,k}^\top P_\perf x_{\complete,k}) $ is positive as well.}
Thus, if we let $ k_{\textup{max}} $ tend to infinity, we 
finally obtain from~\eqref{eqEqProofH2_6}
\begin{align}
	\limsup\limits_{k_{\textup{max}} \to \infty}  \tfrac{1}{k_{\textup{max}}} \sum\limits_{k=0}^{k_{\textup{max}}} \expec( y_{\perf,k}^\top y_{\perf,k} ) \leq \gamma^2,
\end{align}
hence concluding the proof.

	\subsubsection{Proof of~\Cref{lemmaConservativeConvexSynthesisWithPerf}}\label{secProofLemmaConvexSynthesisWithPerf}
	
}{
}

\subsection{Additional material}
\subsubsection{Comparison to~\citet{freeman2018noncausal}}\label{secAppendixDiscussionComparison}
{In the following we compare the set of multipliers derived in~\citet{freeman2018noncausal}
to the multipliers defined in~\eqref{eqDefZamesFalbMultipliers}. For the sake of simplicity we
limit ourselves to the case of monotone uncertainties $ \Delta $; we emphasize that the
same applies to {general slope-restricted uncertainties} with minor adaptations. 
Let $ M = \toep( [m_j]_{j \in {\lbrace-\dimCausal,\dots,\dimAnticausal\rbrace} } ), 
\bar{M} = \toep( [\bar{m}_j]_{j \in {\lbrace-\dimCausal,\dots,\dimAnticausal\rbrace}} ) $, 
{$ \dimAnticausal, \dimCausal \in \nat $,}
be two Toeplitz operators $ M, \bar{M} : {\lf}^1 \to \ell_e^1 $.
{In contrast to~\citet{freeman2018noncausal}, we consider only finite
Toeplitz operators here; the following arguments also apply to infinite operators.}
In~\citet{freeman2018noncausal}, the author considers operators $ \bar{M} $ with
the property $ \bar{M} y = y - h * y $ for all {$ y \in \ell_f^1 $}, where
$ * $ is the convolution operator and $ h(z) = \sum_{j=-{\dimCausal}}^{{\dimAnticausal}} h_j z^{-j} $,
i.e., $ \bar{M} $ defines the transfer function $ E_{\bar{M}}(z) = 1 - h(z) $. 
It is then shown that if 
\begin{subequations}
\begin{align}
	h_j \geq 0 \text{ for } j \in {\lbrace-\dimCausal,\dots,\dimAnticausal\rbrace} \\
	\sum\limits_{j=-{\dimCausal}}^{{\dimAnticausal}} h_j \max \lbrace 1, \rho^{-2j} \rbrace \leq 1,
\end{align}\label{eqConditionsFreemanOriginal}
\end{subequations}
then 
\begin{align}
	\langle \bar{M} y, \Delta(y) \rangle_w = \langle \rho_- \bar{M} y, \rho_-  \Delta(y) \rangle \geq 0
	\label{eqPositivityFreeman}
\end{align}
for all $ y \in \ltworho $, where $ \langle x,y \rangle_w = \langle \rho_- x, \rho_- y \rangle $
denotes the weighted inner product on $ \ltworho $.
Condition~\eqref{eqConditionsFreemanOriginal} is equivalently formulated
in the parameters $ \bar{m}_j $, $ j \in {\lbrace-\dimCausal,\dots,\dimAnticausal\rbrace} $, as
\begin{subequations}
\begin{align}
	\bar{m}_j \leq 0 \text{ for } j \in {\lbrace-\dimCausal,\dots,\dimAnticausal\rbrace} \setminus \lbrace 0 \rbrace \text{ and } \bar{m}_0 \geq -1 \label{eqConditionsFreeman1} \\
	\sum\limits_{j=-{\dimCausal}}^{{\dimAnticausal}} \bar{m}_{-j} \max \lbrace 1, \rho^{-2j} \rbrace \geq 0. \label{eqConditionsFreeman2}
\end{align}\label{eqConditionsFreeman}
\end{subequations}
Similarly, in~\Cref{lemmaZFRhoTimeDomain} we have shown that if 
\begin{subequations}
\begin{align}
	m_j \leq 0 \text{ for } j \in {\lbrace-\dimCausal,\dots,\dimAnticausal\rbrace} \setminus \lbrace 0 \rbrace \label{eqConditionsOur1} \\
	\sum\limits_{j=-{\dimCausal}}^{{\dimAnticausal}}  \rho^{-j} m_j \geq 0 \label{eqConditionsOur2} \\
	\sum\limits_{j=-{\dimCausal}}^{{\dimAnticausal}}  \rho^{j} m_j \geq 0. \label{eqConditionsOur3}
\end{align}\label{eqConditionsOur}
\end{subequations}
then 
\ifthenelse{\boolean{singleColumn}}{
\begin{align}
	\langle M y, \Delta_{\rho}(y) \rangle 
	= \langle M y, \rho_- \Delta (\rho_+ y) \rangle 
	= \langle M \rho_- \rho_+ y, \rho_- \Delta (\rho_+ y) \rangle \geq 0 
\end{align}}{
\begin{align}
	&\langle M y, \Delta_{\rho}(y) \rangle \nonumber \\
	=~& \langle M y, \rho_- \Delta (\rho_+ y) \rangle \nonumber \\
	=~& \langle M \rho_- \rho_+ y, \rho_- \Delta (\rho_+ y) \rangle \geq 0 
\end{align}}
for all $ y \in \ltwo $, i.e., equivalently,  
\begin{align}
	\langle M \rho_- y, \rho_- \Delta(y) \rangle \geq 0
	\label{eqPositivityOur}
\end{align}
for all $ y \in \ltworho $. 
We hence note that the operators $\bar{M}$ and $M$ in~\eqref{eqPositivityFreeman}
and~\eqref{eqPositivityOur} are related by $ \rho_- \bar{M} = M \rho_- $ and, {therefore,}
\begin{align}
	\bar{m}_j = \rho^{-j} m_j \qquad \text{for } j \in {\lbrace-\dimCausal,\dots,\dimAnticausal\rbrace}.
\end{align}
We are now ready to compare the conditions~\eqref{eqConditionsFreeman} from~\cite{freeman2018noncausal}
with the proposed conditions~\eqref{eqConditionsOur}. To this end, we express~\eqref{eqConditionsFreeman}
in terms of $m_j$ instead of $\bar{m}_j$, i.e., we have
\begin{subequations}
\begin{align}
	\rho^{-j} m_j \leq 0 \text{ for } j \in {\lbrace-\dimCausal,\dots,\dimAnticausal\rbrace} \setminus \lbrace 0 \rbrace \text{ and } {m}_0 \geq -1 \\
	\sum\limits_{j=-{\dimCausal}}^{{\dimAnticausal}}  \rho^j m_{-j} \max \lbrace 1,\rho^{-2j} \rbrace \geq 0.
\end{align}
\end{subequations}
Note that with $ \rho \in (0,1) $ this holds if and only if
\begin{subequations}
\begin{align}
	m_j \leq 0 \text{ for } j \in {\lbrace-\dimCausal,\dots,\dimAnticausal\rbrace} \setminus \lbrace 0 \rbrace \label{eqConditionsFreemanReformulated1} \\
	\sum\limits_{j=-{\dimCausal}}^{{\dimAnticausal}}  \rho^j m_{j} + \sum\limits_{j=-{\dimCausal}}^{{\dimAnticausal}}  \rho^{-j} m_{j} = \sum\limits_{j=-{\dimCausal}}^{{\dimAnticausal}}  \rho^{-\vert j \vert} m_j \geq 0. \label{eqConditionsFreemanReformulated2}
\end{align}
\end{subequations}
We first note that~\eqref{eqConditionsFreemanReformulated1} is the same
as~\eqref{eqConditionsOur1} and next show that~\eqref{eqConditionsFreemanReformulated2}
implies~\eqref{eqConditionsOur2}, \eqref{eqConditionsOur3}.
To this end, observe that for any $ j \in \integers \setminus \lbrace 0 \rbrace $ {we have}
\begin{align}
	\rho^{\vert j \vert} m_j \geq m_j \geq \rho^{-\vert j \vert} m_j
\end{align}
since $ \rho \in (0,1) $ and $ m_j \leq 0 $. Consequently,
\begin{subequations}
\begin{align}
	\sum\limits_{j=1}^{{\dimAnticausal}} \rho^{-j} m_{j} &\leq \sum\limits_{j=1}^{{\dimAnticausal}} \rho^j m_j \label{eqInequality1} \\
	\sum\limits_{j=-{\dimCausal}}^{0} \rho^j m_{j} = \sum\limits_{j=-{\dimCausal}}^{0} \rho^{- \vert j \vert} m_{j} &\leq \sum\limits_{j=-{\dimCausal}}^{0} \rho^{\vert j \vert} m_{j} = \sum\limits_{j=-{\dimCausal}}^{0} \rho^{- j } m_{j}. \label{eqInequality2}
\end{align}
\end{subequations}
Thus, \eqref{eqConditionsFreemanReformulated2}
implies~\eqref{eqConditionsOur2} by \eqref{eqInequality2}
and \eqref{eqConditionsOur3} by \eqref{eqInequality1}.
The converse is in general not true as we will show next.
{For the sake of a clearer presentation, we 
suppose that $ \dimCausal = \dimAnticausal $.}
Summing up the two inequalities~\eqref{eqConditionsOur2},
\eqref{eqConditionsOur3}, we obtain
\begin{align}
	m_0 + \tfrac{1}{2} \sum\limits_{i=1}^{{\dimCausal}} (\rho^{i} + \rho^{-i}) (m_i + m_{-i}) \geq 0.
	\label{eqInequality3}
\end{align}
Note further that~\eqref{eqConditionsFreemanReformulated2}
is equivalently formulated as
\begin{align}
	m_0 + \sum\limits_{i=1}^{{\dimCausal}} \rho^{-i} (m_i+m_{-i}) \geq 0.
	\label{eqInequality4}
\end{align}
Suppose now that
\begin{align}
	m_0 = - \tfrac{1}{2} \sum\limits_{i=1}^{{\dimCausal}} (\rho^{i} + \rho^{-i}) (m_i + m_{-i})
\end{align}
which clearly fulfills~\eqref{eqInequality3}. However, it is
\begin{align}
	 - \tfrac{1}{2} \sum\limits_{i=1}^{{\dimCausal}} (\rho^{i} + \rho^{-i}) (m_i + m_{-i})  +  \sum\limits_{i=1}^{{\dimCausal}} \rho^{-i} (m_i+m_{-i}) \ifthenelse{\boolean{singleColumn}}{}{\nonumber\\}
	 = \sum\limits_{i=1}^{{\dimCausal}} (\rho^{-i} - \rho^{i}) (m_i+m_{-i}) \leq 0
\end{align}
and hence $m_0$ does {fulfill~\eqref{eqInequality4} only if
$ m_i = m_{-i} = 0 $ for $ i = 1,2,\dots,\dimCausal $.}
We hence conclude that the conditions on the anticausal parts proposed in the present paper
are less restrictive than {those} derived in~\cite{freeman2018noncausal}. 
We illustrate the potential benefits by means of an example. Since anticausal
multipliers do not yield an improvement of convergence rate guarantees for the 
optimization algorithm analysis problem, we consider
a numerical example unrelated to the problem at hand. 
More precisely, we consider a stable linear system described by the following
transfer function
{
\begin{align}
	G(z) = \frac{z + 0.8111}{z^4 + 1.552 z^3 + 0.6995 z^2 + 0.06042 z - 0.01241}
\end{align}
}
in feedback with a static, slope-restricted uncertainty $ \Delta $ in the sector $ (0,1) $,
i.e., $ \Delta \in \mathbf{\Delta}(0,1) $. 
{Employing~\Cref{lemmaAnalysisConvRate}, we then determine upper bounds on
the exponential convergence rates. The resulting convergence rates are displayed 
in~\Cref{tableComparisonFreeman}. Compared to the multipliers introduced in~\citet{freeman2018noncausal},
the multipliers introduced in the present paper lead to an improvement of
$ 10.85 \perc $.
\begin{table}
	\centering
	\begin{tabular}{cccc}
		\toprule
	                                                           & $ \dimCausal $ & $ \dimAnticausal $ & $ \rho $ \\ \midrule
		Purely causal multipliers                             & $ 1 $          & $ 0 $              & $ 0.977 $ \\
		Anticausal multipliers (\citep{freeman2018noncausal}) & $ 1 $          & $ 3 $              & $ 0.977 $ \\
		Anticausal multipliers (\Cref{lemmaZFIQCrho})         & $ 1 $          & $ 3 $              & $ 0.871 $ \\ \bottomrule
	\end{tabular}
	\caption{Comparison of the resulting convergence rate bounds employing different approaches.}\label{tableComparisonFreeman}
\end{table}}}

{
\subsubsection{State-space realizations}\label{secAppendixStateSpaceRealizations}
In the following we provide explicit state-space realizations of all transfer functions
as required for implementation; an overview is given in~\Cref{tableOverviewSSrealizations}.
\paragraph*{State-space realizations of $ \tfbasisCausal, \tfbasisAnticausal $.}
The state-space realization of $ \tfbasisAnticausal $
as defined in~\eqref{eqDefBasisZF} is given by
\begin{subequations}
\ifthenelse{\boolean{singleColumn}}{
\begin{align}
	\AbasisAnticausal = 
	\left[
	\begin{shortArray}{ccccc}
		0 & 1     & 0 & \cdots \\
		\smash[t]{\vdots} &        \smash[t]{\ddots} & \smash[t]{\ddots} \\
		\smash[t]{\vdots} &               & 0      & 1 \\
		0 & \cdots        & \cdots  & 0
	\end{shortArray}\right], \quad
	\BbasisAnticausal =
	\begin{bmatrix} 0 \\ \smash[t]{\vdots} \\ 0 \\ 1 \end{bmatrix}, \quad
	\CbasisAnticausal &= \left[\begin{shortArray}{cccc} 0 & & \cdots & 1 \\ & & \smash[t]{\iddots} & \\ & 1 & & \\ 1 & 0 & \cdots & 0 \end{shortArray}\right], \quad
	\DbasisAnticausal = \begin{bmatrix} 0 \\ \smash[t]{\vdots} \\ 0 \\ 0 \end{bmatrix}.
\end{align}
}{
\begin{alignat}{3}
	\AbasisAnticausal &= 
	\left[
	\begin{shortArray}{ccccc}
		0 & 1     & 0 & \cdots \\
		\smash[t]{\vdots} &        \smash[t]{\ddots} & \smash[t]{\ddots} \\
		\smash[t]{\vdots} &               & 0      & 1 \\
		0 & \cdots        & \cdots  & 0
	\end{shortArray}\right], \quad
	&
	\BbasisAnticausal &=
	\begin{bmatrix} 0 \\ \smash[t]{\vdots} \\ 0 \\ 1 \end{bmatrix}, \\
	\CbasisAnticausal &= \begin{bmatrix} 0 & \dots & 1 \\ & \iddots & \\ 1 & \dots & 0 \end{bmatrix},
	& 
	\DbasisAnticausal &= \begin{bmatrix} 0 \\ \smash[t]{\vdots} \\ 0 \end{bmatrix}.
\end{alignat}}
\end{subequations}
For $ \tfbasisCausal $ as defined in~\eqref{eqDefBasisZF},
$ \AbasisCausal, \BbasisCausal $ have exactly the same structure 
as $ \AbasisAnticausal, \BbasisAnticausal $ but may be of different size,
while $ \CbasisCausal = I_{\dimCausal} $ and $ \DbasisCausal = 0 $.
\paragraph*{State-space realization of $ \tfMultFac $.}
The state-space realization of $ \tfMultFac $ as defined in~\eqref{eqDefPsiDelta}
is given by
\begin{subequations}
\begin{alignat}{3}
	\AMultFac &= \begin{bmatrix} \AbasisCausal & 0 \\ 0 & \AbasisAnticausal \end{bmatrix} {\otimes I_p},
	\quad &
	\BMultFac &= {\begin{bmatrix} \BbasisCausal \otimes I_p & 0 \\ 0 & \BbasisAnticausal \otimes I_p \end{bmatrix}} {\widehat{\trafoSector}}, \\
	\CMultFac &= 
	\left[\begin{array}{cc}
	0 & 0 \\ 0 & 0 \\ \hline \CbasisCausal & 0 \\ 0 & 0 \\ \hline 0 & 0 \\ 0 & \CbasisAnticausal 
	\end{array}\right] {\otimes I_p},
	\quad &
	\DMultFac &= 
	\left[\begin{array}{cc}
	I_p & 0 \\ 0 & I_p \\ \hline \DbasisCausal & 0 \\ 0 & I_p \\ \hline I_p & 0 \\ 0 & \DbasisAnticausal  
	\end{array}\right].
\end{alignat}
\end{subequations}
\paragraph*{State-space realization of $ \tfcompleteNoPerf $.}
The state-space realization of $ \tfcompleteNoPerf $ as defined in~\eqref{eqSSCompleteSystemNoPerformance}
is given by
\begin{subequations}
\begin{alignat}{3}
	\AcompleteNoPerf(\rho) &=
	\begin{bmatrix}
		A_{\Delta} & \rho^{-1} B_{\Delta} \begin{bmatrix} C \\ 0 \end{bmatrix} \\
		0          & \rho^{-1}\Anom
	\end{bmatrix}, 
	\quad &
	\BcompleteNoPerf &= 
	\begin{bmatrix}
		B_\Delta \begin{bmatrix} 0 \\ I \end{bmatrix} \\
		\Bopt
	\end{bmatrix}, \\
	\CcompleteNoPerf(\rho) &= \begin{bmatrix} C_\Delta & \rho^{-1} D_\Delta \begin{bmatrix} C \\ 0 \end{bmatrix} \end{bmatrix}, 
	\quad &
	\DcompleteNoPerf &= D_\Delta \begin{bmatrix} 0 \\ I \end{bmatrix}
\end{alignat}\label{secAppendixStateSpaceRealizationNoPerf}
\end{subequations}
with $ \Anom = \Aopt + m \Bopt \Copt $.
\paragraph*{State-space realization of $ \tfcomplete $.}
We first note that 
\ifthenelse{\boolean{singleColumn}}{{
\begin{align}
	\left[
	\begin{array}{c|c}
		\Gwtoy & \Gwperftoy \\
		I_p & 0 \\ \hline
		0 & I_{n_{\wperf}} \\
		\Gwtoyperf & \Gwperftoyperf 
	\end{array}
	\right]
	\sim
	\ABCD{A_1}{B_1}{C_1}{D_1}
	=
	\ABCD{\Anom}{\left[\begin{array}{c|c} \Bopt & \Bperf \end{array}\right]}{\left[\begin{array}{c} \Copt \\ 0 \\ \hline 0 \\ \Cperf \end{array}\right]}{\left[\begin{array}{c|c} 0 & 0 \\ I_p & 0 \\ \hline 0 & I_{n_{\wperf}} \\ 0 & 0 \end{array}\right]}
\end{align}}}{
\begin{align}
	\left[
	\begin{array}{c|c}
		\Gwtoy & \Gwperftoy \\
		I_p & 0 \\ \hline
		0 & I_{n_{\wperf}} \\
		\Gwtoyperf & \Gwperftoyperf 
	\end{array}
	\right]
	\sim ( A_1, B_1, C_1, D_1 ),
\end{align}
where
\begin{subequations}
\begin{alignat}{3}
	A_1 &= \Anom, 
	\quad &
	B_1 &= \left[\begin{array}{c|c} \Bopt & \Bperf \end{array}\right], \\
	C_1 &= \left[\begin{array}{c} \Copt \\ 0 \\ \hline 0 \\ \Cperf \end{array}\right],
	\quad & 
	D_1 &= \left[\begin{array}{c|c} 0 & 0 \\ I_p & 0 \\ \hline 0 & I_{n_{\wperf}} \\ 0 & 0 \end{array}\right]
\end{alignat}
\end{subequations}}
with $ \Anom = \Aopt + m \Bopt \Copt $.
Similarly, utilizing that $ \psi_\perf = D_{\psi_\perf} $ is constant, we have
\ifthenelse{\boolean{singleColumn}}{{
\begin{align}
	\left[ 
	\begin{array}{c|c}
		\psi_{\Delta} & 0 \\ \hline 
		0             & \psi_\perf
	\end{array}
	\right]
	\sim 
	\ABCD{A_2}{B_2}{C_2}{D_2}
	= \ABCD{\AMultFac}{\left[\begin{array}{c|c} \BMultFac & 0 \end{array}\right]}{\left[\begin{array}{c} \CMultFac \\ \hline 0 \end{array}\right]}{\left[\begin{array}{c|c} \DMultFac & 0 \\ \hline 0 & D_{\psi_\perf} \end{array}\right]}
\end{align}}}{
\begin{align}
	\left[ 
	\begin{array}{c|c}
		\psi_{\Delta} & 0 \\ \hline 
		0             & \psi_\perf
	\end{array}
	\right]
	\sim 
	( A_2, B_2, C_2, D_2 ),
\end{align}
where 
\begin{subequations}
\begin{alignat}{3}
	A_2 &= \AMultFac,
	\quad &
	B_2 &= \left[\begin{array}{c|c} \BMultFac & 0 \end{array}\right], \\
	C_2 &= \left[\begin{array}{c} \CMultFac \\ \hline 0 \end{array}\right],
	\quad &
	D_2 &= \left[\begin{array}{c|c} \DMultFac & 0 \\ \hline 0 & D_{\psi_\perf} \end{array}\right].
\end{alignat}
\end{subequations}}
Note that $ D_{\psi_\perf} = \begin{bmatrix} 0_{n_{\yperf} \times n_{\wperf}} & I_{n_{\yperf}} \end{bmatrix} $ in the case of $H_2$-performance.
By standard rules for series connections, the state-space 
realization of $ \tfcompleteNoPerf $ as defined in~\eqref{eqDefCompleteTransferFunctionFDI}
is then given by
\ifthenelse{\boolean{singleColumn}}{\ifthenelse{\boolean{saveSpace}}{{
\begin{alignat}{4}
	\Acomplete &= \left[\begin{array}{c|c} A_2 & B_2 C_1 \\ \hline 0 & A_1 \end{array}\right] &&= \AcompleteNoPerf(1) 
	& \quad 
	\Bcomplete &= \left[\begin{array}{c} B_2 D_1 \\ \hline B_1 \end{array}\right] &&= \left[\begin{array}{c|c} \BcompleteNoPerf & \begin{bmatrix} 0 \\ \Bperf \end{bmatrix} \end{array}\right] \nonumber \\
	\Ccomplete &= \begin{bmatrix} C_2 & D_2 C_1 \end{bmatrix} &&= \left[\begin{array}{c} \CcompleteNoPerf(1) \\[0.15em] \hline \\[-0.95em] \begin{bmatrix} 0 & D_{\psi_\perf} \begin{bmatrix} 0 \\ \Cperf \end{bmatrix} \end{bmatrix} \end{array}\right],
	\hspace*{-0.7em} & \Dcomplete &= D_2 D_1 &&= \left[\begin{array}{c|c} \DcompleteNoPerf & 0 \\[0.15em] \hline \\[-0.95em] 0 & D_{\psi_\perf} \begin{bmatrix} I_{n_{\wperf}} \\ 0 \end{bmatrix} \end{array}\right]. \label{eqSSrealizationCompleteSystem}
\end{alignat}}
}{\begin{subequations}
\begin{alignat}{3}
	\Acomplete &= \left[\begin{array}{c|c} A_2 & B_2 C_1 \\ \hline 0 & A_1 \end{array}\right]
	\quad &
	\Bcomplete &= \left[\begin{array}{c} B_2 D_1 \\ \hline B_1 \end{array}\right]  \\
	&= \AcompleteNoPerf(1),
	\quad &
	&= \left[\begin{array}{c|c} \BcompleteNoPerf & \begin{bmatrix} 0 \\ \Bperf \end{bmatrix} \end{array}\right], \nonumber \\
	\Ccomplete &= \begin{bmatrix} C_2 & D_2 C_1 \end{bmatrix}
	\quad &
	\Dcomplete &= D_2 D_1 \\
	&= \left[\begin{array}{c} \CcompleteNoPerf(1) \\ \\[0.15em] \hline \\[-0.95em] \begin{bmatrix} 0 & D_{\psi_\perf} \begin{bmatrix} 0 \\ \Cperf \end{bmatrix} \end{bmatrix} \end{array}\right],
	\quad &
	&= \left[\begin{array}{c|c} \DcompleteNoPerf & 0 \\ \\[0.15em] \hline \\[-0.95em] 0 & D_{\psi_\perf} \begin{bmatrix} I_{n_{\wperf}} \\ 0 \end{bmatrix} \end{array}\right] . \nonumber
\end{alignat}\label{eqSSrealizationCompleteSystem}
\end{subequations}}}{
\begin{subequations}
\begin{alignat}{3}
	\Acomplete &= \left[\begin{array}{c|c} A_2 & B_2 C_1 \\ \hline 0 & A_1 \end{array}\right]
	\quad &
	\Bcomplete &= \left[\begin{array}{c} B_2 D_1 \\ \hline B_1 \end{array}\right]  \\
	&= \AcompleteNoPerf(1),
	\quad &
	&= \left[\begin{array}{c|c} \BcompleteNoPerf & \begin{bmatrix} 0 \\ \Bperf \end{bmatrix} \end{array}\right], \nonumber \\
	\Ccomplete &= \begin{bmatrix} C_2 & D_2 C_1 \end{bmatrix}
	\quad &
	\Dcomplete &= D_2 D_1 \\
	&= \left[\begin{array}{c} \CcompleteNoPerf(1) \\ \\[0.15em] \hline \\[-0.95em] \begin{bmatrix} 0 & D_{\psi_\perf} \begin{bmatrix} 0 \\ \Cperf \end{bmatrix} \end{bmatrix} \end{array}\right],
	\quad &
	&= \left[\begin{array}{c|c} \DcompleteNoPerf & 0 \\ \\[0.15em] \hline \\[-0.95em] 0 & D_{\psi_\perf} \begin{bmatrix} I_{n_{\wperf}} \\ 0 \end{bmatrix} \end{array}\right] . \nonumber
\end{alignat}\label{eqSSrealizationCompleteSystem}
\end{subequations}}
\begin{table*}[t]
\centering
\begin{tabular}{@{}lll@{}}
	\toprule
	Transfer function & Definition & State-space realization  \\ \midrule
	$ \tfbasisAnticausal \in \mathcal{RH}_\infty^{\dimAnticausal \times 1} $ & \eqref{eqDefBasisZF} & $ ( \AbasisAnticausal, \BbasisAnticausal, \CbasisAnticausal, \DbasisAnticausal ) $ \\ 
	$ \tfbasisCausal \in \mathcal{RH}_\infty^{\dimCausal \times 1} $ & \eqref{eqDefBasisZF} & $ ( \AbasisCausal, \BbasisCausal, \CbasisCausal, \DbasisCausal ) $ \\ 
	$ \tfMultFac \in \mathcal{RH}_\infty^{p(\dimCausal+\dimAnticausal+4) \times 2p} $ & \eqref{eqDefPsiDelta} & $ ( \AMultFac, \BMultFac, \CMultFac, \DMultFac ) $  \\
	$ \tfcompleteNoPerf \in \mathcal{RH}_\infty^{p(\dimCausal+\dimAnticausal+4) \times p} $ & \eqref{eqSSCompleteSystemNoPerformance} & $ ( \AcompleteNoPerf, \BcompleteNoPerf, \CcompleteNoPerf, \DcompleteNoPerf ) $  \\
	$ \tfcomplete \in \mathcal{RH}_\infty^{(p(4+\dimCausal+\dimAnticausal)+q) \times 2(p+n_{w\perf})} $ & \eqref{eqDefCompleteTransferFunctionFDI} & $ ( \Acomplete, \Bcomplete, \Ccomplete, \Dcomplete ) $ 
\end{tabular}
\caption{Overview of all required state-space realizations.}
\label{tableOverviewSSrealizations}
\end{table*}

}
\subsubsection{Derivation of~\eqref{eqReformulationConvRateLMIForSynthesis}}\label{secAppendixDerivationH2LMIForSynthesis}
Consider~\eqref{eqLMIConvAnalysis} and note that we may rewrite its
left-hand side as
\begin{align}
	&U(\AcompleteNoPerf(\rho),\BcompleteNoPerf,\CcompleteNoPerf(\rho),\DcompleteNoPerf,P,M_\Delta) \ifthenelse{\boolean{singleColumn}}{}{\nonumber \\
	&\qquad \qquad} - \begin{bmatrix} \star \end{bmatrix}^\top NN^\top P NN^\top \begin{bmatrix} \AcompleteNoPerf(\rho) & \BcompleteNoPerf \end{bmatrix} \prec 0.
\end{align}
Utilizing that $ P_{22} = N^\top P N \succ 0 $ and employing Schur 
complements, \eqref{eqReformulationConvRateLMIForSynthesis} immediately follows.
\subsubsection{Design of {parametrized} optimization algorithms}\label{secAppendixDesignStructuredAlgorithms}
{In this section we consider algorithms of the form~\eqref{eqOptAlgo}
where $ \Aopt, \Bopt $ are parametrized by $ K_i \in \real^{p \times p} $, $ i = 1,2,\dots,n $,
as in~\eqref{eqStructuredAlgoAB} and $ \Copt, \Dopt $ are given by~\eqref{eqDefCoptDdagger}.
Similar to the present paper, we interpret the gradient 
$ \nabla \fObj $ as an uncertainty and embed the design problem into
a robust state-feedback synthesis problem. We follow the procedure from~\citet{mic2014heavy}, 
\citet{mic2016extremum} in a discrete-time setting and extend it to structured objective functions~\eqref{eqExampleStructureGradient}.
To this end, consider an algorithm of the form~\eqref{eqOptAlgo} with 
$ \Aopt, \Bopt $ as in~\eqref{eqStructuredAlgoAB} together with the following state transformation
\begin{align}
	\xTrafo\k{_{,1}} := \nabla \fObj( \Copt \xOpt\k ), \quad \xTrafo\k{_{,i}} = \xOpt\k{_{,i}} \text{ for } i = 2,3,\dots,n.
\end{align}
This yields the transformed dynamics 
\ifthenelse{\boolean{singleColumn}}{
\begin{align}
	\xTrafo\kk
	=~& {\Aopt}_2 \xTrafo\k + \Bopt_1 K \xTrafo\k 
	+ \Bopt_2 \Big( \nabla \fObj\big( C(\Aopt_1+I) \xOpt\k + \Copt \Bopt_1 K \xTrafo\k \big) - \nabla \fObj(\Copt \xOpt\k) \Big),  \label{eqTransformedStructuredAlgo}
\end{align}}{
\begin{align}
	\xTrafo\kk
	=~& {\Aopt}_2 \xTrafo\k + \Bopt_1 K \xTrafo\k \label{eqTransformedStructuredAlgo} \\
	+~& \Bopt_2 \Big( \nabla \fObj\big( C(\Aopt_1+I) \xOpt\k + \Copt \Bopt_1 K \xTrafo\k \big) - \nabla \fObj(\Copt \xOpt\k) \Big), \nonumber
\end{align}}
where
\begin{align}
	{\Aopt}_2 = 
	\left[
	\begin{shortArray}{ccccc}
	I      & 0     & \\
	0      & I     & I & \\
	\smash[t]{\vdots} &       & I & \smash[t]{\ddots} & \\
	\smash[t]{\vdots} &       &   & \smash[t]{\ddots} & I \\
	0      & \dots &   &        & I
	\end{shortArray}
	\right],
	\Bopt_2 = 
	\begin{bmatrix}
		I \\ 0 \\ \smash[t]{\vdots} \\ \smash[t]{\vdots} \\ 0
	\end{bmatrix}.
\end{align}
The main purpose of this transformation is to embed the problem into
a standard robust state-feedback synthesis problem.
In view of~\Cref{subsecParametrizedObjectiveFunctions}, we additionally
suppose that the gradient of the objective function is parametrized as 
\begin{align}
	\nabla \fObj(\optVar) = \fObj_1 \optVar +  \trafoStructured^\top \nabla \fObj_2 (\trafoStructured \optVar), \label{eqExampleStructureGradientAppendix}
\end{align}
where $ \fObj_1 \in \real^{p \times p} $, 
$ \trafoStructured \in \real^{q \times p } $
are known whereas $ \fObj_2 : \real^q \to \real $ is unknown but 
fulfills $ \fObj_2 \in \classObj{m_2}{L_2} $ for some known constants
$ L_2 \geq m_2 > 0 $. Note that this trivially also includes the case when no
structural assumptions are taken on the objective function, simply
by letting $ \fObj_1 = 0 $, $ \trafoStructured = I $. The following
result then follows by using sector bounds for the uncertain terms in~\eqref{eqTransformedStructuredAlgo}. 
\begin{lemma}\label{lemmaSynthesisParametrizedAlgos}
Let $ L_2 \geq m_2 > 0 $, $ n \in \natPos $, $ p \in \natPos $ and $ \fObj_1 \in \real^{p \times p } $,
	$ m_1 I \preceq \fObj_1 \preceq L_1 I $, $ L_1 \geq m_1 \geq 0 $, $\trafoStructured \in \real^{q \times p }, q \in \natPos $, be given.
	Let $ \Copt \in \real^{p\times np} $
	be defined as in~\eqref{eqDefCoptDdagger}.
	Fix $ \rho \in (0,1) $. Suppose there exist $ M \in \real^{p \times np} $, $ Q \in \real^{np \times np} $, such that
	\begin{align}
		\begin{bmatrix}
			Q     & 0     & \bar{\Aopt} Q + \bar{\Bopt} M & \tfrac{1}{2} (L_2 - m_2) {\bar{G}} \\
			\star & I     & \trafoStructured \Copt \Aopt_1 + \trafoStructured \Copt M   & 0 \\
			\star & \star & \rho^2 Q                      & 0 \\
			\star & \star & \star                         & I
		\end{bmatrix}
		\succ 0, \label{eqLMIsynthesisParametrizedAlgos}
	\end{align}
	where
	\begin{subequations}
	\begin{align}
		\bar{\Aopt} &= {\Aopt}_2 + \Bopt_2 \fObj_1 \Copt \Aopt_1 + \tfrac{1}{2}(L_2+m_2) \Bopt_2 \trafoStructured^\top \trafoStructured \Copt \Bopt_1 \\
		\bar{\Bopt} &= \Bopt_1 + \Bopt_2 \fObj_1 \Copt \Bopt_1 + \tfrac{1}{2}(L_2+m_2) \Bopt_2 \trafoStructured^\top \trafoStructured \Copt \Bopt_1 \\
		\bar{G} &= \Bopt_2 \trafoStructured^\top{.}
	\end{align}\label{eqShorthandNotationParametrizedAlgos}
	\end{subequations}
	Then, with $ K = MQ^{-1} $ and $ \Aopt, \Bopt $ given by~\eqref{eqStructuredAlgoAB}
	and $ \Copt, \Dopt $ as in~\eqref{eqDefCoptDdagger}, the equilibrium $ \xOpt^\optSign = \Dopt \optVar^\optSign $
	is globally exponentially stable with rate $ \rho $ for~\eqref{eqOptAlgo} for all
	$ \fObj $ in the form~\eqref{eqExampleStructureGradientAppendix} with $ \fObj_2 \in \classObj{m_2}{L_2} $.
\end{lemma}
\begin{proof}
We prove this result by showing that the equilibrium $ \xTrafo^\optSign = 0 $ 
	is globally exponentially stable with rate $ \rho $ for the transformed dynamics~\eqref{eqTransformedStructuredAlgo}
	for the considered class of objective functions $ \fObj $. To this end,
	we write~\eqref{eqTransformedStructuredAlgo} as a standard robust state-feedback 
	problem and then use quadratic Lyapunov functions together with an $S$-procedure argument. 
	We first calculate for the uncertain terms in~\eqref{eqTransformedStructuredAlgo}
	\begin{align}
		&\nabla \fObj\big( C(\Aopt_1+I) \xOpt\k + \Copt \Bopt_1 K \xTrafo\k \big) - \nabla \fObj(\Copt \xOpt\k) \nonumber \\
		=~& \fObj_1 \big( C(\Aopt_1+I) \xOpt\k + \Copt \Bopt_1 K \xTrafo\k \big) - \fObj_1 \Copt \xOpt\k \nonumber \\
		+~& \trafoStructured^\top \nabla \fObj_2  \big( \trafoStructured C(\Aopt_1+I) \xOpt\k + \trafoStructured \Copt \Bopt_1 K \xTrafo\k \big) - \trafoStructured^\top \nabla \fObj_2 ( \trafoStructured \Copt \xOpt\k ) \nonumber \\
		=~& \fObj_1 ( \Copt \Aopt_1 \xTrafo\k + \Copt \Bopt_1 K \xTrafo\k ) \\
		+~& \trafoStructured^\top \Big( \nabla \fObj_2  \big( \trafoStructured C(\Aopt_1+I) \xOpt\k + \trafoStructured \Copt \Bopt_1 K \xTrafo\k \big) - \nabla \fObj_2 ( \trafoStructured \Copt \xOpt\k ) \Big) \nonumber
	\end{align}
	and then define the uncertainty as
	\ifthenelse{\boolean{singleColumn}}{
	\begin{flalign}
		w\k =~& \nabla \fObj_2  \big( \trafoStructured C(\Aopt_1+I) \xOpt\k + \trafoStructured \Copt \Bopt_1 K \xTrafo\k \big) - \nabla \fObj_2 ( \trafoStructured \Copt \xOpt\k ) 
		- \beta \trafoStructured \Copt ( \Aopt_1 + \Bopt_1 K ) \xTrafo\k \hspace*{-1em} &
	\end{flalign}}{
	\begin{align}
		w\k =~& \nabla \fObj_2  \big( \trafoStructured C(\Aopt_1+I) \xOpt\k + \trafoStructured \Copt \Bopt_1 K \xTrafo\k \big) - \nabla \fObj_2 ( \trafoStructured \Copt \xOpt\k ) \nonumber \\
		&- \beta \trafoStructured \Copt ( \Aopt_1 + \Bopt_1 K ) \xTrafo\k
	\end{align}}
	with $ \beta = \tfrac{1}{2} ( L_2 + m_2 ) $. {The motivation behind this definition
	is that the uncertainty, seen as a function of $ \trafoStructured \Copt ( \Aopt_1 + \Bopt_1 K ) \xTrafo\k $,
	is sector bounded in the sector $ [ -(L_2-\beta), L_2-\beta ] $.
	More precisely, for any $ \xTrafo\k \in \real^{np} $, it fulfills the inequality
	\begin{align}
		\begin{bmatrix} \xTrafo\k \\ w\k \end{bmatrix}^\top 
		\begin{bmatrix} (L_2 - \beta) \Gamma^\top \Gamma & 0 \\ 0 & - \tfrac{1}{L_2 - \beta} I \end{bmatrix}
		\begin{bmatrix} \xTrafo\k \\ w\k \end{bmatrix} \geq 0,
		\label{eqQuadraticIneqUncertainty}
	\end{align}
	where $ \Gamma = \trafoStructured \Copt ( \Aopt_1 + \Bopt_1 K ) $. }
	The transformed dynamics~\eqref{eqTransformedStructuredAlgo} then read as
	\ifthenelse{\boolean{singleColumn}}{
	\begin{align}
		\xTrafo\kk =~& \big( {\Aopt}_2 + \Bopt_2 \fObj_1 \Copt \Aopt_1  + \beta \Bopt_2 \trafoStructured^\top \trafoStructured \Copt \Aopt_1 \big) \xTrafo\k 
		+ ( \Bopt_1 + \Bopt_2 \fObj_1 \Copt \Bopt_1 + \beta \Bopt_2 \trafoStructured^\top \trafoStructured \Copt \Bopt_1 ) K \xTrafo\k \nonumber \\
		&+ \Bopt_2 \trafoStructured^\top w\k. \label{eqTransformedStructuredAlgo2}
	\end{align}}{
	\begin{align}
		\xTrafo\kk =~& \big( {\Aopt}_2 + \Bopt_2 \fObj_1 \Copt \Aopt_1  + \beta \Bopt_2 \trafoStructured^\top \trafoStructured \Copt \Aopt_1 \big) \xTrafo\k \nonumber \\
		&+ ( \Bopt_1 + \Bopt_2 \fObj_1 \Copt \Bopt_1 + \beta \Bopt_2 \trafoStructured^\top \trafoStructured \Copt \Bopt_1 ) K \xTrafo\k \nonumber \\
		&+ \Bopt_2 \trafoStructured^\top w\k. \label{eqTransformedStructuredAlgo2}
	\end{align}}
	Using the shorthand notation~\eqref{eqShorthandNotationParametrizedAlgos},
	\eqref{eqTransformedStructuredAlgo2} can equivalently be written as
	\begin{align}
		\xTrafo\kk = ( \bar{\Aopt} + \bar{\Bopt} K ) \xTrafo\k + \bar{G} w\k.
		\label{eqDynamicsStandardRobustStateFeedback}
	\end{align}
	Hence the problem of finding parameters $K_i$ is a standard
	robust state-feedback problem and the proof is standard from now on.
	Consider a quadratic Lyapunov function candidate $ V : \real^{np} \to \real $
	defined as $ V(\xTrafo) = \xTrafo^\top P \xTrafo $ with $ P \succ 0 $. For 
	exponential convergence with rate $\rho$ of~\eqref{eqDynamicsStandardRobustStateFeedback}
	for all $ w $ that fulfill~\eqref{eqQuadraticIneqUncertainty}
	-- and hence exponential convergence with rate $\rho$ of the original algorithm
	to the minimizer of $\fObj$ for all $ \fObj \in \classObj{m}{L} $ --
	it is sufficient that there exists a $ P \succ 0 $ such that
	\ifthenelse{\boolean{singleColumn}}{
	\begin{align}
		V(\xTrafo\kk) - \rho^2 V(\xTrafo\k) 
		= 
		\begin{bmatrix} \vphantom{\begin{bmatrix} \xTrafo\k \\ w\k \end{bmatrix}} \star \end{bmatrix}^\top
		\begin{bmatrix} \vphantom{
		\begin{bmatrix}
			\bar{\Aopt} + \bar{\Bopt} K & G \\
			I                           & 0
		\end{bmatrix}} \star
		\end{bmatrix}^\top
		\begin{bmatrix} P & 0 \\ 0 & -\rho^2 P \end{bmatrix} 
		\begin{bmatrix} 
			\bar{\Aopt} + \bar{\Bopt} K & G \\
			I                           & 0
		\end{bmatrix}
		\begin{bmatrix} \xTrafo\k \\ w\k \end{bmatrix}
		< 0  
	\end{align}}{
	\begin{align}
		&V(\xTrafo\kk) - \rho^2 V(\xTrafo\k) 
		= \nonumber \\
		&\begin{bmatrix} \vphantom{\begin{bmatrix} \xTrafo\k \\ w\k \end{bmatrix}} \star \end{bmatrix}^\top
		\begin{bmatrix} \vphantom{
		\begin{bmatrix}
			\bar{\Aopt} + \bar{\Bopt} K & {\bar{G}} \\
			I                           & 0
		\end{bmatrix}} \star
		\end{bmatrix}^\top
		\begin{bmatrix} P & 0 \\ 0 & -\rho^2 P \end{bmatrix} 
		\begin{bmatrix} 
			\bar{\Aopt} + \bar{\Bopt} K & G \\
			I                           & 0
		\end{bmatrix}
		\begin{bmatrix} \xTrafo\k \\ w\k \end{bmatrix}
		< 0  
	\end{align}}
	for all {$ \begin{bmatrix} \xTrafo\k^\top & w\k^\top \end{bmatrix}^\top \neq 0 $} that fulfill~\eqref{eqQuadraticIneqUncertainty}.
	By the $S$-procedure it is hence sufficient that there exist $ P \succ 0 $
	and $ \lambda \geq 0 $ such that
	\ifthenelse{\boolean{singleColumn}}{
	\begin{align}
		\left[
		\vphantom{
		\begin{bmatrix}
			\bar{\Aopt} + \bar{\Bopt} K & {\bar{G}} \\
			I                           & 0
		\end{bmatrix}}
		\star
		\right]^\top
		\begin{bmatrix} P & 0 \\ 0 & -\rho^2 P \end{bmatrix} 
		\begin{bmatrix} 
			\bar{\Aopt} + \bar{\Bopt} K & {\bar{G}} \\
			I                           & 0
		\end{bmatrix} 
		+
		\lambda 
		\begin{bmatrix} (L_2 - \beta) \Gamma^\top \Gamma & 0 \\ 0 & - \tfrac{1}{L_2 - \beta} I \end{bmatrix}
		\prec 0.
	\end{align}
	}{
	\begin{align}
		\begin{bmatrix}
			\bar{\Aopt} + \bar{\Bopt} K & {\bar{G}} \\
			I                           & 0
		\end{bmatrix}^\top
		\begin{bmatrix} P & 0 \\ 0 & -\rho^2 P \end{bmatrix} 
		\begin{bmatrix} 
			\bar{\Aopt} + \bar{\Bopt} K & {\bar{G}} \\
			I                           & 0
		\end{bmatrix} \nonumber \\
		+
		\lambda 
		\begin{bmatrix} (L_2 - \beta) \Gamma^\top \Gamma & 0 \\ 0 & - \tfrac{1}{L_2 - \beta} I \end{bmatrix}
		\prec 0.
	\end{align}}
	We exclude the case $ \lambda = 0 $ and thus set $ \lambda = 1 $
	without loss of generality. We can write the above inequality
	equivalently as
	\begin{align}
		\begin{bmatrix}
		\vphantom{
		\begin{bmatrix}
			P ( \bar{\Aopt} + \bar{\Bopt} K ) & P{\bar{G}} \\
			\sqrt{L_2-\beta} \Gamma           & 0
		\end{bmatrix}} \star 
		\end{bmatrix}^\top
		\left[\begin{shortArray}{cc} P^{-1} & 0 \\ 0 & I \end{shortArray}\right]
		\left[ 
		\begin{shortArray}{cc}
			P ( \bar{\Aopt} + \bar{\Bopt} K ) & P{\bar{G}} \\
			\sqrt{L_2-\beta} \Gamma           & 0
		\end{shortArray}\right]
		-
		\left[\begin{shortArray}{cc} \rho^2 P & 0 \\ 0 & \tfrac{1}{L_2 - \beta} I \end{shortArray}\right]
		\prec 0.
	\end{align}
	Since $P$ is assumed to be a positive definite matrix, 
	we can employ Schur complements such that the above matrix inequality
	holds if and only if
	\begin{align}
		\begin{bmatrix}
			P     & 0     & P ( \bar{\Aopt} + \bar{\Bopt} K ) & P{\bar{G}} \\
			\star & I     & \sqrt{L_2-\beta} \Gamma           & 0 \\
			\star & \star & \rho^2 P                          & 0 \\
			\star & \star & \star                             & \tfrac{1}{L_2 - \beta} I 
		\end{bmatrix}
		\succ 0.
	\end{align}
	We multiply from left and right by the positive definite block diagonal
	matrix $ \textup{blkdiag}(\sqrt{L_2-\beta} P^{-1},I, \sqrt{L_2-\beta} P^{-1},\sqrt{L_2-\beta} I ) $
	such that we obtain
	\begin{align}
		\left[
		\begin{shortArray}{cccc}
			(L_2-\beta) P^{-1} & 0     & (L_2-\beta)( \bar{\Aopt} + \bar{\Bopt} K ) P^{-1} & (L_2-\beta) {\bar{G}} \\
			\star              & I     & (L_2-\beta) \Gamma P^{-1}                         & 0 \\
			\star              & \star & \rho^2 (L_2-\beta) P^{-1}                         & 0 \\
			\star              & \star & \star                                             & I
		\end{shortArray}\right]
		\succ 0.
	\end{align}
	We note that the additional scaling is not required
	and included for numerical purpose only.
	The variable transformation $ Q = (L_2-\beta) P^{-1} $, $ M = (L_2-\beta) K P^{-1} $
	then gives~\eqref{eqLMIsynthesisParametrizedAlgos}, thus concluding the proof.
\end{proof}
}


\end{document}